\documentclass{amsart}
\pdfoutput 1
\usepackage[all,hyperref,numberbysection]{bi-discrete}
\usepackage{amsmath}
\usepackage{amssymb}
\usepackage{upgreek}
\usepackage{microtype}
\numberwithin{equation}{section}

\usepackage[backend=biber,maxbibnames=5,maxalphanames=5,style=alphabetic,bibencoding=utf8,giveninits,url=false,isbn=false]{biblatex}
\AtBeginBibliography{\small}
\allowdisplaybreaks

\newtheorem*{ack}{Acknowledgments}
\DeclareMathOperator*{\osc}{osc}

\begin{document}

\author[Diening]{Lars Diening}
\address{Fakult\"at f\"ur Mathematik, Universit\"at Bielefeld, Postfach 100131, D-33501 Bielefeld, Germany}
\email{lars.diening@uni-bielefeld.de}
\author[Kim]{Kyeongbae Kim}
\address{Department of Mathematical Sciences, Seoul National University, Seoul 08826, Korea}
\email{kkba6611@snu.ac.kr}
\author[Lee]{Ho-Sik Lee}
\address{Fakult\"at f\"ur Mathematik, Universit\"at Bielefeld, Postfach 100131, D-33501 Bielefeld, Germany}
\email{ho-sik.lee@uni-bielefeld.de}
\author[Nowak]{Simon Nowak}
\address{Fakult\"at f\"ur Mathematik, Universit\"at Bielefeld, Postfach 100131, D-33501 Bielefeld, Germany}
\email{simon.nowak@uni-bielefeld.de}

\makeatletter
\@namedef{subjclassname@2020}{\textup{2020} Mathematics Subject Classification}
\makeatother


\thanks{Lars Diening and Simon Nowak gratefully acknowledge funding by the Deutsche Forschungsgemeinschaft (DFG, German Research Foundation) - SFB 1283/2 2021 - 317210226. Kyeongbae Kim thanks for funding of the National Research Foundation of Korea (NRF) through IRTG 2235/NRF-2016K2A9A2A13003815 at Seoul National University. Ho-Sik Lee thanks for funding of the Deutsche Forschungsgemeinschaft through GRK 2235/2 2021 - 282638148.
}

\title{Nonlinear nonlocal potential theory at the gradient level}

	\begin{abstract}
		The aim of this work is to establish numerous interrelated gradient estimates in the nonlinear nonlocal setting. First of all, we prove that weak solutions to a class of homogeneous nonlinear nonlocal equations of possibly arbitrarily low order have H\"{o}lder continuous gradients.
		Using these estimates in the homogeneous case, we then prove sharp higher differentiability as well as pointwise gradient potential estimates for nonlinear nonlocal equations of order larger than one in the presence of general measure data. Our pointwise estimates imply that the first-order regularity properties of such nonlinear nonlocal equations coincide with the sharp ones of the fractional Laplacian.
	\end{abstract}

	\maketitle
	
	\tableofcontents
	

	\section{Introduction}
	\subsection{Aim and scope}
	The primary goal of this paper is to establish for the first time pointwise gradient estimates in terms of Riesz potentials of the data for solutions to equations that are both nonlinear and nonlocal in nature.
	
	To be more precise, we study nonlinear nonlocal equations of the type
	\begin{equation}
		\label{pt : eq.main}
		\setLcal u=\mu\quad\text{in }\Omega \subset \mathbb{R}^n,
	\end{equation}
	where $n \geq 2$ and the nonlinear nonlocal operator $\setLcal$ is formally defined by
	\begin{equation} \label{nonlocalop}
		\setLcal u(x)=(1-s) \hspace{1mm} \mathrm{P}.\mathrm{V}.\int_{\RRn}{{\Phi}}\left(\frac{u(x)-u(y)}{|x-y|^{s}}\right)\frac{\,dy}{{|x-y|^{n+s}}}.
	\end{equation}
	Here $s\in(0,1)$ is a fixed parameter, $\mu$ belongs to the class $\mathcal{M}(\mathbb{R}^n)$ of signed Radon measures on $\mathbb{R}^n$ with finite total mass, while the nonlinearity $\Phi$ satisfies the following Lipschitz and monotonicity assumptions:

	\begin{assumption} \label{assump}
		We assume that $\Phi:\bbR\to\bbR$ is an odd function such that for all $t,t^\prime \in \mathbb{R}$ and some $\Lambda \geq 1$, we have
		\begin{equation}
			\label{pt : assmp.phi}
			|{{\Phi}}(t)-\Phi(t^\prime)|\leq \Lambda|t-t^\prime|\quad\text{and}\quad ({\Phi}(t)-{\Phi}(t^\prime))(t-t^\prime)\geq\Lambda^{-1}|t-t^\prime|^{2}.
		\end{equation}
	\end{assumption}
	
	We note that \eqref{pt : eq.main} is the Euler-Lagrange equation of the nonlocal energy functional
	\begin{equation} \label{energy}
	u \mapsto (1-s) \iint_{(\Omega^c \times \Omega^c)^c}{{\Psi}}\left(\frac{|u(x)-u(y)|}{|x-y|^{s}}\right)\frac{\,dx\,dy}{{|x-y|^{n}}} - \int_{\Omega} u \, d \mu
	\end{equation}
	in case $\mu$ is sufficiently regular, where the Young function $\Psi:[0,\infty) \to [0,\infty)$ is given by $\Psi(t):=\int_0^t \Phi(\tau)d\tau$.
	
	Studying the regularity of solutions to nonlinear nonlocal equations has become a very active research area in recent years, see e.g.\ \cite{KassCalcVar,CSCPAM,CSannals,CCV,KMS1,SchikorraMA,RSDuke,DKP,BL,BLS,CCDuke,DFPJDE,CKWCalcVar,BKJMA,DNCZ,GL23} for a non-exhaustive list of fundamental contributions in this direction.
	
	This rapid development of nonlocal regularity theory was in particular facilitated by the large number of applications of nonlocal models in many areas of pure and applied mathematics such as for instance stochastic processes of jump-type (see e.g.\ \cite{bertoin,Fukushima}), classical harmonic analysis (see e.g.\ \cite{landkof}), conformal geometry (see e.g.\ \cite{GZ,CG,Case-Chang}), phase transitions (see e.g.\ \cite{fife}), physics of materials and relativistic models (see e.g.\ \cite{LiebYau2}), fluid dynamics (see e.g.\ \cite{NInvent,CaffVass}) and kinetic theory (see e.g.\ \cite{ImbSil}). Moreover, nonlocal energy functionals of the particular type \eqref{energy} arise in image processing (see e.g.\ \cite{GOsh}). 
	In addition, for embedding properties of energies of the type \eqref{energy} we refer to \cite{CianchiFrac}.

	\subsubsection{Potential estimates for local elliptic equations}
	Before stating our results in a precise fashion, let us provide some context and motivation. To do this, we begin by considering on the whole space the classical Poisson equation given by
	\begin{equation}
		\label{pt : PE}
		-\Delta u=\mu\quad\text{in } \mathbb{R}^n.
	\end{equation}
	It is well known that if $u$ decays to zero at infinity, then at least formally we have the representation formula
	\begin{equation} \label{Greenrep}
		u(x)= \int_{\mathbb{R}^n} G(x,y)\,d\mu(y),
	\end{equation}
	where for $n \geq 3$ the Green function/fundamental solution $G$ of the Laplacian satisfies
	$$
	G(x,y) \approx
	|x-y|^{2-n}.
	$$
	This concrete representation formula directly yields the following zero-order pointwise estimate in terms of Riesz potentials \smallskip 
	\begin{equation} \label{PotEst1}
	|u(x)| \lesssim \int_{\bbR^n} \frac{d|\mu|(y)}{|x-y|^{n-2}} =: I_2^{|\mu|}(x),
	\end{equation}
	as well as the gradient potential estimate
	\begin{equation} \label{PotEst2}
	|\nabla u(x)| \lesssim \int_{\bbR^n} \frac{d|\mu|(y)}{|x-y|^{n-1}} =: I_1^{|\mu|}(x).
	\end{equation}
	Since the explicit representation of solutions given by \eqref{Greenrep} is based on convolution and therefore relies on the linearity of the Poisson equation \eqref{pt : PE}, a highly nontrivial question that subsequently arose is if similar pointwise potential estimates remain valid for nonlinear generalizations of \eqref{pt : PE}. In fact, in their seminal paper \cite{KM}, Kilpeläinen and Mal\'y succeeded in recovering pointwise zero-order potential estimates of the type \eqref{PotEst1} for nonlinear second-order elliptic equations of the type
	\begin{equation}
		\label{pt : LNE}
		-\textnormal{div}(\textbf{a}(\nabla u))=\mu
	\end{equation}
	under very general growth assumptions on the vector field $\textbf{a}$. An alternative proof of this result was later given by Trudinger and Wang in \cite{TWAJM}.
	In another pioneering paper \cite{Min}, Mingione then proved that also the gradient potential estimate \eqref{PotEst2} remains true for local nonlinear equations of the type \eqref{pt : LNE} under assumptions on $\textbf{a}$ that are similar to our assumptions on $\Phi$ given by \eqref{pt : assmp.phi}. 
	
	Since the mapping properties of Riesz potentials are classically known in a large variety of function spaces measuring even very refined scales, the potential estimates obtained in \cite{KM} and \cite{Min} essentially imply that the zero-order and first-order regularity properties of solutions to nonlinear equations of the type \eqref{pt : LNE} coincide with the known sharp ones of the Laplacian, linearizing the theory up to the gradient level. In a large number of subsequent papers, similar potential estimates were then established also for even more general second-order elliptic and parabolic equations and systems, see for instance \cite{DM1,DM2,DKSBMO,KuMiARMA1,KuMiG,Baroni,KuMiV,BCDKS,CM1,BYP,BCDS,NNARMA,DongJEMS,DFJMPA,CKW23}.
	
	\subsubsection{Potential estimates for nonlocal equations}
	Let us now turn to discussing similar results for nonlocal equations. Let us once more begin by discussing the linear case, this time given by the fractional Laplacian, which for $s \in (0,1)$ can be defined by
	$$(-\Delta)^s u(x) := c_{n,s} \hspace{1mm} \mathrm{P}.\mathrm{V}. \int_{\mathbb{R}^n} \frac{u(x)-u(y)}{|x-y|^{n+2s}}\,dy,$$
	where $c_{n,s}$ is a positive constant that guarantees that $(-\Delta)^s$ converges to the local Laplacian $-\Delta$ as $s \to 1$ in a suitable sense. In analogy to the local setting, for $n \geq 2$ the Green function or fundamental solution $G_s$ of $(-\Delta)^s$ in $\mathbb{R}^n$ satisfies
	$$G_s(x,y) \approx |x-y|^{2s-n},$$
	see for instance \cite[Theorem 8.4]{Garofalo}. Therefore, in a similar way as in the local setting for solutions to the fractional Poisson equation
	\begin{equation}
		\label{pt : FPE}
		(-\Delta)^s u=\mu\quad\text{in } \mathbb{R}^n,
	\end{equation}
	at least formally we have the explicit representation formula
	\begin{equation} \label{Greenrepfrac}
		u(x)= \int_{\mathbb{R}^n} G_s(x,y)\,d\mu(y),
	\end{equation}
	provided $u$ decays to zero at infinity. As before, formula \eqref{Greenrepfrac} directly yields the pointwise zero-order Riesz potential estimate
	\begin{equation} \label{PotEst1f}
		|u(x)| \lesssim \int_{\mathbb{R}^n} \frac{d|\mu|(y)}{|x-y|^{n-2s}} =: I_{2s}^{|\mu|}(x),
	\end{equation}
	and provided that $s \in (1/2,1)$, also the gradient potential estimate
	\begin{equation} \label{PotEst2f}
		|\nabla u(x)| \lesssim \int_{\mathbb{R}^n} \frac{d|\mu|(y)}{|x-y|^{n-2s+1}} =: I_{2s-1}^{|\mu|}(x).
	\end{equation}
	
	In light of the pointwise estimates \eqref{PotEst1f} and \eqref{PotEst2f} and the well-established nonlinear potential theory in the local setting, a natural question arising at this point is if similar potential estimates can be proved also for more general and in particular nonlinear nonlocal equations of the type \eqref{pt : eq.main}-\eqref{nonlocalop}. In fact, in \cite{KuuMinSir15} Kuusi, Mingione and Sire showed that zero-order estimates similar to \eqref{PotEst1f} remain valid for a large class of nonlinear nonlocal equations that in particular contains the type of equations we treat in the present paper. In a number of subsequent papers, similar zero-order potential estimates were then established for even more general nonlinear nonlocal equations, see e.g.\ \cite{KLL,Minhyun1,NNSW}.

	Concerning gradient potential estimates similar to \eqref{PotEst2f} for nonlinear nonlocal equations, to the best of our knowledge no previous results prior to the present work seem to exist. Nevertheless, in \cite{KuuSimYan22} together with Kuusi and Sire the last-named author established gradient potential estimates for solutions to a class of linear nonlocal equations of order larger than one that arise when the fractional Laplacian is perturbed by H\"older continuous coefficients, which already turned out to be delicate despite the linearity of the equation. This is mainly since in contrast to local second-order equations, the gradient is not naturally associated with nonlocal equations of the type we have in mind, leading to substantial additional technical difficulties at essentially any stage of the proof already in the linear case. 
	
	However, this still left open the central question if it is possible to obtain gradient potential estimates for nonlinear nonlocal equations. In fact, in the present work we finally answer this question affirmatively for equations of the type \eqref{pt : eq.main} under the natural Lipschitz and monotonicity assumptions on the nonlinearity $\Phi$ given by \eqref{pt : assmp.phi} that are in similar to the assumptions made in the local case in \cite{Min}. More precisely, we provide estimates similar to \eqref{PotEst2f} for nonlinear nonlocal equations posed on the whole space as well as for equations in bounded domains. On the whole space, we indeed recover the estimate \eqref{PotEst2f} in its exact form, as our first main result shows.
	\begin{theorem}[Gradient potential estimates in $\mathbb{R}^n$]
		\label{el : thm.main3}
		Let $s \in (1/2,1)$, $\mu \in \mathcal{M}(\mathbb{R}^n)$ and let $u \in W^{s,2}(\mathbb{R}^n)$ be a weak solution of \eqref{pt : eq.main} in $\mathbb{R}^n$. Moreover, suppose that $\Phi$ satisfies Assumption \ref{assump} for some $\Lambda \geq 1$. Then for almost every $x_0 \in \mathbb{R}^n$, we have the pointwise estimate
		\begin{equation} \label{PERn}
			|\nabla u(x_{0})| \leq cI_{2s-1}^{|\mu|}(x_{0})
		\end{equation}
		for some constant $c=c(n,s,\Lambda)$. In addition, for any fixed $s_0\in(1/2,1)$, the constant $c$ depends only on $n,s_{0}$, and $\Lambda$ whenever $s\in [s_0,1)$.
	\end{theorem}
	
	For the case of bounded domains, we refer to Theorem \ref{el : thm.main2}. In this case we prove our gradient potential estimates for a very general class of solutions called SOLA (=\hspace{0.2mm} solutions obtained by limiting approximations), which are known to always exist even in the presence of general measure data as proved in \cite{KuuMinSir15}.

	Moreover, since we provide estimates that are stable as $s \to 1$ and at least heuristically nonlocal equations converge to local ones of second-order as the order of the equation approaches two (see e.g.\ \cite{BBM1,FKV20} for some rigorous results in this direction), our gradient potential estimates can indeed be considered to be nonlocal analogues of the ones obtained in \cite{Min}.
	
	Concerning applications, the pointwise gradient bounds we prove in particular lead to sharp gradient regularity results in a large variety of function spaces measuring fine scales. In other words, we obtain that essentially the entire first-order regularity theory of solutions to \eqref{pt : eq.main} coincides with the sharp one of the fractional Poisson equation \eqref{pt : FPE}, providing an in principle complete linearization of the regularity theory of solutions to the equation \eqref{pt : eq.main} at the gradient level. 
	
	In addition, along the way we shall also obtain further interrelated new estimates for such equations that are also interesting in their own right, in the sense that they provide control of the oscillations rather than the size of the gradient of solutions in H\"older and fractional Sobolev spaces.

	\subsection{Setting and further main results}
	Before stating our other main results, we have to introduce our setup in a more rigorous fashion.
	
	In order to control the growth of solutions at infinity, we consider the tail space
	$$L^1_{2s}(\mathbb{R}^n):= \left \{u \in L^1_{\loc}(\mathbb{R}^n) \mathrel{\Big|} \int_{\mathbb{R}^n} \frac{|u(y)|}{(1+|y|)^{n+2s}}\,dy < \infty \right \}$$
	
	introduced in \cite{existence}. We note that a function $u \in L^1_{\loc}(\mathbb{R}^n)$ belongs to the space $L^{1}_{2s}(\mathbb{R}^n)$ if and only if the \emph{nonlocal tails} of $u$ given by
	$$ \textnormal{Tail}(u;B_R(x_0)):= (1-s) R^{2s} \int_{\mathbb{R}^n \setminus B_{R}(x_0)} \frac{|u(y)|}{|x_0-y|^{n+2s}}\,dy$$
	are finite for all $R>0$ and $x_0 \in \mathbb{R}^n$.
	
	For notational convenience, we also consider the standard local excess functional as well as the following nonlocal excess functional which was introduced in \cite{KuuMinSir15}.
	\begin{definition}[Local and nonlocal excess functional] \label{def:excess}
		Fix $x_0 \in \mathbb{R}^n$ and $R>0$. For any $u \in L^1(B_R(x_0))$ and any $q \in [1,\infty)$, we define the local excess functional $E_{\mathrm{loc}}^q(u;B_R(x_0))$ of $u$ by
		$$ E_{\mathrm{loc}}^q(u;B_R(x_0)):=\left ( \dashint_{B_{R}(x_0)} |u-(u)_{B_R(x_0)}|^q \,dx \right )^\frac{1}{q}.$$
		We also set $E_{\mathrm{loc}}(u;B_R(x_0)):=E_{\mathrm{loc}}^1(u;B_R(x_0))$. Moreover, given $s \in (0,1)$ and $u \in L^{1}_{2s}(\mathbb{R}^n)$, we define the nonlocal excess functional $E(u;B_R(x_0))$ of $u$ by
		$$
		E(u;B_R(x_0)):= E_{\mathrm{loc}}(u;B_R(x_0)) + \textnormal{Tail}(u-(u)_{B_R(x_0)};B_R(x_0)).
		$$
	\end{definition}
	
	We now define standard energy-type weak solutions to \eqref{pt : eq.main} as follows.
	\begin{definition}[Weak solutions]
		Let $\Omega \subset \mathbb{R}^n$ be a domain. Given $\mu\in \mathcal{M}(\mathbb{R}^n)$, we say that $u\in W^{s,2}_{\mathrm{loc}}(\Omega)\cap L^{1}_{2s} (\bbR^{n})$
		is a weak solution to \eqref{pt : eq.main} if
	\begin{equation} \label{distsol}
		\begin{aligned}
		(1-s) & \int_{\bbR^{n}}\int_{\bbR^{n}}{\Phi}\left(\frac{u(x)-u(y)}{|x-y|^{s}}\right)\frac{\psi(x)-\psi(y)}{|x-y|^{s}}\frac{\,dx\,dy}{|x-y|^{n}} \\& =\int_{\Omega} \psi\,d\mu \quad \forall \psi \in C_c^\infty(\Omega).
		\end{aligned}
	\end{equation}
		\end{definition}
	For the definition of the standard local fractional Sobolev spaces $W^{s,2}_{\loc}(\Omega)$, we refer to Section \ref{sec:2}.
	We also note that if the measure $\mu$ is given by a function belonging to $L^{\frac{2n}{n+2s}}_{\mathrm{loc}}(\Omega)$ and thus belongs to the dual of $W^{s,2}$, then in view of the assumptions \eqref{pt : assmp.phi} it is possible to prove the existence of weak solutions belonging to the energy space $W^{s,2}$ by standard monotonicity methods, see e.g.\ \cite[Remark 3]{existence}, \cite[Appendix A]{Kyeongbae2} or \cite[Proposition 4.1]{MeN}.

	We are now in the position to state our first main result, which is concerned with H\"{o}lder estimates for the gradient of weak solutions to nonlinear nonlocal equations of the type \eqref{pt : eq.main} in the homogeneous case when $\mu \equiv 0$ that are stable as $s \to 1$.
	\begin{theorem}[$C^{1,\alpha}$ regularity]
		\label{el : thm.main1}
		Let $\Omega \subset \mathbb{R}^n$ be a domain and $s \in (0,1)$. Moreover, suppose that $\Phi$ satisfies Assumption \ref{assump} for some $\Lambda \geq 1$ and assume that $u\in W^{s,2}_{\mathrm{loc}}(\Omega)\cap L^{1}_{2s} (\bbR^{n})$ is a weak solution of $$\setLcal u =0 \quad \textnormal{in } \Omega.$$
		Then there exists a constant $\alpha=\alpha(n,s,\Lambda)\in(0,1)$ such that $u\in C^{1,\alpha}_{\mathrm{loc}}(\Omega)$. Moreover, for any $x_0 \in \Omega$ and any $R>0$ with ${B}_{R}(x_{0})\Subset \Omega$, we have the estimate
		\begin{equation} \label{Holdest}
			\|\nabla u\|_{L^{\infty}(B_{R/2}(x_{0}))}+R^{\alpha} [\nabla u]_{C^\alpha(B_{R/2}(x_0))}
			\leq{c}E(u/{R};B_{R}(x_{0}))
		\end{equation}
		for some constant $c=c(n,s,\Lambda)$. In addition, for any fixed $s_0\in(0,1)$, the constants $c$ and $\alpha$ depend only on $n,s_{0}$, and $\Lambda$ whenever $s\in [s_0,1)$.
	\end{theorem}

	\begin{remark} \label{Alls} \normalfont
		An interesting feature of Theorem \ref{el : thm.main1} is that it yields $C^{1,\alpha}$ regularity in the whole range $s \in (0,1)$ and thus for nonlocal equations of possibly arbitrarily low order. In other words, the order of the equation does not obstruct gradient regularity even in the nonlinear regime as long as the right-hand side is regular.
	\end{remark}
	
	Since standard weak solutions to \eqref{pt : eq.main} might in general not exist under general measure data, in this case we consider the following notion of SOLA solutions introduced in \cite{KuuMinSir15}.
	
	\begin{definition} \label{SOLA}
		Consider a bounded domain $\Omega \subset \mathbb{R}^n$, let $\mu \in \mathcal{M}(\mathbb{R}^n)$ and $g \in W^{s,2}_{\loc}(\mathbb{R}^n) \cap L^1_{2s}(\mathbb{R}^n)$. We say that a function $u$ is a SOLA of the Dirichlet problem 
		\begin{equation} \label{NonlocalDir}
			\begin{cases} \normalfont
				\setLcal u = \mu & \text{ in } \Omega \\
				u=g & \text{ a.e. in } \mathbb{R}^n \setminus \Omega,
			\end{cases}
		\end{equation}
		if the following conditions are satisfied:
		\begin{itemize}
			\item[(1)] $u \in W^{h,q}(\Omega)$ for any $h \in (0,s)$ and any $q \in \big [1,\frac{n}{n-s} \big )$.
			\item[(2)] $u$ is a distributional solution of $\setLcal u=\mu$ in $\Omega$ in the sense that \eqref{distsol} holds.
			\item[(3)] $u=g$ a.e.\ in $\mathbb{R}^n \setminus \Omega$.
			\item[(4)] There exist sequences of functions $\{u_j\}_{j=1}^\infty \subset W^{s,2}(\mathbb{R}^n)$, $\{\mu_j\}_{j=1}^\infty \subset C_0^\infty(\mathbb{R}^n)$, $\{g_j\}_{j=1}^\infty \subset C_0^\infty(\mathbb{R}^n)$ such that each $u_j$ weakly solves the Dirichlet problem
			$$
			\begin{cases} \normalfont
				\setLcal u_j = \mu_j & \text{ in } \Omega \\
				u_j=g_j & \text{ a.e. in } \mathbb{R}^n \setminus \Omega.
			\end{cases}
			$$
			\item[(5)] $u_j$ converges to $u$ a.e.\ in $\mathbb{R}^n$ and locally in $L^q(\mathbb{R}^n)$.
			\item[(6)] The sequence $\{\mu_j\}_{j=1}^\infty$ converges weakly to $\mu$ in the sense of measures in $\Omega$ and additionally satisfies
			\begin{equation*}
				\limsup_{j \to \infty} |\mu_j|(B) \leq |\mu|(\overline B) \text{ for any ball } B.
			\end{equation*}
			\item[(7)] The sequence $\{g_j\}_{j=1}^\infty$ converges to $g$ in the following sense: For any $z \in \mathbb{R}^n$ and any $r>0$, we have
			\begin{equation*} 
				g_j \to g \text{ in } W^{s,2}(B_r(z)), \quad \lim_{j \to \infty} \int_{\setR^n\setminus B_r(z)} \frac{|g_j(y)-g(y)|}{|z-y|^{n+2s}} \,dy=0.
			\end{equation*}
		\end{itemize}
	\end{definition}

	\begin{remark}[Existence] \label{SOLAExistence} \normalfont
	We note that for any $\mu \in \mathcal{M}(\mathbb{R}^n)$ and any $g \in W^{s,2}_{\loc}(\mathbb{R}^n) \cap L^1_{2s}(\mathbb{R}^n)$, there in fact always exists a SOLA of (\ref{NonlocalDir}) whenever $\Phi$ satisfies Assumption \ref{assump}, which can be proved in exactly the same way as in \cite[Theorem 1.1]{KuuMinSir15}.
	\end{remark}

	A crucial additional difficulty when trying to prove gradient potential estimates in the nonlocal setting is that the gradient of a SOLA or even weak solution is a priori not even known to be integrable, but merely exists as a distribution. This is in sharp contrast to local elliptic equations of the type \eqref{pt : LNE}, for which the gradient of any SOLA as defined e.g.\ in \cite{DM2} is locally integrable by definition. In our nonlocal setting, establishing gradient integrability is instead already a highly nontrivial task and follows from the following higher differentiability result.

	\begin{theorem}[Higher differentiability under measure data] \label{HD}
		Let $s \in (1/2,1)$, $\mu \in \mathcal{M}(\mathbb{R}^n)$, $g \in W^{s,2}_{\loc}(\mathbb{R}^n) \cap L^1_{2s}(\mathbb{R}^n)$ and let $u$ be a SOLA of (\ref{NonlocalDir}) in a bounded domain $\Omega \subset \mathbb{R}^n$. Moreover, suppose that $\Phi$ satisfies Assumption \ref{assump} for some $\Lambda \geq 1$. Then for all 
		\begin{align}
  		\label{el : choi.sigq}
			\sigma\in(0,2s-1)\quad\text{and}\quad q\in\left[1,\frac{n}{n-2s+1+\sigma}\right),
		\end{align}
		we have $ u\in W^{1+\sigma,q}_{\loc}(\Omega)$. Moreover, for any $x_0 \in \Omega$ and any $R>0$ with $B_R(x_0) \Subset \Omega$,
		\begin{align*}
		& \left (\dashint_{B_{R/2}(x_0)} |\nabla u|^q \,dx \right )^{1/q} + R^{\sigma }\left (\dashint_{B_{R/2}(x_0)} \int_{B_{R/2}(x_0)} \frac{|\nabla u(x)-\nabla u(y)|^q}{|x-y|^{n+\sigma q}} \,dx\,dy \right )^{1/q} \\ &\leq {c}E(u/{R};B_R(x_0))+c\frac{|\mu|(B_{R}(x_{0}))}{R^{n-2s+1}}
		\end{align*}
		for some constant $c=c(n,s,\Lambda,\sigma,q)$. In addition, for any fixed $s_0\in(1/2,1)$ if $\sigma $ and $q$ are in \eqref{el : choi.sigq} with $s=s_{0}$, then the constant $c$ depends only on $n,s_{0},\Lambda,\sigma$ and $q$ whenever $s \in [s_0,1)$.
	\end{theorem}

	\begin{remark}[Sharpness] \normalfont
		We note that in Theorem \ref{HD}, the limit case $\sigma=2s-1$, $q=1$ is in general unattainable, since the fundamental solution of the fractional Laplacian $(-\Delta)^s$ does not belong to $W^{2s,1}_{\loc}(\Omega)$, see \cite[Remark 1.7]{KuuSimYan22}.
	\end{remark}
	
	In view of the Lebesgue differentiation theorem, Theorem \ref{HD} in particular implies that the set of non-Lebesgue points of $\nabla u$ has Lebesgue measure zero and therefore Hausdorff dimension smaller or equal than $n$.
	However, due to the differentiability gain on $\nabla u$ in Theorem \ref{HD}, we are also able to obtain the following corollary of Theorem \ref{HD}, which yields a more refined bound on the Hausdorff dimension of the set of non-Lebesgue points of $\nabla u$.
    \begin{corollary}\label{el : cor.sin}
    Let $s \in (1/2,1)$, $\mu \in \mathcal{M}(\mathbb{R}^n)$, $g \in W^{s,2}_{\loc}(\mathbb{R}^n) \cap L^1_{2s}(\mathbb{R}^n)$ and let $u$ be a SOLA of (\ref{NonlocalDir}) in a bounded domain $\Omega \subset \mathbb{R}^n$. Moreover, suppose that $\Phi$ satisfies Assumption \ref{assump} for some $\Lambda \geq 1$. In addition, we define the set
    \begin{align*}
        \Upsigma_{u}\coloneqq\left\{x\in\Omega\,:\,\liminf_{\rho\to0}E_{\mathrm{loc}}(\nabla u;B_{\rho}(x))>0\quad\text{or}\quad \limsup_{\rho\to0}|(\nabla u)_{B_{\rho}(x)}|=\infty\right\}.
    \end{align*}
	Then we have 
    \begin{equation*}
        \mathrm{dim}(\Upsigma_{u})\leq n-2s+1,
    \end{equation*}
    where we denote by $\mathrm{dim}(\Upsigma_{u})$ the Hausdorff dimension of the set $\Upsigma_{u}$.
    \end{corollary}

	For $\mu \in \mathcal{M}(\mathbb{R}^n)$, recall that the classical Riesz potential of order $\beta \in (0,n)$ of $|\mu|$ is defined by
	$$
	I^{|\mu|}_{\beta}(x_0):= \int_{\mathbb{R}^n} \frac{d|\mu|(y)}{|x_0-y|^{n-\beta}}, \quad x_0 \in \mathbb{R}^n.
	$$
	
	Since we are concerned with equations posed in general domains, for our purposes it is convenient to also define a truncated version of the classical Riesz potentials.
	\begin{definition} \label{Riesz}
		Let $\mu \in \mathcal{M}(\mathbb{R}^n)$. For any $x_0 \in \mathbb{R}^n$ and any $R>0$, we define the truncated Riesz potential of order $\beta \in (0,n)$ of $|\mu|$ by
		$$
		I^{|\mu|}_{\beta}(x_0,R):= \int_0^R \frac{|\mu|(B_t(x_0))}{t^{n-\beta}}\frac{dt}{t}.
		$$
	\end{definition}
	It is not difficult to see that the above truncated version of the Riesz potential is consistent with the classical one in the sense that for any $R>0$ and any $x_0 \in \mathbb{R}^n$,
	\begin{equation} \label{truncated}
		I^{|\mu|}_{\beta}(x_0,R) \leq c I^{|\mu|}_{\beta}(x_0)
	\end{equation}
	for some constant $c=c(n)$.
	We are finally in the position to state our main result concerning gradient potential estimates in bounded domains.

	\begin{theorem}[Gradient potential estimates in bounded domains]
		\label{el : thm.main2}
		Let $s \in (1/2,1)$, $\mu \in \mathcal{M}(\mathbb{R}^n)$, $g \in W^{s,2}_{\loc}(\mathbb{R}^n) \cap L^1_{2s}(\mathbb{R}^n)$ and let $u$ be a SOLA of (\ref{NonlocalDir}) in a bounded domain $\Omega \subset \mathbb{R}^n$. Moreover, suppose that $\Phi$ satisfies Assumption \ref{assump} for some $\Lambda \geq 1$. Then for almost every $x_0 \in \Omega$ and any $R>0$ such that $B_R(x_0) \subset \Omega$, we have the pointwise estimate
		\begin{equation} \label{NNPE}
			|\nabla u(x_{0})| \leq {c}E(u/{R};B_R(x_0))+cI_{2s-1}^{|\mu|}(x_{0},R)
		\end{equation}
		for some constant $c=c(n,s,\Lambda)$. In addition, for any fixed $s_0\in(1/2,1)$, the constants $c$  depends only on $n,s_{0}$, and $\Lambda$ whenever $s \in [s_0,1)$.
	\end{theorem}

Let us conclude this section by discussing applications of the obtained gradient potential estimates to regularity theory. In fact, it is well-known that in the realm of Calder\'on-Zygmund-type estimates, ``passing through potentials" enables us to detect finer scales that are difficult to reach by more traditional methods.

More concretely, the following fine regularity results in terms of the well-known Lorentz spaces $L^{p,q}(\Omega)$ that refine the classical $L^p$ spaces follow directly from the estimate \eqref{NNPE} by taking into account \eqref{truncated} and the mapping properties of the Riesz potential $I^{|\mu|}_{2s-1}$ given by \cite[Proposition 2.8]{Cianchi}.

\begin{corollary} \label{Lorentzreg}
	Let $s \in (1/2,1)$ and let $u$ be a SOLA of (\ref{NonlocalDir}) in a bounded domain $\Omega \subset \mathbb{R}^n$, where $\Phi$ satisfies Assumption \ref{assump} for some $\Lambda \geq 1$.
	\begin{itemize}
		\item We have the implication 
		\begin{equation} \label{Lorentz1}
			\mu \in \mathcal{M}(\mathbb{R}^n) \implies \nabla u\in L^{\frac{n}{n-2s+1},\infty}_{\loc}(\Omega).
		\end{equation}
		\item For any $p \in \left (1,\frac{n}{2s-1} \right )$ and any $q \in (0,\infty]$, we have the implication 
		\begin{equation} \label{Lorentz2}
			\mu \in L^{p,q}(\Omega) \implies \nabla u\in L^{\frac{np}{n-(2s-1)p},q}_{\loc}(\Omega).
		\end{equation}
		\item We have the Lipschitz criterion
		\begin{equation} \label{Lorentz3}
			\mu \in L^{\frac{n}{2s-1},1}(\Omega) \implies \nabla u\in L^{\infty}_{\loc}(\Omega).
		\end{equation}
	\end{itemize}
\end{corollary}
For a precise definition of Lorentz spaces and a discussion of the relations between them, we refer to \cite[Section 1.3]{KuuSimYan22}.

We shall just note that the implication \eqref{Lorentz1} sharpens Theorem \ref{HD} in the case when $\gamma=0$, which only yields the slightly weaker conclusion $\nabla u\in L^p_{\loc}(\Omega)$ for any $p <\frac{n}{n-2s+1}$.
Moreover, we note that the implication \eqref{Lorentz2} applied with $p=q$ in particular yields the following slightly coarser implication in the standard $L^p$ spaces: For any $p \in \left (1,\frac{n}{2s-1} \right )$, we have
\begin{equation} \label{Lpest}
	\mu \in L^{p}(\Omega) \implies \nabla u\in L^{\frac{np}{n-(2s-1)p}}_{\loc}(\Omega).
\end{equation}

More generally, it is well-known that gradient potential estimates of the type \eqref{NNPE} imply gradient estimates in any function space in which the mapping properties of the Riesz potentials are known, which is the case also in many other commonly used rearrangement invariant function spaces such as Orlicz spaces, see \cite{Cianchi}.

\subsection{Related results} \label{sec:pr}

Let us now compare our main results with the previous literature in a more comprehensive fashion. Concerning our $C^{1,\alpha}$ regularity result Theorem \ref{el : thm.main1}, a similar result for nonlinear nonlocal equations defined in terms of a $C^1$ nonlinearity was proved by Caffarelli, Chan and Vasseur in \cite{CCV} under the assumption that the equation holds on the whole space $\mathbb{R}^n$. Together with the $C^1$ assumption on the nonlinearity instead of a Lipschitz assumption as in our case, this allows the authors of \cite{CCV} to directly differentiate the equation and then deduce $C^{1,\alpha}$ regularity by applying a De Giorgi-type result to the first-order quotients of the solution. In our setting of equations given by Lipschitz nonlinearities posed in bounded domains, this strategy needs to be augmented by a delicate localization argument and an iteration scheme via fractional De Giorgi classes, so that Theorem \ref{el : thm.main1} can be considered to be a nontrivial extension of the $C^{1,\alpha}$ regularity result on the whole space obtained in \cite{CCV}.

Concerning our higher differentiability result under measure data given by Theorem \ref{HD}, a similar result was obtained in \cite[Theorem 1.6]{KuuSimYan22} in the case of linear nonlocal equations with H\"older coefficients. Moreover, under the stronger assumption that $\mu$ belongs to an appropriate fractional Sobolev space, in \cite{BL} higher differentiability above the gradient level was proved in the case of the fractional $p$-Laplacian with $p \geq 2$, that arises when taking $\Phi(t)=|t|^{p-2}t$ in \eqref{nonlocalop}. Nevertheless, Theorem \ref{HD} seems to be the first higher differentiability result above the gradient level for nonlinear nonlocal equations with general measure data.

Finally, Theorem \ref{el : thm.main3} and Theorem \ref{el : thm.main2} seem to be the first gradient potential estimates for nonlinear nonlocal equations recorded in the literature. Nevertheless, as already mentioned, in \cite{KuuSimYan22} gradient potential estimates were established for a large class of linear nonlocal equations that arise when the fractional Laplacian is perturbed by H\"older coefficients. Moreover, zero-order potential estimates for different types of nonlinear nonlocal equations were previously obtained in e.g.\ \cite{KuuMinSir15,KMSsurvey,KLL,Minhyun1}, while similar potential estimates on the solution itself for a class of nonlocal drift-diffusion equations related to the surface quasi-geostrophic (SQG) equation from fluid dynamics were recently proved in \cite{NNSW}. In addition, similar fine pointwise estimates of any order strictly smaller than one in terms of certain nonlocal fractional sharp maximal functions were recently provided by two of the authors in \cite{DNCZ} for the fractional $p$-Laplacian.

\subsection{Approach}

Let us give a brief heuristic overview of our approach to proving our main results and ultimately the gradient potential estimates given by Theorem \ref{el : thm.main3} and Theorem \ref{el : thm.main2}, highlighting the novelties in comparison to previous approaches to obtain gradient potential estimates for local and nonlocal equations.

In most papers concerned with obtaining gradient potential estimates in the local regime, either one of the following two methods is applied:

\begin{itemize}
	\item A De Giorgi-type iteration inspired by \cite{KM} combined with difference quotient arguments inspired by Littlewood-Paley theory allowing to differentiate measure data problems in a fractional sense (see e.g.\ \cite{Min}).
	\item A potential-theoretic Campanato-type iteration in terms of the local excess functional of the gradient (see e.g.\ \cite{DM2,DM1,BCDKS}).
\end{itemize}

In our nonlinear nonlocal setting, trying to adapt any of the above two approaches to prove gradient potential estimates leads to severe complications. This is mainly for two reasons. First of all, the unavoidable appearance of nonlocal tail terms in our setting makes it difficult to naively differentiate the equation. Secondly, the lower order of the nonlocal equations we consider leads to an absence of obvious energy estimates on the gradient, leading in particular to a lack of obvious comparison estimates at the gradient level. 

In the linear nonlocal setting with coefficients treated in \cite{KuuSimYan22}, these difficulties were surmounted by combining elements of both of the mentioned approaches. More precisely, a Campanato-type iteration in terms of affine functions inspired by the fully nonlinear setting (see e.g.\ \cite{CaffFNAnnals,CSARMA}) was combined with delicate difference quotient arguments inspired by corresponding methods developed in the context of local problems (see e.g.\ \cite{KrMin,KrMin1,MinCZMD,Min,AKMARMA,DFMInvent}) in order to overcome the mentioned difficulties arising due to the nonlocality and the lower order of the equation in comparison to the local setting.

In our setting of nonlinear nonlocal equations of the type \eqref{pt : eq.main}, it seems difficult to apply Campanato-type techniques in terms of affine functions in a similar fashion. This is because given a solution to \eqref{pt : eq.main} and an arbitrary affine function $\ell$, in contrast to the linear setting $u-\ell$ does no longer seem to solve a corresponding nonlocal equations exhibiting sufficiently strong regularity properties to be of use to obtaining precise gradient estimates. 

For this reason, in this paper we apply an even more involved hybrid approach of the mentioned two approaches originating in the local setting, which is combined with delicate localization and interpolation arguments that enable us to differentiate the equation and avoid the use of affine functions. To be more precise, our proof strategy to prove our gradient potential estimates and our other main results can be structured into the following three steps.

\textbf{Step 1:} \emph{Gradient regularity and first-order oscillation decay for homogeneous equations via fractional De Giorgi classes and localization}.

This first step consists of proving precise first-order estimates for solutions of \eqref{pt : eq.main} in the homogeneous case when $\mu=0$, as this lays the foundation to being able to prove gradient potential estimates under general measure data. In order to accomplish this, we differentiate the equation in a discrete sense, proving that truncations of difference quotients of solutions $v$ to $\setLcal v=0$ given by $\left (\frac{\delta_h v}{|h|^\beta} -k \right )_+$ satisfy certain Caccioppoli inequalities with nonlocal tails for any $\beta \in (0,1]$ and any $k \in \mathbb{R}$. This implies that these difference quotients of such solutions $v$ belong to a fractional De Giorgi class in the sense of \cite{CozziJFA}, which in view of the theory established by Cozzi in \cite{CozziJFA} implies that for any $\beta \in (0,1]$, such quotients satisfy H\"older estimates with respect to some small exponent $\alpha$ that does not depend on $\beta$. Therefore, these H\"older estimates can be iterated to obtain that the first-order quotients of $v$ given by $\frac{\delta_h v}{|h|}$ satisfy such a H\"older estimate. Since this H\"older estimate involves also first-order quotients in the nonlocal tail terms which cannot be differentiated in a traditional sense, in order to obtain $C^{1,\alpha}$ regularity and thus Theorem \ref{el : thm.main1}, the obtained H\"older estimate on quotients has to be combined with a delicate localization argument at each step of the iteration. In fact, we prove that for suitable cutoff functions $\psi$, despite the nonlinear nature of the operator \eqref{nonlocalop}, the product $v\psi$ satisfies a nonlocal equation of the type $\setLcal v=f$, where $f$ is at least bounded and in addition inherits the regularity of $v$ in H\"older spaces, making it suitable for our bootstrapping argument via fractional De Giorgi classes. In addition to implying Theorem \ref{el : thm.main1} in the whole range $s \in (0,1)$, this approach is flexible enough to obtain first-order oscillation decay estimates of the type
\begin{equation} \label{oscdecay}
	E\left(\frac{\delta_{h}{v}}{|h|};B_{\rho R}(x_{0})\right) \lesssim \rho^\alpha E\left(\frac{\delta_{h}{v}}{|h|};B_{R}(x_{0})\right)
\end{equation}
for all $\rho \in (0,1/4]$ and all increments $h$ with $|h|$ small enough. The main advantage of the estimate \eqref{oscdecay} in comparison to the estimate \eqref{Holdest} from Theorem \ref{el : thm.main1} is that in \eqref{oscdecay} oscillations are controlled by oscillations of the quotients rather than by their size, making the estimate \eqref{oscdecay} more suitable for running a Campanato-type iteration to obtain first-order potential estimates.

\textbf{Step 2:} \emph{Higher differentiability and first-order comparison estimates under measure data via a perturbative difference quotient argument and interpolation}.

Next, we turn to considering solutions $u$ to equations of the type \eqref{pt : eq.main} with general measure data in the case when $s \in (1/2,1)$. The starting point in this case is to use harmonic replacement throughout the scales, comparing $u$ at each scale to the solution $v$ of $\setLcal v=0$ with complement data given by $u$. More precisely, for any ball $B \Subset \Omega$, we consider the solution $v$ of
\begin{equation*}
	\left\{
	\begin{alignedat}{3}
		\setLcal v&= 0&&\quad \mbox{in  $B$}, \\
		v&=u&&\quad  \mbox{a.e. in }\RRn\setminus B.
	\end{alignedat} \right.
\end{equation*}
While as mentioned, at least initially no first-order comparison estimate between $u$ and $v$ are available, zero-order comparison estimates at the level of the solutions themselves are readily available due to \cite{KuuMinSir15}, in the sense that
\begin{equation*}
	\dashint_{B}|u-v|\,dx\lesssim r_B^{2s-n}|\mu|(B),
\end{equation*}
where $r_B$ denotes the radius of $B$. The key idea in order to make this comparison estimate compatible with the first-order oscillation decay estimate \eqref{oscdecay} is to localize the analysis to balls that depend on the increment $|h|$ from \eqref{oscdecay} itself. This is possible since the $|h|$-dependence of the balls can in the end always be removed by a quantitative covering lemma yielding control on the overlap of the balls, which in the nonlocal context was first realized in \cite{KuuSimYan22}. First of all, this approach allows us to prove our higher differentiability result given by Theorem \ref{HD}.

Combining further novel variations of these ideas with \eqref{oscdecay} also enables us to prove gradient oscillation decay estimates incorporating complement data given by
\begin{equation} \label{oscdecaym}
	E(\nabla v;\rho B ) \lesssim \rho^\alpha E_\loc (\nabla v;\tfrac{1}{2} B) + E(\nabla u;\tfrac{1}{2} B) + r_B^{2s-1-n}|\mu|(B)
\end{equation}
for all $\rho \in (0,1/4]$, and with the aid of an interpolation inequality of Gagliardo-Nirenberg-type also comparison estimates at the gradient level of the form
\begin{equation} \label{gcomp}
	\dashint_{\tfrac12 {B}} |\nabla u-\nabla v| \,dx \leq r_B^{2s-1-n} |\mu|(B) +\left( r_B^{2s-1-n} |\mu|(B)\right)^{1-\theta}E(\nabla u;B)^{\theta}
\end{equation}
for some $\theta \in (0,1)$.

We note that the precise form of the right-hand side of \eqref{oscdecaym} is crucial as it respects the lack of control on $\nabla v$ close to the boundary of $B$. In fact, while we can assume that $u \in W^{1,1}(\mathbb{R}^n)$ in view of Theorem \ref{HD} and localization arguments, we cannot simply apply a similar localization argument to $v$, which is because the precise structure of the estimate \eqref{oscdecay} is not invariant under localization.

\textbf{Step 3:} \emph{Gradient potential estimates via a potential-theoretic Campanato-type iteration and localization}.

Combining the estimates \eqref{oscdecaym} and \eqref{gcomp} then enables us to prove excess decay estimates of the form
\begin{equation} \label{edc}
	\begin{aligned}
	E(\nabla u;\rho B) & \lesssim \rho^{\alpha} E(\nabla u;B)\\
	& \quad + \rho^{-n} \left(\frac{|\mu|(B)}{r_B^{n-2s+1}}\right)^{1-\theta} E\left(\nabla u;B \right)^\theta \\ & \quad + \rho^{-n} \frac{|\mu|(B)}{r_B^{n-2s+1}}.
	\end{aligned}
\end{equation}

Observing that the contribution of the excess of $\nabla u$ in the second term on the right-hand side of \eqref{edc} can be made arbitrarily small in view of Young's inequality, by an adaption of the potential-theoretic Campanato-type iteration introduced in \cite{DM2}, for any ball $B_R(x_0) \subset \Omega$ we arrive at the pointwise gradient estimate
\begin{align*}
	|\nabla u(x_0)| \lesssim E(\nabla u;B_R(x_0))+I^{|\mu|}_{2s-1}(x_{0},R)
\end{align*}
whenever $x_0 \in \Omega$ is a Lebesgue point of $\nabla u$.
The proof of Theorem \ref{el : thm.main2} is now finished in view of Theorem \ref{HD} and another application of our localization lemma in order to remove the gradient in the tail that appears in the above estimate, while Theorem \ref{el : thm.main3} follows simply by letting $R \to \infty$ in Theorem \ref{el : thm.main2}.

\subsection{Some open questions}

Let us conclude our introduction by discussing some open questions related to this work that we consider to be interesting.

\textbf{Nonlinear nonlocal equations with coefficients}: Since in \cite{KuuSimYan22} gradient potential estimates were proved for a related class of linear nonlocal equations with coefficients, a natural question is if our main results remain valid also if our nonlinear nonlocal operator \eqref{nonlocalop} is perturbed by coefficients. This is in particular since the higher regularity theory of nonlinear nonlocal equations with coefficients below the gradient level is by now well understood, see for instance \cite{Sob,FallCalcVar,MSY,MeH,MeV,MeI,MeD,FMSYPDEA,DNCZ,Kyeongbae1,Kyeongbae2,Kyeongbae3} for a non-exhaustive list of noteworthy contributions in this direction.

\textbf{Nonlinear nonlocal equations with $p$-growth}: Another interesting question is if there are counterparts of our main results in the case when the nonlinearity $\Phi$ in \eqref{nonlocalop} more generally satisfies suitable growth assumptions of the type $\Phi(t) \approx t^{p-1}$ for some $p \in (1,\infty)$. This concerns in particular the model case of the fractional $p$-Laplacian given by the case $\Phi(t)=|t|^{p-2}t$. To the best of our knowledge, analogues of our main results are not known in each case unless $p=2$.

On the other hand, H\"older regularity below the gradient level for the fractional $p$-Laplacian was proved in \cite{BLS,GL23}, while fine higher regularity results below the gradient level were obtained in \cite{DNCZ}.
Moreover, the classical $C^{1,\alpha}$ estimates for the local $p$-Laplacian were proved in \cite{NU68,KU77}, while Riesz potential estimates for the local $p$-Laplacian were established in \cite{KuMiARMA1}. \bigskip



\subsection{Outline}

The organization of this paper is as follows. In Section \ref{sec:2}, we provide basic notation, function spaces and auxiliary lemmas which will be used frequently throughout the paper. Section \ref{sec:3} is devoted to prove $C^{1,\alpha}$-regularity for weak solutions of \eqref{pt : eq.main} in the homogeneous case when $\mu\equiv 0$. In Section \ref{sec:4} we then prove gradient oscillation decay for homogeneous equations as well as higher differentiability of the gradient under measure data.
Finally, in Section \ref{sec:5} we establish comparison estimates at the gradient level for solutions to \eqref{pt : eq.main} and prove potential estimates for the gradient of any SOLA to \eqref{pt : eq.main}.

\begin{ack} \normalfont
	We thank Tuomo Kuusi and Yannick Sire for useful discussions concerning the topic of the present work.
\end{ack}

\section{Preliminaries}\label{sec:2}

\subsection{Some notation}
First of all, throughout this paper $c$ denotes general positive constants which could vary line by line. Moreover, we use parentheses to emphasize relevant dependencies on parameters, so that for example, $c=c(n,s,\Lambda)$ means that $c$ depends only on $n,s$ and $\Lambda$. 

For $x_0\in\setR^n$ and $r>0$, we denote the open ball $B_r(x_0)=\{y\in\setR^n:|y-x_0|<r\}$. We shall omit the center of the ball and simply write $B_r$ if $x_0=0$. 

For $U\subset\setR^n$, we define the indicator function of $U$ as
\begin{align*}
\bfchi_{U}(x):=
\begin{cases}
1\quad\text{if }x\in U\\
0\quad\text{if }x\in \setR^n\setminus U.
\end{cases}
\end{align*}
Given a measurable function $g:\setR^n\rightarrow\setR$, we write
\begin{align*}
g_{\pm}(x):= \max\{\pm g(x),0\}.
\end{align*}
If $g$ is integrable over a measurable set $U\subset\setR^n$ with $U$ having positive measure, i.e., $0<|U|<\infty$, then we denote by the integral average of $g$ over $U$ 
\begin{align*}
(g)_U:=\dashint_{U}g\,dx=\dfrac{1}{|U|}\int_{U}g\,dx.
\end{align*}

In addition, given a signed Radon measure $\mu$ on $\mathbb{R}^n$, as usual we define the variation of $\mu$ as the measure defined by
$$ |\mu|(E):=\mu^+(E) + \mu^-(E), \quad E \subset \mathbb{R}^n \text{ measurable},$$
where $\mu^+$ and $\mu^-$ are the positive and negative parts of $\mu$, respectively. In the case when $|\mu|(\mathbb{R}^n)<\infty$, then we say that $\mu$ has finite total variation or finite total mass.

Finally, given a domain $\Omega \subset \mathbb{R}^n$, throughout the paper we conceptualize functions $g \in L^1(\Omega)$ as signed Radon measures on $\mathbb{R}^n$ by extending $g$ by $0$ to $\mathbb{R}^n$ if necessary and denoting
$$ g(E):=\int_E g dx, \quad  E \subset \mathbb{R}^n \text{ measurable}.$$
Note that in this case for any measurable set $E \subset \mathbb{R}^n$, we have $$|g|(E)=\int_E |g| dx.$$

In addition, for any domain $\Omega \subset\setR^n$, $s\in(0,1)$ and $p\geq 1$, the fractional Sobolev space $W^{s,p}(\Omega)$ is defined as the set of all functions $g:\Omega\rightarrow\setR^n$ with
\begin{align*}
\|g\|_{W^{s,p}(\Omega)} &\coloneqq \|g\|_{L^p(\Omega)}+[g]_{W^{s,p}(\Omega)}\\
&:=\left(\int_{\Omega}|g|^p\,dx\right)^{\frac{1}{p}}+\left(\int_{\Omega}\int_{\Omega}\dfrac{|g(x)-g(y)|^p}{|x-y|^{n+sp}}\,dxdy\right)^{\frac{1}{p}}<\infty.
\end{align*}
We also consider the corresponding local spaces given by 
\begin{align*}
W^{s,p}_{\loc}(\Omega):=\{g\in L^p_{\loc}(\Omega): g\in W^{s,p}(K)\,\,\text{for any compact } K\subset \Omega\}.
\end{align*}
For more details concerning fractional Sobolev spaces, see for instance \cite{DinPalVal12}.

\subsection{Some elementary estimates} 

We also frequently use the following straightforward lemma.
\begin{lemma}
	\label{el : lem.ave}
	Let $g\in L^{1}(B_{R}(x_{0}))$. For any $c\in\bbR$, we have 
	\begin{align*}
		\int_{B_{R}(x_{0})}|g-(g)_{B_{R}(x_{0})}|\,dx\leq 2\int_{B_{R}(x_{0})}|g-c|\,dx.
	\end{align*}
\end{lemma}

By following the same lines as in the proof of \cite[Lemma 2.9]{KuuSimYan22}, we obtain the following tail estimate which is frequently used in the remaining sections.
\begin{lemma}
\label{el : lem.tail}
Let $g\in L^{1}_{2s}(\bbR^{n})$. Then there is a constant $c=c(n)$ such that
\begin{align*}
&\mathrm{Tail}(g-(g)_{B_{R}(x_{0})};B_{R}(x_{0}))\\
&\quad\leq cs^{-1}\sum_{i=0}^{i_{0}}2^{-2si}E_{\mathrm{loc}}(g;B_{2^{i}R}(x_{0}))+c2^{-2si_{0}}\mathrm{Tail}(g-(g)_{B_{2^{i_{0}}R}(x_{0})};B_{2^{i_{0}}R}(x_0)).
\end{align*}
\end{lemma}

\subsection{Some embedding results}

We start this subsection with providing the following fractional Sobolev-Poincar\'e inequality (see \cite[Theorem 6.7]{DinPalVal12}).
\begin{lemma}
	\label{el : lem.emb.poi}
	Let $g\in W^{\gamma,q}(B_{R}(x_{0}))$ with $\gamma\in(0,1)$ and $q\geq 1$ with $n>q\gamma$. Then we have 
	\begin{align*}
		\left(        \dashint_{B_{R}(x_{0})}|g-(g)_{B_{R}(x_{0})}|^{\frac{nq}{n-q\gamma}}\,dx\right)^{\frac{n-q\gamma}{nq}}\leq c[g]_{W^{\gamma,q}(B_{R}(x_{0}))}
	\end{align*}
	for some constant $c=c(n,\gamma,q)$.
\end{lemma}

Next, we deduce the following kind of Sobolev embedding lemma.
\begin{lemma}
	\label{el : lem.embdiff}
	Let $g\in W^{\widetilde{\gamma},1}(B_{R}(x_{0}))$ with $\widetilde{\gamma}\in(0,1)$.
	Let us fix ${\gamma}\in(0,\widetilde{\gamma})$ and $q\in[1,\infty)$ such that 
	\begin{equation}
 \label{el : ineq.embdi}
		{\gamma}-n/q\leq\widetilde{\gamma}-n.
	\end{equation}
	Then we have 
 \begin{equation}\label{el : emb1.embd}
	\begin{aligned}
		R^{-n/q}\|g\|_{L^{q}(B_{R}(x_{0}))}+ R^{-n/q+{\gamma}}[g]_{W^{{\gamma},q}(B_{R}(x_{0}))}&\leq cR^{-n}\|g\|_{L^{1}(B_{R}(x_{0}))}\\
		&\quad+cR^{-n+\widetilde{\gamma}}[g]_{W^{{\widetilde{\gamma}},1}(B_{R}(x_{0}))}
	\end{aligned}
 \end{equation}
	for some constant $c=c(n,\gamma,\widetilde{\gamma},q)$.
\end{lemma}
\begin{proof}
Since the above inequality is scaling invariant, we may assume $R=1$ and $x_0=0$. Let us fix $\gamma\in(0,\widetilde{\gamma})$ and $q\geq1$ satisfying \eqref{el : ineq.embdi}. Then there is $\epsilon\geq0$ such that $\gamma+\epsilon-n/q=\widetilde{\gamma}-n$.
We first note from \cite[Theorem 5.4]{DinPalVal12} that there exists a function $\widetilde{g}\in W^{\widetilde{\gamma},1}(\bbR^{n})$ such that $\widetilde{g}(x)=g(x)$ for all $x\in B_{1}$ and 
\begin{equation*}
    \|\widetilde{g}\|_{W^{\widetilde{\gamma},1}(\bbR^{n})}\leq c\|g\|_{W^{\widetilde{\gamma},1}(B_{1})}
\end{equation*}
for some constant $c=c(n,\widetilde{\gamma})$. By \cite[Proposition in Section 2.1.2]{RunSic96} and \cite[Equation (1.301)]{Tr3}, we observe 
\begin{align*}
    \|\widetilde{g}\|_{W^{{\gamma}+\epsilon,q}(\bbR^{n})}\approx \|\widetilde{g}\|_{B^{{\gamma}+\epsilon}_{q,q}(\bbR^{n})}\leq \|\widetilde{g}\|_{B^{\widetilde{\gamma}}_{1,1}(\bbR^{n})}\approx  \|\widetilde{g}\|_{W^{\widetilde{\gamma},1}(\bbR^{n})}\leq c\|g\|_{W^{\widetilde{\gamma},1}(B_{1})}
\end{align*}
for some constant $c=c(n,\widetilde{\gamma},\gamma,q)$, where we denote $B^{\gamma}_{q,q}(\bbR^{n})$ the standard Besov space given in \cite{RunSic96,Tr3}. Here, for any constant $a,b\geq0$, $a\approx b$ means that there is a constant $c\geq1$ such that $a/c\leq b\leq ca$. Using the above inequality along with the fact that $\widetilde{g}(x)=g(x)$ for all $x\in B_{1}$, we have
\begin{equation*}
    \|{g}\|_{W^{{\gamma}+\epsilon,q}(B_{1})}\leq c\|g\|_{W^{\widetilde{\gamma},1}(B_{1})}.
\end{equation*}
We note that
\begin{equation}
\label{el : ineq2.embdi}
\begin{aligned}
    [{g}]^{q}_{W^{\gamma,q}(B_{1})}=\int_{B_{1}}\int_{B_1}\frac{|{g}(x)-{g}(y)|^{q}}{|x-y|^{n+q\gamma}}\,dx\,dy&\leq 2^{q\epsilon}\int_{B_{1}}\int_{B_1}\frac{|{g}(x)-{g}(y)|^{q}}{|x-y|^{n+q(\gamma+\epsilon)}}\,dx\,dy\\
    &\leq c [{g}]^{q}_{W^{\gamma+\epsilon,q}(B_{1})}
\end{aligned}
\end{equation}
for some $c=c(q)$, where we have used that $|x-y|\leq2$ for any $x,y\in B_1$. Combining the above three inequalities, we obtain \eqref{el : emb1.embd}, which completes the proof.
\end{proof}



We also give an embedding lemma from fractional Sobolev spaces to H\"older spaces (see \cite[Proposition 2.2]{MeI}).
\begin{lemma}
	\label{el : lem.emb.hol}
	Let $g\in W^{\gamma,q}(B_{R}(x_{0}))$ with $\gamma\in(0,1)$ and $q\geq1$ with $\gamma-\frac{n}{q}>0$. Then we have 
	\begin{align*}
		[g]_{C^{0,\gamma-\frac{n}{q}}\left(B_{R}(x_{0})\right)}\leq c[g]_{W^{\gamma,q}(B_{R}(x_{0}))}
	\end{align*}
	for some constant $c=c(n,q,\gamma)$.
\end{lemma}

\subsection{Some properties of difference quotients}
For any measurable function $g:\setR^n\rightarrow\setR$ and $h\in\setR^n$, let us write
\begin{align*}
g_{h}(x)=g(x+h),\quad\delta_h g(x)=g_{h}(x)-g(x), \quad \delta^2_h g:=\delta_h(\delta_h g).
\end{align*}
Then we have the following lemma that will be useful in Section \ref{sec:4}.
\begin{lemma}
	\label{el : thm.main4}
	Let $B_{R}(x_{0})\subset \bbR^{n}$. Let $g\in W^{1,1}_{\mathrm{loc}}(\bbR^{n})$ with $\nabla g\in L^{1}_{2s}(\bbR^{n},\bbR^{n})$. Then for any $h\in\bbR$ with $0<|h|<\frac{R}{4}$, we have 
	\begin{align}\label{eq:osc.est1}
	E_{\mathrm{loc}}\left(\frac{\delta_{h}g}{|h|};B_{R}(x_{0})\right)\leq  c E_{\mathrm{loc}}\left(\nabla g;B_{R+2|h|}(x_{0})\right)
	\end{align}
	and 
	\begin{align}\label{eq:osc.est2}
	\mathrm{Tail}\left(\frac{\delta_{h}g}{|h|}-\left(\frac{\delta_{h}g}{|h|}\right)_{B_{R}(x_{0})};B_{R}(x_{0})\right)&\leq \frac{c}{s}E\left(\nabla g;B_{R+2|h|}(x_{0})\right)
	\end{align}
	for some constant $c=c(n)$.
\end{lemma}

\begin{proof}
	It suffices to give the proof when $x_{0}=0$ and $R=1$, since by translation and scaling we can obtain \eqref{eq:osc.est1} and \eqref{eq:osc.est2} for all $R>0$ and $x_0\in\setR^n$. We also assume that $g\in C_{\mathrm{loc}}^{1}(\bbR^{n})$ by standard approximation arguments. Let us write $D_{e}g(x)=\nabla g(x)\cdot e$, where $e=\frac{h}{|h|}$ is any unit vector. 
	We now prove \eqref{eq:osc.est1}. We observe that in view of the fundamental theorem of calculus and Fubini's theorem,
	\begin{equation} \label{pt : fundamental}
	\begin{aligned}
	&E_{\mathrm{loc}}\left(\frac{\delta_{h}g}{|h|};B_1\right) \\&=\dashint_{B_{1}}\left|\int_{0}^{1}\nabla g(x+th)\cdot \frac{h}{|h|}\,dt-\left(\int_{0}^{1}\nabla g(\cdot+th)\cdot \frac{h}{|h|}\,dt\right)_{B_{1}}\right|\,dx\\
	&=\dashint_{B_{1}}\left|\int_{0}^{1}D_{e}g(x+th)\,dt-\left(\int_{0}^{1}D_{e}g(\cdot+th)\,dt\right)_{B_{1}}\right|\,dx\\
	&=\dashint_{B_{1}}\left|\int_{0}^{1}D_{e}g(x+th)-\left(D_{e}g(\cdot+th)\right)_{B_{1}}\,dt\right|\,dx.
	\end{aligned}
	\end{equation}
	Using Fubini's theorem again and a change of variables, we obtain
	\begin{align*}
	E_{\mathrm{loc}}\left(\frac{\delta_{h}g}{|h|};B_1\right)&=\dashint_{B_{1}}\left|\int_{0}^{1}D_{e}g(x+th)-\left(D_{e}g\right)_{B_{1}(th)}\,dt\right|\,dx\\
	&\leq\int_{0}^{1}\dashint_{B_{1}}\left|D_{e}g(x+th)-\left(D_{e}g\right)_{B_{1}(th)}\right|\,dx\,dt\\
	&\leq \int_{0}^{1}\dashint_{B_{1}(th)}\left|D_{e}g(x)-\left(D_{e}g\right)_{B_{1}(th)}\right|\,dx\,dt\\
	&\leq c\dashint_{B_{1+2|h|}}\left|D_{e}g(x)-\left(D_{e}g\right)_{B_{1+2|h|}}\right|\,dx\\
	&\leq cE_{\mathrm{loc}}(\nabla g;B_{1+2|h|})
	\end{align*}
	for some $c=c(n)$,
	where we have used that $B_{1}(th)\subset B_{1+2|h|}$ for any $t\in[0,1]$ and 
	\begin{equation}
	\label{el : ineq1.basic}
	\begin{aligned}
	\dashint_{B_{1}(th)}\left|D_{e}g(x)-\left(D_{e}g\right)_{B_{1}(th)}\right|\,dx&\leq \dashint_{B_{1}(th)}\left|D_{e}g(x)-\left(D_{e}g\right)_{B_{1+2|h|}}\right|\,dx\\
	&\quad +|\left(D_{e}g\right)_{B_{1}(th)}-\left(D_{e}g\right)_{B_{1+2|h|}}|\\
	&\leq c\dashint_{B_{1+2|h|}}\left|D_{e}g(x)-\left(D_{e}g\right)_{B_{1+2|h|}}\right|\,dx.
	\end{aligned}
	\end{equation}

	We now prove \eqref{eq:osc.est2}. Observe that
	\begin{align*}
	I&\coloneqq \mathrm{Tail}\left(\frac{\delta_{h}g}{|h|}-\left(\frac{\delta_{h}g}{|h|}\right)_{B_{1}};B_{1}\right)\\
	&=(1-s)\int_{\mathbb{R}^{n}\setminus B_{1}}\left|\int_{0}^{1}D_{e}g(y+th)-\left(D_{e}g\right)_{B_{1}(th)}\,dt\right|\frac{\,dy}{|y|^{n+2s}}\\
	&\leq (1-s)\int_{0}^{1}\int_{\mathbb{R}^{n}\setminus B_{1}(th)}|D_{e}g(y)-\left(D_{e}g\right)_{B_{1}(th)}|\frac{\,dy\,dt}{|y-th|^{n+2s}},
	\end{align*}
	where we have used Fubini's theorem and a change of variables. We further use the fact that $B_{{1}/{2}}\subset B_{1}(th)$ and
	\begin{equation*}
	|y-th|\geq \frac{|y|}{4}\quad\text{for any }y\in\bbR^{n}\setminus B_{{1}/{2}}
	\end{equation*}
	to see that
	\begin{align}\label{el : ineq1.basic2}
	\begin{split}
	I&\leq c(1-s)\int_{0}^{1}\int_{\mathbb{R}^{n}\setminus B_{{1}/{2}}}|D_{e}g(y)-\left(D_{e}g\right)_{B_{1}(th)}|\frac{\,dy\,dt}{|y|^{n+2s}}\\
	&\leq c(1-s)\int_{\mathbb{R}^{n}\setminus B_{{1}/{2}}}|D_{e}g(y)-\left(D_{e}g\right)_{B_{1+2|h|}}|\frac{\,dy\,dt}{|y|^{n+2s}}\\
	&\quad+c(1-s)\int_{0}^{1}\int_{\mathbb{R}^{n}\setminus B_{{1}/{2}}}|\left(D_{e}g\right)_{B_{1+2|h|}}-\left(D_{e}g\right)_{B_{1}(th)}|\frac{\,dy\,dt}{|y|^{n+2s}}.
	\end{split}
	\end{align}
	Next, using $B_{1}(th)\subset B_{1+2|h|}$ we deduce that
	\begin{align}\label{el : ineq1.basic3}
	\begin{split}
	&\int_{0}^{1}\int_{\mathbb{R}^{n}\setminus B_{{1}/{2}}}|\left(D_{e}g\right)_{B_{1+2|h|}}-\left(D_{e}g\right)_{B_{1}(th)}|\frac{\,dy\,dt}{|y|^{n+2s}}\\
	&\quad\leq c\int_{0}^{1}\int_{\mathbb{R}^{n}\setminus B_{{1}/{2}}}\dashint_{B_{1+2|h|}}|D_eg-\left(D_{e}g\right)_{B_{1+2|h|}}|\,dx\frac{\,dy\,dt}{|y|^{n+2s}}\\
	&\quad\leq 
	\dfrac{c(n)}{s}\dashint_{B_{1+2|h|}}|D_eg-\left(D_{e}g\right)_{B_{1+2|h|}}|\,dx,
	\end{split}
	\end{align}
	where for the last inequality we have used the fact that 
	\begin{equation*}
	\int_{\bbR^{n}\setminus B_{{1}/{2}}}\frac{\,dy}{|y|^{n+2s}}\leq \frac{c(n)}{s}.
	\end{equation*}
	Summing up \eqref{el : ineq1.basic2} and \eqref{el : ineq1.basic3}, we have 
	\begin{align*}
	I&\leq \dfrac{c}{s}\dashint_{B_{1+2|h|}}|D_{e}g-(D_{e}g)_{B_{1+2|h|}}|\,dx+c\mathrm{Tail}(D_{e}g-(D_{e}g)_{B_{1+2|h|}};B_{\frac{1}{2}})\\
	&\leq \dfrac{c}{s}E(\nabla g;B_{1+2|h|})
	\end{align*}
	for some constant $c=c(n)$, which completes the proof.
\end{proof}

Next, we mention an embedding lemma related to first-order difference quotients (see \cite[Proposition 2.7]{BL}) as well as an embedding lemma for second-order difference quotients (see \cite[Proposition 2.4]{BL}).
\begin{lemma}
	\label{el : lem.firq}
	Let $q\in[1,\infty)$ and let us fix $h_{0}\in(0,\infty)$. Let $g\in L^{q}(\bbR^{n})$ satisfy 
	\begin{equation*}
		\sup_{0<|h|<h_{0}}\int_{\bbR^{n}}\frac{|\delta_{h}g|^{q}}{|h|^{q\gamma }}\,dx<\infty
	\end{equation*}
	for some constant $\gamma\in(0,1)$. For any $\widetilde{\gamma}\in(0,\gamma)$, we get
	\begin{align*}
		[g]^{q}_{W^{\widetilde{\gamma},q}(\bbR^{n})}\leq c\left(\frac{h_{0}^{q(\gamma-\widetilde{\gamma})}}{\gamma-\widetilde{\gamma}}\sup_{0<|h|<h_{0}}\int_{\bbR^{n}}\frac{|\delta_{h}g|^{q}}{|h|^{q\gamma }}\,dx+\frac{h_{0}^{-q\widetilde{\gamma}}}{\widetilde{\gamma}}\|g\|^{q}_{L^{q}(\bbR^{n})}\right)
	\end{align*}
	for some constant $c=c(n,q)$.
\end{lemma}

\begin{lemma}
	\label{el : lem.sirq}
	Let $q\in[1,\infty)$ and $\gamma\in(0,1)$. Suppose that $g\in L^{q}(\bbR^{n})$ satisfy 
	\begin{align*}
		\sup_{0<|h|<\infty}\int_{\bbR^{n}}\frac{|\delta_{h}^{2}g|^{q}}{|h|^{q(1+\gamma)}}\,dx<\infty.
	\end{align*}
	Then we have 
	\begin{equation*}
		\|\nabla g\|_{L^{q}(\bbR^{n})}\leq c\|g\|_{L^{q}(\bbR^{n})}+\frac{c}{\gamma}\left(\sup_{0<|h|<\infty}\int_{\bbR^{n}}\frac{|\delta_{h}^{2}g|^{q}}{|h|^{q(1+\gamma)}}\,dx\right)^{\frac{1}{q}}
	\end{equation*}
	and
	\begin{equation*}
		\sup_{0<|h|<\infty}\left\|\frac{\delta_{h}(\nabla g)}{|h|^{\gamma}}\right\|_{L^{q}(\bbR^{n})}\leq \frac{c}{\gamma(1-\gamma)}\left(\sup_{0<|h|<\infty}\int_{\bbR^{n}}\frac{|\delta_{h}^{2}g|^{q}}{|h|^{q(1+\gamma)}}\,dx\right)^{\frac{1}{q}}
	\end{equation*}
	for some constant $c=c(n,q)$.
\end{lemma}



With the above two lemmas, we can prove the following ones which will be useful in Section \ref{sec:4}.

\begin{lemma}
	\label{el : lem.emb}
	Let us fix $q\in[1,\infty)$, $R>0$ and $h_{0}\in(0,R)$.
	Let $g\in W^{1,q}(B_{R+6h_{0}}(x_{0}))$ satisfy
	\begin{equation}
		\label{el : ass.emb}
		\left(\sup_{0<|h|<h_{0}}\int_{B_{R+4h_{0}}(x_{0})}\frac{|\delta^{2}_{h}g|^{q}}{|h|^{q(1+\gamma)}}\,dx\right)^{\frac{1}{q}}<M
	\end{equation}
	for some constants $M>0$ and $\gamma\in(0,1)$.
	Then for any $\widetilde{\gamma}\in(0,\gamma)$, we have 
	\begin{align*}
		[\nabla g]^{q}_{W^{\widetilde{\gamma},q}\left(B_{R}(x_{0})\right)}&\leq \frac{ch_{0}^{q(\gamma-\widetilde{\gamma})}M^{q}}{(\gamma-\widetilde{\gamma})\gamma^{q}(1-\gamma)^{q}}\\
		&\quad+\frac{ch_{0}^{q(\gamma-\widetilde{\gamma})}}{\widetilde{\gamma}(\gamma-\widetilde{\gamma})\gamma^{q}(1-\gamma)^{q}}\frac{(R+4h_{0})^{q+n}}{h_{0}^{q(1+{\gamma})}}E_{\mathrm{loc}}^{q}(\nabla g;B_{R+4h_0}(x_{0}))^{q},
	\end{align*}
	where $c=c(n,q)$.
\end{lemma}
\begin{proof}
	We can assume $x_{0}=0$. In addition, using standard approximation arguments, we may assume $g\in C^{1}(B_{R+6h_0})$.  We next take $\xi \in C_{c}^{\infty}(B_{R+h_{0}})$ with $\xi\equiv 1$ on $B_{R+h_{0}/2}$ satisfying
	\begin{equation}
		\label{el : test.emb}
		|\nabla \xi|\leq \frac{c}{h_{0}}\quad\text{and}\quad |\nabla^{2}\xi|\leq \frac{c}{h_{0}^{2}}
	\end{equation}
	for some constant $c=c(n)$, and let $G(x)=g(x)-(g)_{B_{R+4h_{0}}}-(\nabla g)_{B_{R+4h_{0}}}\cdot x$ so that $\delta_{h}^{2}g=\delta_{h}^{2}G$.
 By Poincar\'e's inequality, we have 
\begin{equation}\label{el : gpoin.emb}	
 \begin{aligned}
		&\|G\|^{q}_{L^{q}(B_{R+4h_{0}})}\\
  &\quad=\int_{B_{R+4h_{0}}}|g(x)-(\nabla g)_{B_{R+4h_{0}}}\cdot x-(g(x)-(\nabla g)_{B_{R+4h_{0}}}\cdot x)_{B_{R+4h_{0}}}|^{q}\,dx\\
		&\quad\leq c(R+4h_{0})^{n+q}\dashint_{B_{R+4h_{0}}}|\nabla g-(\nabla g)_{B_{R+4h_{0}}}|^{q}\,dx
	\end{aligned}
 \end{equation}
	for some constant $c=c(n,q)$. Also, we estimate
\begin{equation}\label{el : dgpoin.emb}	
 \begin{aligned}
		\norm{\nabla G}^q_{L^q(B_{R+4h_0})}&= \int_{B_{R+4h_0}}\abs{\nabla g-(\nabla g)_{B_{R+4h_0}}}^q\,dx\\
		&\leq c(R+h_0)^n\dashint_{B_{R+4h_0}}|\nabla g-(\nabla g)_{B_{R+4h_0}}|^q\,dx\\
		&\leq c\dfrac{(R+4h_0)^{n+q}}{h_0^q}\dashint_{B_{R+4h_0}}|\nabla g-(\nabla g)_{B_{R+4h_0}}|^q\,dx.
	\end{aligned}
 \end{equation}
	We next observe 
	\begin{align*}
		\delta_{h}^{2}(G\xi)=\delta_{h}(\delta_{h}G \xi_{h}+G\delta_{h}\xi)&=\delta_{h}^{2}G\xi_{2h}+\delta_{h}G\delta_{h}\xi_{h}+\delta_{h}G(\delta_{h}\xi)_{h}+G\delta_{h}^{2}\xi\\
		&= \delta_{h}^{2}G\xi_{2h}+2\delta_{h}G(\xi_{2h}-\xi_{h})+G\delta_{h}^{2}\xi.
	\end{align*}
	Therefore, we have
	\begin{equation}
		\label{el : ineq1.emb}
		\begin{aligned}
			\sup_{0<|h|<\infty}\int_{\bbR^{n}}\frac{|\delta_{h}^{2}(G\xi)|^{q}}{|h|^{q(1+\gamma)}}\,dx&\leq\sup_{0<|h|<h_0}\int_{\bbR^{n}}\frac{|\delta_{h}^{2}(G\xi)|^{q}}{|h|^{q(1+\gamma)}}\,dx\\
			&\quad+3\sup_{h_0<|h|<\infty}\int_{\bbR^{n}}\frac{|\delta_{h}^{2}(G\xi)|^{q}}{|h_0|^{q(1+\gamma)}}\,dx\eqqcolon I_{1}+I_{2}.
		\end{aligned}
	\end{equation}
    Using the above observation, we further estimates $I_{1}$ as
    \begin{align*}
        I_{1}&\leq c\sup_{0<|h|<h_{0}}\int_{\bbR^{n}}\frac{|(\delta_{h}^{2}G)(x)\xi(x+2h)|^{q}}{|h|^{q(1+\gamma)}}\,dx+c\sup_{0<|h|<h_{0}}\int_{\bbR^{n}}\frac{|\delta_{h}G(\xi_{2h}-\xi_{h})|^{q}}{|h|^{q(1+\gamma)}}\,dx\\
        &\quad +c\sup_{0<|h|<h_{0}}\int_{\bbR^{n}}\frac{|G\delta_{h}^{2}\xi|^{q}}{|h|^{q(1+\gamma)}}\,dx\eqqcolon I_{1,1}+I_{1,2}+I_{1,3}.
    \end{align*}
    By \eqref{el : ass.emb}, \eqref{el : test.emb} and the fact that $\delta_{h}^{2}g=\delta_{h}^{2}G$, we have 
    \begin{equation*}
        I_{1,1}+I_{1,3}+I_{2}\leq cM^{q}+ch_{0}^{-q(1+\gamma)}\|G\|^{q}_{L^{q}(B_{R+4h_0})}
    \end{equation*}
    for some constant $c=c(n,q)$. On the other hand, by \eqref{el : test.emb} and the fundamental theorem of calculus as in \eqref{pt : fundamental}, we have
    \begin{align*}
        I_{1,2}\leq ch_{0}^{-q}\sup_{0<|h|<h_{0}}\int_{B_{R+3h_0}}\frac{|\delta_{h}G|^{q}}{|h|^{q\gamma}}\,dx\leq ch_{0}^{-q\gamma}\int_{B_{R+4h_0}}|\nabla G|^{q}\,dx
    \end{align*}
	for some constant $c=c(n,q)$. 
 Plugging the above estimates $I_{1,1},I_{1,2},I_{1,3}$ and $I_{2}$ along with \eqref{el : gpoin.emb} and \eqref{el : dgpoin.emb} into \eqref{el : ineq1.emb}, we get
 \begin{align}
 \label{el : almost.emb}
     \sup_{0<|h|<\infty}\int_{\bbR^{n}}\frac{|\delta_{h}^{2}(G\xi)|^{q}}{|h|^{q(1+\gamma)}}\,dx&\leq cM^{q}+c\frac{(R+h_0)^{n+q}}{h_{0}^{q(1+\gamma)}}E_{\mathrm{loc}}^{q}(\nabla g;B_{R+4h_0})^{q}
 \end{align}
 for some constant $c=c(n,q)$.  Using this along with Lemma \ref{el : lem.sirq}, we see
	\begin{align*}
		\sup_{0<|h|<\infty}\int_{\bbR^{n}}\frac{|\delta_{h}\nabla (G\xi)|^{q}}{|h|^{q\gamma}}\,dx&\leq \sup_{0<|h|<\infty}\frac{c}{\gamma^{q}(1-\gamma)^{q}}\left(\int_{\bbR^{n}}\frac{|\delta^{2}_{h}(G\xi)|^{q}}{|h|^{q(1+\gamma)}}\,dx\right)^{\frac{1}{q}}\\
		&\leq \frac{c}{\gamma^{q}(1-\gamma)^{q}}M^{q}\\
  &\quad+\frac{c}{\gamma^{q}(1-\gamma)^{q}}\frac{(R+h_0)^{n+q}}{h_{0}^{q(1+\gamma)}}E_{\mathrm{loc}}^{q}(\nabla g;B_{R+4h_0})^{q}
	\end{align*}
	for some constant $c=c(n,q)$. 
  We now employ Lemma \ref{el : lem.firq} to get
	\begin{align*}
		[\nabla g]^{q}_{W^{\widetilde{\gamma},q}\left(B_{R}\right)}=[\nabla G]^{q}_{W^{\widetilde{\gamma},q}\left(B_{R}\right)}&\leq [\nabla (G\xi)]^{q}_{W^{\widetilde{\gamma},q}(\bbR^{n})}\\
		&\leq c\frac{h_{0}^{q(\gamma-\widetilde{\gamma})}}{\gamma-\widetilde{\gamma}}\sup_{0<|h|<h_{0}}\int_{\bbR^{n}}\frac{|\delta_{h}\nabla (G\xi)|^{q}}{|h|^{q\gamma}}\,dx\\
		&\quad+c\frac{h_{0}^{-q\widetilde{\gamma}}}{\widetilde{\gamma}}\|\nabla (G\xi)\|^{q}_{L^{q}(\bbR^{n})}
	\end{align*}
	for some constant $c=c(n,q)$. By combining the above two estimates along with \eqref{el : gpoin.emb} and \eqref{el : dgpoin.emb}, we  obtain the desired estimate.
\end{proof}


\begin{lemma}\label{el : lem.emb1}
Let  us fix $R>0$ and $h_{0}\in(0,R)$. Let us assume $g\in L^{1}(B_{R+6h_{0}}(x_{0}))$ with
\begin{equation}
		\label{el : ass.emb1}
		h_{0}^{-1}\sup_{0<|h|<h_{0}}\int_{B_{R+4h_{0}}(x_{0})}\frac{|\delta_{h}g|}{|h|^{\gamma}}\,dx+\sup_{0<|h|<h_{0}}\int_{B_{R+4h_{0}}(x_{0})}\frac{|\delta^{2}_{h}g|}{|h|^{1+\gamma}}\,dx<M
	\end{equation}   
for some constants $M>0$ and $\gamma\in(0,1)$. Then we have $g\in W^{1,1}(B_{R+2h_0}(x_0))$ with the estimate
\begin{align*}
   \|\nabla g\|_{L^{1}(B_{R+h_{0}/2})}\leq cM+c(h_{0}^{-1-\gamma}+1)\|g\|_{L^{1}(B_{R+4h_0})}
\end{align*}
for some constant $c=c(n,\gamma)$.
\end{lemma}
\begin{proof}
    Let us assume $x_{0}=0$. We next take $\xi \in C_{c}^{\infty}(B_{R+h_{0}})$ with $\xi\equiv 1$ on $B_{R+h_{0}/2}$ satisfying \eqref{el : test.emb}. As in the above estimates of $I_{1,1},I_{1,2},I_{1,3}$ and $I_{2}$ in Lemma \ref{el : lem.emb} with $G$ replaced by $g$, we have 
    \begin{align*}
    \sup_{0<|h|<\infty}\int_{\bbR^{n}}\frac{|\delta_{h}(g\xi)|}{|h|^{1+\gamma}}\,dx&\leq cM+ch_{0}^{-1}\sup_{0<|h|<h_{0}}\int_{B_{R+3h_0}}\frac{|\delta_{h}g|}{|h|^{\gamma}}\,dx\\
    &\quad+ch_{0}^{-1-\gamma}\|g\|_{L^{1}(B_{R+4h_0})}\\
    &\leq cM+h_{0}^{-1-\gamma}\|g\|_{L^{1}(B_{R+4h_0})}
    \end{align*}
    for some constant $c=c(n)$. By Lemma \ref{el : lem.sirq} along with the fact that $\xi\equiv1$ on $B_{R+h_0/2}$, we obtain the desired result.
\end{proof}


\subsection{A covering lemma}
We shall regularly use the following simple lemma.
\begin{lemma}
	\label{el : lem.besi}
	Let $x_0 \in \mathbb{R}^n$, $R>0$ and $r\in(0,R/2)$ be given. Then there is a constant $c=c(n)$, a finite index set $I$ and a sequence $\{z_i\}_{i \in I} \subset B_R(x_0)$ such that
	\begin{equation*}
		B_{R}(x_0)\subset \bigcup_{i\in I}B_{r}(z_{i}), \quad \sup_{x \in \mathbb{R}^n} \sum_{i\in I}\bfchi_{B_{2^{k}r}(z_{i})}(x)\leq c2^{nk}, \quad |I| \leq c \frac{R^n}{r^n},
	\end{equation*}
	where we denote by $|I|$ the number of elements in the set $I$.
\end{lemma}
\begin{proof}
	We note that there is a mutually disjoint covering $\{Q_{r/\sqrt{n}}(z_i)\}_{i\in I}$ of $B_{R}(x_{0})$ such that $Q_{r/\sqrt{n}}(z_i)\subset B_{2R}(x_{0})$, where we denote by $Q_{r/\sqrt{n}}(z_i)$ a cube with center $z_i$ and radius $2r/\sqrt{n}$. Then we observe that $\{B_{r}(z_i)\}_{i\in I}$ is a covering of $B_{R}(x_0)$ and $\{B_{r/n}(z_i)\}_{i\in I}$ is a mutually disjoint set. Therefore, we have
	\begin{align}
		\label{el : ineq1.besi}
		|I||B_{r}|=\sum_{i\in I}|B_{r}|\leq n^{n}\sum_{i\in I}|B_{r/n}(z_i)|\leq n^{n}|B_{2R}|,
	\end{align}
 which gives the third inequality of the lemma.
	We are now in the position to prove 
	\begin{align}
		\label{el : ineq00.bdry}
		 \sup_{x \in \mathbb{R}^n}\sum_{i\in I}\bfchi_{B_{2^{k}r}(z_{i})}(x)\leq n^{n}2^{n(k+1)}.
	\end{align}
	Suppose there is a point $x_{0}\in \bbR^{n}$ such that $\sum_{i\in I}\bfchi_{B_{2^{k}r}(z_{i})}(x_{0})>n^{n} 2^{n(k+1)}$.
	We now denote ${I}_{0}$ the set $\{i\in I\,:\, \bfchi_{B_{2^{k}r}(z_{i})}(x_{0})=1\}$.
	Then we observe 
	\begin{equation*}
		\bigcup_{i\in{I}_{0}} B_{2^{k}r}(z_{i})  \subset B_{2^{k+1}r}(x_{0}),
	\end{equation*}
	which implies 
	\begin{align*}
		(n^{n}2^{n(k+1)}\!+\!1)|B_{2^{k}r}|\leq\sum_{i\in I_{0}}|B_{2^{k}r}(z_i)|&=2^{kn}n^{n}\sum_{i\in I_{0}}|B_{r/n}(z_i)|\leq 2^{kn}n^{n}|B_{2^{k+1}r}(x_{0})|
	\end{align*}
	where we have used the fact that $\{B_{r/n}(z_{i})\}_{i\in I}$ is a mutually disjoint set. This is a contradiction. Thus we obtain \eqref{el : ineq00.bdry}, which completes the proof.
\end{proof}

\subsection{A classical result in potential theory}
We give the following lemma which describes the Hausdorff dimension of non-Lebesgue points of regular vector-valued function (see e.g.\ \cite{AdamsHedberg,Mingione} or \cite[Proposition 2.4]{MinCZMD}).
\begin{lemma}
\label{el : lem.sin}
    Let $G\in W^{\gamma,q}_{\mathrm{loc}}(\Omega,\bbR^{n})$ with $\gamma\in(0,1)$ and $q\in[1,\infty)$ satisfying $q\gamma<n$. Then the Hausdorff dimension of 
    \begin{align*}
        \Upsigma_{G}\coloneqq\left\{x\in\Omega\,:\,\liminf_{\rho\to0}E_{\mathrm{loc}}(G;B_{\rho}(x))>0\quad\text{or}\quad \limsup_{\rho\to0}|(G)_{B_{\rho}(x)}|=\infty\right\}
    \end{align*}
    is less or equal than $n-q\gamma$.
\end{lemma}

\subsection{Fractional De Giorgi classes}

We define upper and lower level sets of any function $g\in L^{1}(\Omega)$ on $B_{\rho}(x_{0})\subset\Omega$ by 
\begin{equation}
	\label{el : defn.levelset}
	A_{\pm}(g,x_{0},\rho;k)=\left\{x\in B_{\rho}(x_{0})\,:\,(g-k)_{\pm}(x)>0\right\},
\end{equation}
where $k\in\bbR$. 

Given $g\in L^{1}_{2s}(\bbR^{n})$, we also denote
\begin{align}
\label{el : defn.widee}
    \widetilde{E}(g;B_{R}(x_{0}))=\dashint_{B_{R}(x_{0})}|g|\,dx+\mathrm{Tail}(g;B_{R}(x_{0})).
\end{align}
In particular, we observe 
\begin{equation}
\label{el : obs.e}
    \widetilde{E}(g-(g)_{B_{R}(x_{0})};B_{R}(x_{0}))=E(g;B_{R}(x_{0})).
\end{equation}

We end this section with the following local boundedness and H\"older regularity result when a given function $v:\setR^n\rightarrow\setR$ satisfies a specific condition, i.e., $v$ is in a certain fractional De Giorgi class as defined in the paper \cite{CozziJFA}.
\begin{lemma}
\label{el : lem.degio}
Let $ v\in W^{s,2}(B_{R}(x_{0}))\cap L^{1}_{2s} (\bbR^{n})$ and two constants $M\geq 0$ and $F\geq0$ are given. Assume that for any $x_{1}\in B_{{R}/{2}}(x_{0})$, $0<r<\rho\leq {R}/{4}$ and $k\in\bbR$,
	\begin{equation}
	\label{el : ineq.degiorene}
	\begin{aligned}
	&(1-s)\int_{B_{r}(x_{1})}\int_{B_{r}(x_{1})}\frac{|(v-k)_{\pm}(x)-(v-k)_{\pm}(y)|^{2}}{|x-y|^{n+2s}}\,dy\,dx\\
	&\quad+(1-s)\int_{B_{r}(x_{1})}(v-k)_{\pm}\left(x\right)\left(\int_{B_{\rho}(x_{1})}\frac{(v-k)_{\mp}\left(y\right)}{|x-y|^{n+2s}}\,dy\right)\,dx\\
	&\quad\quad\leq \frac{M\rho^{2-2s}}{(\rho-r)^{2}}\int_{B_{\rho}(x_{1})}(v-k)_{\pm}^{2}(x)\,dx\\
	&\quad\quad\quad +\frac{M(1-s)\rho^{n+2s}}{\left(\rho-r\right)^{n+2s}}\int_{\bbR^{n}\setminus B_{\rho}(x_{1})}\int_{B_{\rho}(x_{1})}(v-k)_{\pm}(y)\frac{(v-k)_{\pm}(x)}{|y-x_{1}|^{n+2s}}\,dx\,dy\\
	&\quad\quad\quad+MF^{2}\rho^{2s}|A_{\pm}(v,x_{1},\rho;k)|
	\end{aligned}
	\end{equation}
	holds. Then we have
	\begin{align}
 \label{el : res.difrac}
	&\|v\|_{L^{\infty}(B_{R/2}(x_{0}))}+R^{\gamma}[v]_{C^{0,\gamma}(B_{R/2}(x_{0}))}\leq c\widetilde{E}(v;B_{R}(x_{0}))+cR^{2s}F
	\end{align}
	for some constants $\gamma=\gamma(n,s,M)\in(0,1)$ and $c=c(n,s,M)$.
In particular, for any fixed $s_{0}\in(0,1)$, the constants $c$ and $\gamma$ mentioned above depend only on $n,s_{0}$ and $M$ whenever $s\geq s_{0}$.
\end{lemma}
\begin{proof}
Let us fix $s_{0}\in(0,s]$.
    For any $x_{1}\in B_{R/2}(x_{0})$, we have
    \begin{align*}
	\|v\|_{L^{\infty}(B_{R/8}(x_{1}))}+R^{\gamma}[v]_{C^{0,\gamma}(B_{R/8}(x_{1}))}&\leq c\left(\dashint_{B_{R/4}(x_{1})}|v|^{2}\,dx\right)^{\frac{1}{2}}\\
 &\quad+c\mathrm{Tail}(v;B_{R/4}(x_{1}))+cR^{2s}F,
	\end{align*}
 where $\gamma=\gamma(n,s_0,M)\in(0,1)$ and $c=c(n,s_0,M)$ (see \cite{CozziJFA,CKWCalcVar}).
 On the other hand, by following a standard iteration argument as in \cite[Lemma 2.1]{KuuMinSir15}, we obtain
 \begin{align*}
    \sup_{B_{R/4}(x_{1})}|v|\leq c\widetilde{E}(v;B_{R/2}(x_{1}))+cR^{2s}F
\end{align*}
for some constant $c=c(n,s_0,M)$. Combining above two inequalities along with standard covering arguments, we obtain the desired estimate.
\end{proof}


\section{Gradient H\"{o}lder regularity for homogeneous equations}\label{sec:3}
In this section, we provide estimates that imply local boundedness and H\"older continuity of the gradient of solutions to \eqref{pt : eq.main} and are uniform in $s$. Indeed, in order to obtain results that are stable as $s \to 1$, throughout this section we fix some parameter $s_0 \in (0,1)$ and assume that
\begin{equation}
\label{el : choi.s0.sec3}
s\in [s_0,1).
\end{equation}

Moreover, for the rest of this paper we assume that $\Phi$ satisfies Assumption \ref{assump}.

\begin{lemma}
\label{pt : lem.hol.basis}
Let $u\in W^{s,2}_{\mathrm{loc}}(\Omega)\cap L^{1}_{2s}(\bbR^{n})$ be a weak solution to \eqref{pt : eq.main} with $\mu\equiv 0$. Then for any $B_{R}(x_{0})\subset\Omega$, we have 
\begin{align*}
    \|u\|_{L^{\infty}(B_{R/2}(x_{0}))}+R^{\gamma_{0}}[u]_{C^{0,\gamma_{0}}(B_{R/2}(x_{0}))}\leq c\widetilde{E}(u;B_{R}(x_{0}))
\end{align*}
with $c=c(n,s_0,\Lambda)$ and $\gamma_{0}=\gamma_{0}(n,s_0,\Lambda)$.
\end{lemma}
\begin{proof}
Fix $B_{R}(x_{0})\subset\Omega$. By \cite[Theorem 6.2]{CKWCalcVar} with $f'=\Phi$, we deduce that there is a constant $c=c(n,\Lambda)$ such that \eqref{el : ineq.degiorene} holds with $M=c$ and $F=0$. In light of Lemma \ref{el : lem.degio}, we obtain the desired estimate.
\end{proof}

The following key lemma, which implies that any localized solution of \eqref{pt : eq.main} satisfies an equation with regular right-hand side, is an essential ingredient to prove our main results.
\begin{lemma}[Localization]
\label{el : lem.loc}
Let $B_{5R}(x_{0})\subset\Omega$ and let $u\in W^{s,2}_{\mathrm{loc}}(\Omega)\cap L^{1}_{2s}(\bbR^{n})$ be a weak solution to \eqref{pt : eq.main} with $\mu \in L^\frac{2n}{n+2s}_{\mathrm{loc}}(\Omega)$. Let us fix a cut off function $\xi \in C_{c}^{\infty}\left(B_{4R}(x_{0})\right)$ with $\xi\equiv 1$ on $B_{3R}(x_{0})$. Then we have that $w\coloneqq u\xi\in W^{s,2}(B_{5R}(x_{0}))\cap L^{1}_{2s}(\bbR^{n})$ is a weak solution to
\begin{equation}\label{eq : lem.loc}
    \mathcal{L}w=f+\mu \quad\text{in }B_{2R}(x_{0})
\end{equation}
for some $f\in L^{\infty}\left(B_{2R}(x_{0})\right)$ with the estimate
\begin{align}
\label{el : ineq.tail}
    \sup_{x\in B_{{2R}}(x_{0})}|f(x)|&\leq cR^{-2s}\mathrm{Tail}(u;B_{3R}(x_{0})),
\end{align}
where $c=c(n,\Lambda)$. Moreover, if $u$ additionally belongs to $C^{0,\beta}(B_{3R}(x_{0}))$ for some $\beta\in(0,1]$, then $f\in C^{0,\beta}(B_{2R}(x_0))$ with the estimate
\begin{equation}
\label{el : lem.loc.est}
\begin{aligned}
[f]_{C^{0,\beta}(B_{2R}(x_0))}&\leq cR^{-2s}[u]_{C^{0,\beta}(B_{3R}(x_{0}))}\\
&\quad+c{R^{-(2s+\beta)}}\left[{\|u\|_{L^{\infty}(B_{3R}(x_{0}))}}+{\mathrm{Tail}(u;B_{3R}(x_{0}))}\right]
\end{aligned}
\end{equation}
for some constant $c=c(n,s_0,\Lambda)$, where the constant $s_0$ is determined in \eqref{el : choi.s0.sec3}.
\end{lemma}

\begin{remark}
Note that if $u\in W^{s,2}_{\mathrm{loc}}(\Omega)\cap L^{1}_{2s}(\bbR^{n})$ is a weak solution of \eqref{pt : eq.main} with $\mu \equiv 0$, then $u$ satisfies $u\in C_{\mathrm{loc}}^{0,\gamma_{0}}(\Omega)$ for sufficiently small $\gamma_{0}\in(0,1)$ by Lemma \ref{pt : lem.hol.basis}. Here we state the above lemma for any $\beta\in(0,1]$ since it is needed for every range $(0,1]$ of $\beta$ when the iteration is applied later in Lemma \ref{el : lem.lochol}. 
\end{remark}
\begin{proof}[Proof of Lemma \ref{el : lem.loc}]
Let us fix a test function $\psi\in W^{s,2}(B_{2R}(x_{0}))$ which has compact support in $B_{2R}(x_{0})$. Then we have 
\begin{equation}
\label{el : lem.loc.eq1}
\begin{aligned}
    &(1-s)\int_{\bbR^{n}}\int_{\bbR^{n}}{\Phi}\left(\frac{w(x)-w(y)}{|x-y|^{s}}\right)\frac{\psi(x)-\psi(y)}{|x-y|^{s}}\frac{\,dx\,dy}{|x-y|^{n}} - \int_{\Omega} \mu \psi \,dx\\
    &=(1-s)\int_{\bbR^{n}}\int_{\bbR^{n}}{\Phi}\left(\frac{w(x)-w(y)}{|x-y|^{s}}\right)\frac{\psi(x)-\psi(y)}{|x-y|^{s}}\frac{\,dx\,dy}{|x-y|^{n}}\\
    &\quad-(1-s)\int_{\bbR^{n}}\int_{\bbR^{n}}{\Phi}\left(\frac{u(x)-u(y)}{|x-y|^{s}}\right)\frac{\psi(x)-\psi(y)}{|x-y|^{s}}\frac{\,dx\,dy}{|x-y|^{n}}\eqqcolon J
\end{aligned}
\end{equation}
as $u$ is a weak solution to \eqref{pt : eq.main}. Here, since $w(x)=u(x)$ in $B_{3R}(x_{0})$, $\psi\equiv 0$ on $\mathbb{R}^{n}\setminus B_{2R}(x_{0})$ and $\Phi$ is an odd function, we obtain
\begin{align*}
    J&=2(1-s)\int_{B_{2R}(x_{0})}\int_{\mathbb{R}^{n}\setminus B_{3R}(x_{0})}{\Phi}\left(\frac{w(x)-w(y)}{|x-y|^{s}}\right)\frac{\psi(x)}{|x-y|^{n+s}}\,dy\,dx\\
    &\quad-2(1-s)\int_{B_{2R}(x_{0})}\int_{\mathbb{R}^{n}\setminus B_{3R}(x_{0})}{\Phi}\left(\frac{u(x)-u(y)}{|x-y|^{s}}\right)\frac{\psi(x)}{|x-y|^{n+s}}\,dy\,dx.
\end{align*}
\noindent
As a result, the equality \eqref{el : lem.loc.eq1} can be rewritten as follows
\begin{align*}
   (1-s)\int_{\bbR^{n}}\int_{\bbR^{n}}{\Phi}\left(\frac{w(x)-w(y)}{|x-y|^{s}}\right)\frac{\psi(x)-\psi(y)}{|x-y|^{s}}\frac{\,dx\,dy}{|x-y|^{n}}=\int_{B_{2R}(x_{0})}{(f+\mu)}\psi\,dx,
\end{align*}
where 
\begin{align*}
    f(x)&=2(1-s)\int_{\mathbb{R}^{n}\setminus B_{3R}(x_{0})}{\Phi}\left(\frac{w(x)-w(y)}{|x-y|^{s}}\right)\frac{\,dy}{|x-y|^{n+s}}\\
    &\quad-2(1-s)\int_{\mathbb{R}^{n}\setminus B_{3R}(x_{0})}{\Phi}\left(\frac{u(x)-u(y)}{|x-y|^{s}}\right)\frac{\,dy}{|x-y|^{n+s}},
\end{align*}
which implies \eqref{eq : lem.loc} by considering \eqref{distsol}.

We are now in the position to prove $f\in C^{0,\beta}\left(B_{{2R}}(x_{0})\right)$. To do this, note that
\begin{align}
\label{el : lem.loc.ineq1}
|x-y|\geq \frac{|y-x_{0}|}{6}
\end{align}
for any $x\in B_{{2R}}(x_{0})$ and $y\in \bbR^{n}\setminus B_{3R}(x_{0})$.
Using this along with \eqref{pt : assmp.phi}, the fact that $\xi(x)=1$ on $x\in B_{3R}(x_{0})$ and \eqref{el : lem.loc.ineq1}, we have 
\begin{align*}
|f(x)|&\leq 2\Lambda(1-s)\int_{\bbR^{n}\setminus B_{3R}(x_{0})}\left|\frac{w(x)-w(y)-(u(x)-u(y))}{|x-y|^{s}}\right|\frac{\,dy}{|x-y|^{n+s}}\\
&\leq c(1-s)\int_{\bbR^{n}\setminus B_{3R}(x_{0})}\frac{|u(y)|}{|x_{0}-y|^{n+2s}}\,dy
\end{align*}
for any $x\in B_{{2R}}(x_0)$, where $c=c(n,\Lambda)$. So we have \eqref{el : ineq.tail} with $c=c(n,\Lambda)$. For the H\"older regularity of $f$, let us fix the points $x_{0},x_{1}\in B_{{2R}}(x_{0})$. Then we observe
\begin{equation}
\label{el : lem.loc.eq2}
\begin{aligned}
    (1-s)^{-1}(f(x_{1})-f(x_{2}))&=2\int_{\mathbb{R}^{n}\setminus B_{3R}(x_{0})}{\Phi}\left(\frac{w(x_{1})-w(y)}{|x_{1}-y|^{s}}\right)\frac{\,dy}{|x_{1}-y|^{n+s}}\\
    &\quad-2\int_{\mathbb{R}^{n}\setminus B_{3R}(x_{0})}{\Phi}\left(\frac{w(x_{2})-w(y)}{|x_{2}-y|^{s}}\right)\frac{\,dy}{|x_{2}-y|^{n+s}}\\
    &\quad-2\int_{\mathbb{R}^{n}\setminus B_{3R}(x_{0})}{\Phi}\left(\frac{u(x_{1})-u(y)}{|x_{1}-y|^{s}}\right)\frac{\,dy}{|x_{1}-y|^{n+s}}\\
    &\quad+2\int_{\mathbb{R}^{n}\setminus B_{3R}(x_{0})}{\Phi}\left(\frac{u(x_{2})-u(y)}{|x_{2}-y|^{s}}\right)\frac{\,dy}{|x_{2}-y|^{n+s}}.
\end{aligned}
\end{equation}
We first estimate the term $J_{1}$ which is given by
\begin{equation}
\label{el : defn.j1}
\begin{aligned}
    J_{1}&\coloneqq\int_{\mathbb{R}^{n}\setminus B_{3R}(x_{0})}{\Phi}\left(\frac{w(x_{1})-w(y)}{|x_{1}-y|^{s}}\right)\frac{\,dy}{|x_{1}-y|^{n+s}}\\
    &\quad-\int_{\mathbb{R}^{n}\setminus B_{3R}(x_{0})}{\Phi}\left(\frac{w(x_{2})-w(y)}{|x_{2}-y|^{s}}\right)\frac{\,dy}{|x_{2}-y|^{n+s}}.
\end{aligned}
\end{equation}
Let us write
\begin{align*}
    J_{1}&=\int_{\mathbb{R}^{n}\setminus B_{3R}(x_{0})}\left[{\Phi}\left(\frac{w(x_{1})-w(y)}{|x_{1}-y|^{s}}\right)-{\Phi}\left(\frac{w(x_{2})-w(y)}{|x_{1}-y|^{s}}\right)\right]\frac{\,dy}{|x_{1}-y|^{n+s}}\\
    &\quad+\int_{\mathbb{R}^{n}\setminus B_{3R}(x_{0})}\left[{\Phi}\left(\frac{w(x_{2})-w(y)}{|x_{1}-y|^{s}}\right)-{\Phi}\left(\frac{w(x_{2})-w(y)}{|x_{2}-y|^{s}}\right)\right]\frac{\,dy}{|x_{1}-y|^{n+s}}\\
    &\quad+\int_{\mathbb{R}^{n}\setminus B_{3R}(x_{0})}{\Phi}\left(\frac{w(x_{2})-w(y)}{|x_{2}-y|^{s}}\right)\left[\frac{\,dy}{|x_{1}-y|^{n+s}}-\frac{\,dy}{|x_{2}-y|^{n+s}}\right] \eqqcolon \sum_{i=1}^{3}J_{1,i}.
\end{align*}
We now estimate each term $J_{1,i}$ for each $i=1,2$ and 3. For the remaining argument of the proof, we first note
\begin{equation}
\label{el : ineq3.lem.loc}
|x_2-y|,|x_1-y|\geq\frac{1}{6} |x_0-y|\quad\text{for any }y\in \bbR^{n}\setminus B_{3R}(x_{0}),
\end{equation}
which follows from \eqref{el : lem.loc.ineq1}.

\textbf{Estimate of $J_{1,1}$.} Observe that
\begin{align*}
    |J_{1,1}|&\leq \int_{\mathbb{R}^{n}\setminus B_{3R}(x_{0})}\frac{|w(x_{1})-w(x_{2})|}{|x_{1}-y|^{s}}\frac{\,dy}{|x_{1}-y|^{n+s}}\\
    &\leq c[w]_{C^{\beta}(B_{3R}(x_{0}))}|x_{1}-x_{2}|^{\beta}\int_{\bbR^{n}\setminus B_{3R}(x_{0})}\frac{1}{|x_{0}-y|^{n+2s}}\,dy\\
    &\leq cs^{-1}[w]_{C^{\beta}(B_{3R}(x_{0}))}\frac{|x_{1}-x_{2}|^{\beta}}{R^{2s}}
\end{align*}
for some constant $c=c(n,\Lambda)$, where we have used \eqref{pt : assmp.phi} and \eqref{el : ineq3.lem.loc}.

\textbf{Estimate of $J_{1,2}$.}
Using again \eqref{pt : assmp.phi} and \eqref{el : ineq3.lem.loc}, it follows that
\begin{align*}
    |J_{1,2}|&\leq c\int_{\bbR^{n}\setminus B_{3R}(x_{0})}|w(x_{2})-w(y)|\left|\frac{1}{|x_{1}-y|^{s}}-\frac{1}{|x_{2}-y|^{s}}\right|\frac{\,dy}{|x_{1}-y|^{n+s}}\\
    &\leq c\int_{\bbR^{n}\setminus B_{3R}(x_{0})}\left|\frac{1}{|x_{1}-y|^{s}}-\frac{1}{|x_{2}-y|^{s}}\right|\frac{\|w\|_{L^{\infty}(B_{3R}(x_{0}))}+|w(y)|}{|x_{0}-y|^{n+s}}\,dy.
\end{align*}
Here, we note
\begin{align}
\label{el : lem.loc.ineq2}
    \left|\frac{1}{|x_{1}-y|^{s}}-\frac{1}{|x_{2}-y|^{s}}\right|=\left|\int_{|x_{1}-y|}^{|x_{2}-y|}\frac{s}{t^{s+1}}\,dt\right|\leq c\frac{|x_{2}-x_{1}|}{|x_{0}-y|^{s+1}},
\end{align}
where we have used \eqref{el : ineq3.lem.loc}. Combining the above two estimates with the fact that $\frac{|x_{2}-x_{1}|}{|y-x_{0}|}\leq\frac{|x_{2}-x_{1}|^{\beta}}{(2R)^{\beta}}$ for any $y\in \bbR^{n}\setminus B_{3R}(x_{0})$, we get
\begin{align*}
    |J_{1,2}|&\leq c|x_{2}-x_{1}|\int_{\bbR^{n}\setminus B_{3R}(x_{0})}\frac{\|w\|_{L^{\infty}(B_{3R}(x_{0}))}+|w(y)|}{|x_{0}-y|^{n+2s+1}}\,dy\\
    &\leq c\frac{|x_{2}-x_{1}|^{\beta}}{R^{2s+\beta}}\left[{s^{-1}\|w\|_{L^{\infty}(B_{3R}(x_{0}))}}+(1-s)^{-1}{\mathrm{Tail}(w;B_{3R}(x_{0}))}\right]
\end{align*}
for some constant $c=c(n,\Lambda)$.

\textbf{Estimate of $J_{1,3}$.} 
In light of \eqref{pt : assmp.phi}, \eqref{el : ineq3.lem.loc} and \eqref{el : lem.loc.ineq2} with $s$ replaced by $n+s$, we find
\begin{align*}
    |J_{1,3}|&\leq \int_{\mathbb{R}^{n}\setminus B_{3R}(x_{0})}\left(\frac{\|w\|_{L^{\infty}(B_{3R}(x_{0}))}+|w(y)|}{|x_{0}-y|^{s}}\right)\left|\frac{1}{|x_{1}-y|^{n+s}}-\frac{1}{|x_{2}-y|^{n+s}}\right|\,dy\\
    &\leq c\frac{|x_{2}-x_{1}|^{\beta}}{R^{\beta}}\int_{\mathbb{R}^{n}\setminus B_{3R}(x_{0})}\left(\frac{\|w\|_{L^{\infty}(B_{3R}(x_{0}))}+|w(y)|}{|x_{0}-y|^{s}}\right)\frac{\,dy}{|x_{0}-y|^{n+s}},
\end{align*}
which implies 
\begin{align*}
    |J_{1,3}|\leq c\frac{|x_{2}-x_{1}|^{\beta}}{R^{2s+\beta}}\left[s^{-1}{\|w\|_{L^{\infty}(B_{3R}(x_{0}))}}+(1-s)^{-1}{\mathrm{Tail}(w;B_{3R}(x_{0}))}\right]
\end{align*}
for some constant $c=c(n,\Lambda)$. We now combine all the estimates $J_{1,1},J_{1,2}$ and $J_{1,3}$ to see that
\begin{align*}
    |J_{1}|& \leq c[w]_{C^{\beta}(B_{3R}(x_{0}))}\frac{|x_{1}-x_{2}|^{\beta}}{R^{2s}}+c\frac{|x_{1}-x_{2}|^{\beta}}{R^{2s+\beta}}{\|w\|_{L^{\infty}(B_{3R}(x_{0}))}}\\
    &\quad +\dfrac{c}{1-s}\frac{|x_{1}-x_{2}|^{\beta}}{R^{2s+\beta}}{\mathrm{Tail}(w;B_{3R}(x_{0}))}
\end{align*}
for some constant $c=c(n,s_0,\Lambda)$, as $s\geq s_0$ by recalling the choice of the constant $s_{0}$ determined in \eqref{el : choi.s0.sec3}. On the other hand, define
\begin{align*}
J_2&:=-\int_{\mathbb{R}^{n}\setminus B_{3R}(x_{0})}{\Phi}\left(\frac{u(x_{1})-u(y)}{|x_{1}-y|^{s}}\right)\frac{\,dy}{|x_{1}-y|^{n+s}}\\
&\quad+\int_{\mathbb{R}^{n}\setminus B_{3R}(x_{0})}{\Phi}\left(\frac{u(x_{2})-u(y)}{|x_{2}-y|^{s}}\right)\frac{\,dy}{|x_{2}-y|^{n+s}} =\frac{1}{2}(f(x_1)-f(x_2)-2J_1).
\end{align*}
By following the same lines as in the estimate of $J_{1}$ with $w$ replaced by $u$, we have 
\begin{align*}
    |J_{2}|& \leq c[u]_{C^{\beta}(B_{3R}(x_{0}))}\frac{|x_{1}-x_{2}|^{\beta}}{R^{2s}}+c\frac{|x_{1}-x_{2}|^{\beta}}{R^{2s+\beta}}{\|u\|_{L^{\infty}(B_{3R}(x_{0}))}}\\
    &\quad +\dfrac{c}{1-s}\frac{|x_{1}-x_{2}|^{\beta}}{R^{2s+\beta}}{\mathrm{Tail}(u;B_{3R}(x_{0}))}
\end{align*}
for some constant $c=c(n,s_0,\Lambda)$. Finally, recalling \eqref{el : lem.loc.eq2}, $w=u\xi$, where $\xi\equiv1$ on $B_{3R}(x_{0})$, we have \eqref{el : lem.loc.est}, where the constant $c$ depends only on $n,s_0$ and $\Lambda$. Therefore, the proof is complete.
\end{proof}

We now give Caccioppoli-type estimates for the $\beta$-order difference quotients of solutions.
\begin{lemma}\label{lem:reg.cacc}
Let $w\in W^{s,2}(B_{R}(x_{0}))\cap L^{1}_{2s}(\bbR^{n})$ be a weak solution to 
\begin{equation}
\label{el : eq.locregnor}
    \mathcal{L}w=f\quad\text{in }B_{R}(x_{0}),
\end{equation}  
where $f\in L^{\infty}\left(B_{2R}(x_{0})\right)$.
Let us fix $0< r<\rho\leq{R}/{4}$ and $x_{1}\in B_{R/2}(x_{0})$.  We
choose a cut-off function $\psi\in C_{c}^{\infty}\left(B_{\frac{r+\rho}{2}}(x_{1})\right)$ with $\psi\equiv 1$ on $B_{r}(x_{1})$. Then there is a constant $c=c(n,\Lambda)$ such that for any $\beta \in (0,1]$, we have
\begin{equation}
\label{pt : ineq.cca.diffq}
\begin{aligned}
&(1-s)\int_{B_{\rho}(x_{1})}\int_{B_{\rho}(x_{1})}\frac{|((\widetilde{w}_{\beta}-k)_{\pm}\psi)(x)-((\widetilde{w}_{\beta}-k)_{\pm}\psi)(y)|^{2}}{|x-y|^{n+2s}}\,dx\,dy\\
    &\quad+(1-s)\int_{B_{\rho}(x_{1})}
    ((\widetilde{w}_{\beta}-k)_{\pm}\psi^{2})(x)\left(\int_{\setR^n\setminus B_{\rho}(x_{1})}\frac{(\widetilde{w}_{\beta}-k)_{\mp}\left(y\right)}{|x-y|^{n+2s}}\,dy\right)\,dx\\
    &\leq \frac{c\rho^{2(1-s)}}{(\rho-r)^{2}}\int_{B_{\rho}(x_{1})}(\widetilde{w}_{\beta}-k)_{\pm}^{2}\,dx\\
    &\quad +\frac{c(1-s)\rho^{n+2s}}{\left(\rho-r\right)^{n+2s}}\int_{\bbR^{n}\setminus B_{\rho}(x_{1})}\int_{B_{\rho}(x_{1})}(\widetilde{w}_{\beta}-k)_{\pm}(y)\frac{((\widetilde{w}_{\beta}-k)_{\pm}\psi^{2})(x)}{|y-x_{1}|^{n+2s}}\,dx\,dy\\
    &\quad+c\int_{B_{\rho}(x_{1})}|\widetilde{f}_{\beta}(\widetilde{w}_{\beta}-k)_{\pm}\psi^{2}|\,dx
\end{aligned}
\end{equation}
for any $k\in\bbR$ and $0<|h|<R/8 $, where $\widetilde{w}_{\beta}\coloneqq\frac{\delta_{h}w}{|h|^{\beta}}$ and $\widetilde{f}_{\beta}\coloneqq\frac{\delta_{h}f}{|h|^{\beta}}$.
\end{lemma}
\begin{proof}
Let us fix $k\in\bbR$ and denote
\begin{align*}
\widehat{w}:=\widetilde{w}_{\beta}-k
\end{align*}
throughout the proof. Since we choose $h\in B_{{R}/{8}}\setminus\{0\}$, we observe that for any $g\in W^{s,2}(B_{7R/8}(x_{0}))$ with compact support in $B_{{7R}/{8}}(x_{0})$, 
\begin{align*}
g_{-h}(x):=g(x-h)\in W^{s,2}(B_{R}(x_{0}))
\end{align*}
has compact support in $B_{R}(x_{0})$. We first prove \eqref{pt : ineq.cca.diffq} for $\widehat{w}_{+}$. Setting ${\kappa}:=\widehat{w}_{+}\psi^{2}$, we test \eqref{el : eq.locregnor} with $\delta_{-h}{\kappa}$, which yields
\begin{align}
\label{pt : lem.entildew.eq1}
    &\underbrace{\int_{B_{R}(x_{0})}\frac{\widetilde{f}_{\beta}\widehat{w}_{+}\psi^{2}}{1-s}\,dx}_{:=I}\nonumber\\
    &=\underbrace{\int_{\bbR^{n}}\int_{\bbR^{n}}\left({\Phi}\left(\frac{w_{h}(x)-w_{h}(y)}{|x-y|^{s}}\right)-{\Phi}\left(\frac{w(x)-w(y)}{|x-y|^{s}}\right)\right)\frac{{\kappa}(x)-{\kappa}(y)}{|h|^{\beta}|x-y|^{n+s}}\,dx\,dy}_{=:J}.
\end{align}
We now estimate each term $I$ and $J$.

\noindent
\textbf{Estimate of $I$.}
Since $\psi\equiv0$ on $\bbR^{n}\setminus B_{\rho}(x_{1})$, $I$ is bounded as
\begin{align*}
    \left|I\right|\leq \frac{1}{1-s}\int_{B_{\rho}(x_{1})}\left|\widetilde{f}_{\beta}\widehat{w}_{+}\psi^{2}\right|\,dx.
\end{align*}

\noindent
\textbf{Estimate of $J$.}
To estimate $J$, let us rewrite $J$ as
\begin{align*}
    J&=
    \int_{B_{\rho}(x_{1})}\int_{B_{\rho}(x_{1})} J_{1}\,dx\,dy+\int_{\bbR^{n}\setminus B_{\rho}(x_{1})}\int_{B_{\rho}(x_{1})} J_{2}\,dx\,dy\\
    &\quad+\int_{B_{\rho}(x_{1})}\int_{\bbR^{n}\setminus B_{\rho}(x_{1})} J_{3}\,dx\,dy,
\end{align*}
where $$J_1:=\left({\Phi}\left(\frac{w_{h}(x)-w_{h}(y)}{|x-y|^{s}}\right)-{\Phi}\left(\frac{w(x)-w(y)}{|x-y|^{s}}\right)\right)\frac{{\kappa}(x)-{\kappa}(y)}{|h|^{\beta}|x-y|^{n+s}},$$
$$J_2:=\left({\Phi}\left(\frac{w_{h}(x)-w_{h}(y)}{|x-y|^{s}}\right)-{\Phi}\left(\frac{w(x)-w(y)}{|x-y|^{s}}\right)\right)\frac{{\kappa}(x)}{|h|^{\beta}|x-y|^{n+s}},$$
$$J_3:=\left({\Phi}\left(\frac{w_{h}(x)-w_{h}(y)}{|x-y|^{s}}\right)-{\Phi}\left(\frac{w(x)-w(y)}{|x-y|^{s}}\right)\right)\frac{{\kappa}(y)}{|h|^{\beta}|x-y|^{n+s}}$$
and we have used the fact that $\kappa=0$ in $\bbR^{n}\setminus B_{\rho}(x_{1})$. We now estimate each of the terms $J_{1},J_{2}$ and $J_{3}$. 

\textbf{Estimate of $J_{1}$.} We may assume $\delta_{h}w(x)>\delta_{h}w(y)$, since if $\delta_{h}w(x)\leq\delta_{h}w(y)$, then we similarly estimate the term $J$ by changing the role of $\delta_{h}w(x)$ and $\delta_{h}w(y)$. Note that
\begin{equation*}
\delta_{h}w(x)>\delta_{h}w(y)\implies w_{h}(x)-w_{h}(y)>w(x)-w(y).
\end{equation*}
Thus, we observe from \eqref{pt : assmp.phi} that 
\begin{equation*}
    {\Phi}\left(\frac{w_{h}(x)-w_{h}(y)}{|x-y|^{s}}\right)-{\Phi}\left(\frac{w(x)-w(y)}{|x-y|^{s}}\right)\geq \frac{\delta_{h}w(x)-\delta_{h}w(y)}{\Lambda|x-y|^{s}}\geq0,
\end{equation*}
which will be frequently used in the remainder of the proof. We then consider the following two cases (a) and (b).
\begin{enumerate}
    \item In case of $\frac{\delta_{h}w}{|h|^{\beta}}(x)\geq\frac{\delta_{h}w}{|h|^{\beta}}(y)\geq k$: We first deal with the case that $\psi(x)\geq\psi(y)$. By \eqref{pt : assmp.phi} and $w_{h}(x)-w_{h}(y)>w(x)-w(y)$ together with the fact that
    \begin{equation*}
    \widehat{w}_+(x)=(\widetilde{w}_{\beta}-k)_{+}(x)=\widetilde{w}_{\beta}(x)-k\geq\widehat{w}_+(y)=(\widetilde{w}_{\beta}-k)_{+}(y)=\widetilde{w}_{\beta}(y)-k,
    \end{equation*} 
    there holds
    \begin{align}
    \label{el : lem.reg.ineq1}
        J_{1}&\geq \left({\Phi}\left(\frac{w_{h}(x)-w_{h}(y)}{|x-y|^{s}}\right)-{\Phi}\left(\frac{w(x)-w(y)}{|x-y|^{s}}\right)\right)\frac{(\widetilde{w}_{\beta}(x)-\widetilde{w}_{\beta}(y))\psi^{2}(x)}{|h|^{\beta}|x-y|^{n+s}}\nonumber\\
        &\geq\frac{1}{\Lambda}\frac{|\widehat{w}_{+}(x)-\widehat{w}_{+}(y)|^{2}\psi^{2}(x)}{|x-y|^{n+2s}}\nonumber\\
        &\geq \frac{1}{2\Lambda}\frac{|(\widehat{w}_{+}\psi)(x)-(\widehat{w}_{+}\psi)(y)|^{2}}{|x-y|^{n+2s}}-\frac{1}{\Lambda}\frac{|\widehat{w}_{+}(y)|^{2}|\psi(y)-\psi(x)|^{2}}{|x-y|^{n+2s}}.
    \end{align}
    We now assume $\psi(x)\leq\psi(y)$. Then we rewrite $J_{1}$ as 
    \begin{align*}
        J_{1}&=\left({\Phi}\left(\frac{w_{h}(x)-w_{h}(y)}{|x-y|^{s}}\right)-{\Phi}\left(\frac{w(x)-w(y)}{|x-y|^{s}}\right)\right)\frac{(\widehat{w}(x)-\widehat{w}(y))\psi^{2}(y)}{|h|^{\beta}|x-y|^{n+s}}\\
        &\quad+\left({\Phi}\left(\frac{w_{h}(x)-w_{h}(y)}{|x-y|^{s}}\right)-{\Phi}\left(\frac{w(x)-w(y)}{|x-y|^{s}}\right)\right)\frac{\widehat{w}(x)(\psi^{2}(x)-\psi^{2}(y))}{|h|^{\beta}|x-y|^{n+s}}\\
        &\eqqcolon J_{1,1}+J_{1,2}.
    \end{align*}
As in \eqref{el : lem.reg.ineq1}, we estimate $J_{1,1}$ as 
\begin{align*}
    J_{1,1}\geq\frac{1}{\Lambda}\frac{|\widehat{w}_{+}(x)-\widehat{w}_{+}(y)|^{2}\psi^{2}(y)}{|x-y|^{n+2s}}&\geq \frac{1}{2\Lambda}\frac{|(\widehat{w}_{+}\psi)(x)-(\widehat{w}_{+}\psi)(y)|^{2}}{|x-y|^{n+2s}}\\
    &\quad-\frac{1}{\Lambda}\frac{|\widehat{w}_{+}(x)|^{2}|\psi(y)-\psi(x)|^{2}}{|x-y|^{n+2s}}.
\end{align*}
Using \eqref{pt : assmp.phi}, $\psi(x)\leq\psi(y)$, Young's inequality and $\widehat{w}_+(x)\geq\widehat{w}_+(y)$ from $\frac{\delta_{h}w}{|h|^{\beta}}(x)\geq\frac{\delta_{h}w}{|h|^{\beta}}(y)\geq k$, we next estimate $J_{1,2}$ as 
\begin{align*}
    J_{1,2}&\geq -\Lambda\frac{|\delta_{h}w(x)-\delta_{h}w(y)|\widehat{w}(x)|\psi^{2}(x)-\psi^{2}(y)|}{|h|^{\beta}|x-y|^{n+2s}}\\
    &\geq -2\Lambda\frac{|\delta_{h}w(x)-\delta_{h}w(y)|\widehat{w}(x)\psi(y)|\psi(x)-\psi(y)|}{|h|^{\beta}|x-y|^{n+2s}}\\
    &\geq -\frac{1}{8\Lambda}\frac{|\widehat{w}_{+}(x)-\widehat{w}_{+}(y)|^{2}}{|x-y|^{n+2s}}\psi^{2}(y)-16\Lambda^{2}\frac{|\widehat{w}_{+}(x)|^{2}|\psi(x)-\psi(y)|^{2}}{|x-y|^{n+2s}}\\
    &\geq -\frac{1}{4\Lambda}\frac{|(\widehat{w}_{+}\psi)(x)-(\widehat{w}_{+}\psi)(y)|^{2}}{|x-y|^{n+2s}}-32\Lambda^{2}\frac{|\widehat{w}_{+}(x)|^{2}|\psi(x)-\psi(y)|^{2}}{|x-y|^{n+2s}}.
\end{align*}

\item In case of $\frac{\delta_{h}w}{|h|^{\beta}}(x)\geq k\geq \frac{\delta_{h}w}{|h|^{\beta}}(y)$: Due to \eqref{pt : assmp.phi}, we estimate
\begin{align*}
    J_{1}&=\left({\Phi}\left(\frac{w_{h}(x)-w_{h}(y)}{|x-y|^{s}}\right)-{\Phi}\left(\frac{w(x)-w(y)}{|x-y|^{s}}\right)\right)\frac{\widehat{w}_{+}(x)\psi^{2}(x)}{|h|^{\beta}|x-y|^{n+s}}\\
    &\geq \frac{1}{\Lambda|h|^{\beta}}\frac{(\delta_{h}w(x)-\delta_{h}w(y))\widehat{w}_{+}(x)\psi^{2}(x)}{|x-y|^{n+2s}}\\
    &= \frac{1}{\Lambda}\frac{(\widehat{w}_{+}(x)+\widehat{w}_{-}(y))\widehat{w}_{+}(x)\psi^{2}(x)}{|x-y|^{n+2s}}.
\end{align*}
\end{enumerate}
Since $\widehat{w}_+(y)=0$ in this case, we have 
\begin{align*}
    J_{1}&\geq \frac{1}{\Lambda}\frac{(\widehat{w}_{+}(x)-\widehat{w}_{+}(y))^{2}\psi^{2}(x)}{|x-y|^{n+2s}}+\frac{1}{\Lambda}\frac{\widehat{w}_{-}(y)\widehat{w}_{+}(x)\psi^{2}(x)}{|x-y|^{n+2s}}\\
    &\geq \frac{1}{2\Lambda}\frac{|(\widehat{w}_{+}\psi)(x)-(\widehat{w}_{+}\psi)(y)|^{2}}{|x-y|^{n+2s}}-\frac{1}{\Lambda}\frac{|\widehat{w}_{+}(y)|^{2}|\psi(y)-\psi(x)|^{2}}{|x-y|^{n+2s}}\\
    &\quad+\frac{1}{\Lambda}\frac{\widehat{w}_{-}(y)\widehat{w}_{+}(x)\psi^{2}(x)}{|x-y|^{n+2s}}.
\end{align*}
Combining all the estimates, we have 
\begin{align*}
    J_{1}&\geq \frac{1}{4\Lambda}\frac{|(\widehat{w}_{+}\psi)(x)-(\widehat{w}_{+}\psi)(y)|^{2}}{|x-y|^{n+2s}}+\frac{1}{\Lambda}\frac{\widehat{w}_{-}(y)\widehat{w}_{+}(x)\psi^{2}(x)}{|x-y|^{n+2s}}\\
    &\quad-32\Lambda^{2}\frac{(|\widehat{w}_{+}(x)|^{2}+|\widehat{w}_{+}(y)|^{2})|\psi(x)-\psi(y)|^{2}}{|x-y|^{n+2s}}.
\end{align*}

\textbf{Estimate of $J_{2}$.} We also divide this part into two cases.
\begin{enumerate}
    \item In case of $\delta_{h}w(x)>\delta_{h}w(y)$: With \eqref{pt : assmp.phi}, we observe 
    \begin{align*}
        J_{2}&=\frac{1}{|h|^{\beta}}\left({\Phi}\left(\frac{w_{h}(x)-w_{h}(y)}{|x-y|^{s}}\right)-{\Phi}\left(\frac{w(x)-w(y)}{|x-y|^{s}}\right)\right)\frac{\widehat{w}_{+}(x)\psi^{2}(x)}{|x-y|^{n+s}}\\
        &\geq \frac{1}{|h|^{\beta}}\left(\dfrac{w_h(x)-w_h(y)-w(x)+w(y)}{|x-y|^s}\right)\frac{\widehat{w}_{+}(x)\psi^{2}(x)}{|x-y|^{n+s}}\\
        &\geq 0.
    \end{align*}
    \item In case of $\delta_{h}w(x)\leq\delta_{h}w(y)$: We estimate $J_{2}$ as
    \begin{align*}
        J_{2}&\geq \frac{\Lambda}{|h|^{\beta}}\left(\delta_{h}w(x)-\delta_{h}w(y)\right)\frac{\widehat{w}_{+}(x)\psi^{2}(x)}{|x-y|^{n+2s}}\\
        &={\Lambda}\left(\frac{\delta_{h}w(x)}{{|h|^{\beta}}}-k+\left(k-\frac{\delta_{h}w(y)}{{|h|^{\beta}}}\right)\right)\frac{\widehat{w}_{+}(x)\psi^{2}(x)}{|x-y|^{n+2s}}\\
        &\geq {\Lambda}\left(k-\frac{\delta_{h}w(y)}{{|h|^{\beta}}}\right)\frac{\widehat{w}_{+}(x)\psi^{2}(x)}{|x-y|^{n+2s}},
        \end{align*}
    where we have used \eqref{pt : assmp.phi} and the fact that $\frac{\delta_{h}w(x)}{{|h|^{\beta}}}-k\geq0$ if $\widehat{w}(x)\geq 0$.
    Therefore, we get
    \begin{align*}
        J_{2}\geq -\Lambda\widehat{w}_{+}(y)\frac{\widehat{w}_{+}(x)\psi^{2}(x)}{|x-y|^{n+2s}}.
    \end{align*}
\end{enumerate}
Similarly to $J_{2}$, we estimate $J_{3}$ as 
\begin{align*}
     J_{3}\geq -\Lambda\widehat{w}_{+}(x)\frac{\widehat{w}_{+}(y)\psi^{2}(y)}{|x-y|^{n+2s}}.
\end{align*}
Plugging all the above estimates of $J_{1},J_{2}$ and $J_{3}$ into \eqref{pt : lem.entildew.eq1}, we arrive at
\begin{align*}
    &(1-s)\int_{B_{\rho}(x_{1})}\int_{B_{\rho}(x_{1})}\frac{|(\widehat{w}_{+}\psi)(x)-(\widehat{w}_{+}\psi)(y)|^{2}}{|x-y|^{n+2s}}\,dx\,dy\\
    &\quad+(1-s)\int_{B_{\rho}(x_{1})}\int_{B_{\rho}(x_{1})}\frac{\widehat{w}_{-}(y)\widehat{w}_{+}(x)\psi^{2}(x)}{|x-y|^{n+2s}}\,dx\,dy\\
    &\leq c(1-s)\int_{B_{\rho}(x_{1})}\int_{B_{\rho}(x_{1})}\frac{|\widehat{w}_{+}(x)|^{2}|\psi(x)-\psi(y)|^{2}}{|x-y|^{n+2s}}\,dx\,dy\\
    &\quad+c(1-s)\int_{\bbR^{n}\setminus B_{\rho}(x_{1})}\int_{B_{\rho}(x_{1})}\widehat{w}_{+}(y)\frac{\widehat{w}_{+}(x)\psi^{2}(x)}{|x-y|^{n+2s}}\,dx\,dy 
    +c\int_{B_{\rho}(x_{1})}|\widetilde{f}_{\beta}\widehat{w}_{+}\psi^{2}|\,dx
\end{align*}
for some $c=c(n,\Lambda)$. After a few simple calculations together with the fact that 
\begin{equation*}
    |y-x|\geq \frac{\rho-r}{2\rho}|y-x_{1}|\quad\text{for any }x\in B_{\frac{\rho+r}{2}}(x_{1})\text{ and }y\in\bbR^{n}\setminus B_{\rho}(x_{1})
\end{equation*}
and 
\begin{equation*}
    |\psi(x)-\psi(y)|\leq \frac{2\rho}{\rho-r}|x-y|,
\end{equation*}
we obtain \eqref{pt : ineq.cca.diffq} for $\widehat{w}_{+}$. The proof of \eqref{pt : ineq.cca.diffq} in case of $\widehat{w}_{-}$ is similar.
\end{proof}

Next, we obtain more refined energy estimates by analyzing the last term appearing in \eqref{pt : ineq.cca.diffq}.

\begin{lemma}
\label{pt : lem.cca.refined}
Let $w\in W^{s,2}(B_{R}(x_{0}))\cap L^{1}_{2s}(\bbR^{n})$ be a weak solution to 
\begin{equation*}
    \mathcal{L}w=f\quad\text{in }B_{R}(x_{0}),
\end{equation*}  
where $f\in L^{\infty}\left(B_{2R}(x_{0})\right)$. For all $x_{1}\in B_{{R}/{2}}(x_0)$, $0<r<\rho\leq{R}/{4}$ and all $\beta \in (0,1]$, we have 
\begin{equation}
\label{el : ineq.ccarefined}
\begin{aligned}
    &(1-s)\int_{B_{r}(x_{1})}\int_{B_{r}(x_{1})}\frac{|(\widetilde{w}_{\beta}-k)_{\pm}(x)-(\widetilde{w}_{\beta}-k)_{\pm}(y)|^{2}}{|x-y|^{n+2s}}\,dx\,dy\\
    &\quad+(1-s)\int_{B_{r}(x_{1})}\int_{B_{\rho}(x_{1})}(\widetilde{w}_{\beta}-k)_{\pm}(x)\frac{(\widetilde{w}_{\beta}-k)_{\mp}(y)}{|x-y|^{n+2s}}\,dy\,dx\\
    &\leq \frac{c\rho^{2-2s}}{(\rho-r)^{2}}\int_{B_{\rho}(x_{1})}(\widetilde{w}_{\beta}-k)_{\pm}^{2}\,dx+ c{\|\widetilde{f}_{\beta}\|^{2}_{L^{\infty}(B_{R}(x_{0}))}}R^{2s}\left|A_{\pm}\left(\widetilde{w}_{\beta},x_{1},\rho;k\right)\right|\\
    &\quad+ c(1-s)\left(\frac{\rho}{\rho-r}\right)^{n+2s}\int_{\bbR^{n}\setminus B_{\rho}(x_{1})}\int_{B_{\rho}(x_{1})}(\widetilde{w}_{\beta}-k)_{\pm}(x)\frac{(\widetilde{w}_{\beta}-k)_{\pm}(y)}{|y-x_{1}|^{n+2s}}\,dx\,dy
\end{aligned}
\end{equation}
for some constant $c=c(n,\Lambda)$.
\end{lemma}
\begin{proof}
Using \eqref{pt : ineq.cca.diffq}, we have 
\begin{align*}
&(1-s)\int_{B_{\rho}(x_{1})}\int_{B_{\rho}(x_{1})}\frac{|((\widetilde{w}_{\beta}-k)_{\pm}\psi)(x)-((\widetilde{w}_{\beta}-k)_{\pm}\psi)(y)|^{2}}{|x-y|^{n+2s}}\,dy\,dx\\
    &\quad+(1-s)\int_{B_{\rho}(x_{1})}(\widetilde{w}_{\beta}-k)_{\pm}\psi^{2}\left(x\right)\left(\int_{B_{\rho}(x_{1})}\frac{(\widetilde{w}_{\beta}-k)_{\mp}\left(y\right)}{|x-y|^{n+2s}}\,dy\right)\,dx\\
    &\leq \frac{c\rho^{2-2s}}{(\rho-r)^{2}}\int_{B_{\rho}(x_{1})}(\widetilde{w}_{\beta}-k)_{\pm}^{2}\,dx\\
    &\quad +\frac{c(1-s)\rho^{n+2s}}{\left(\rho-r\right)^{n+2s}}\int_{\bbR^{n}\setminus B_{\rho}(x_{1})}\int_{B_{\rho}(x_{1})}(\widetilde{w}_{\beta}-k)_{\pm}(y)\frac{(\widetilde{w}_{\beta}-k)_{\pm}\psi^{2}(x)}{|y-x_{1}|^{n+2s}}\,dx\,dy\\
    &\quad+c\int_{B_{\rho}(x_{1})}|\widetilde{f}_{\beta}(\widetilde{w}_{\beta}-k)_{\pm}\psi^{2}|\,dx
\end{align*}
with $c=c(n,\Lambda)$, where $\psi(x)\in C_{c}^{\infty}\left(B_{\frac{\rho+r}{2}}(x_{1})\right)$ is a cut off function with $\psi\equiv 1$ on $B_{\rho}(x_{1})$.
We first use H\"older's inequality and Cauchy's inequality to see that
\begin{align*}
    &\int_{B_{\rho}(x_{1})}|\widetilde{f}_{\beta}(\widetilde{w}_{\beta})_{+}\psi^{2}|\,dx\\
    &\quad\leq \left(\int_{B_{\rho}(x_{1})}|\widetilde{f}_{\beta}\bfchi_{\{(\widetilde{w}_{\beta}-k)_{\pm}>0\}}|^{2}\,dx\right)^{\frac{1}{2}}\left(\int_{B_{\rho}(x_{1})}|(\widetilde{w}_{\beta}-k)_{\pm}\psi|^{2}\,dx\right)^{\frac{1}{2}}\\
    &\quad\leq c{\|\widetilde{f}_{\beta}\|^{2}_{L^{\infty}(B_{R}(x_{0}))}}\rho^{2s}|A_{\pm}(\widetilde{w}_{\beta},x_{1},\rho;k)|+{c}{\rho^{-2s}}\int_{B_{\rho}(x_{1})}|(\widetilde{w}_{\beta}-k)_{\pm}\psi|^{2}\,dx.
\end{align*}
Combining the above two inequalities along with the fact that $1\leq\rho/(\rho-r)$ and $\psi\equiv 1$ on $B_{\rho}(x_{1})$, we obtain the desired result.
\end{proof}

Using the above lemma together with Lemma \ref{el : lem.degio}, we obtain the following lemma which is a main tool in order to apply a bootstrap argument.
\begin{lemma}
\label{pt : lem.hol.ind}
Let $w\in W^{s,2}(B_{R}(x_{0}))\cap L^{1}_{2s}(\bbR^{n})$ be a weak solution to 
\begin{equation*}
    \mathcal{L}w=f\quad\text{in }B_{R}(x_{0}),
\end{equation*}  
where $f\in L^{\infty}(B_{2R}(x_{0}))$. Then there are constants $\gamma_{1}=\gamma_{1}(n,s_0,\Lambda)\in(0,1)$ and $c=c(n,s_0,\Lambda)$ such that for any $\beta \in (0,1]$, $\frac{\delta_{h}{w}}{|h|^{\beta}}\in C_{\mathrm{loc}}^{0,\gamma_{1}}(B_{R}(x_{0}))$ with the estimate
\begin{align*}
    \left\|\frac{\delta_{h}{w}}{|h|^{\beta}}\right\|_{L^{\infty}(B_{R/2}(x_{0}))}+R^{\gamma_{1}}\left[\frac{\delta_{h}{w}}{|h|^{\beta}}\right]_{C^{0,\gamma_{1}}(B_{R/2}(x_{0}))}&\leq c\widetilde{E}\left(\frac{\delta_{h}{w}}{|h|^{\beta}};B_{R}(x_{0})\right)\\
    &\quad+c\left\|\frac{\delta_{h}{f}}{|h|^{\beta}}\right\|_{L^{\infty}(B_{R}(x_{0}))},
\end{align*}
where the constant $s_0$ is determined in \eqref{el : choi.s0.sec3}.
\end{lemma}
\begin{proof}
By Lemma \ref{pt : lem.cca.refined}, we observe that $\frac{\delta_{h}w}{|h|^{\beta}}$ satisfies \eqref{el : ineq.degiorene} with $M=M(n,\Lambda)>0$ and $F=\left\|\frac{\delta_{h}{f}}{|h|^{\beta}}\right\|_{L^{\infty}(B_{R}(x_{0}))}$ in Lemma \ref{el : lem.degio}. Thus we obtain the desired result.
\end{proof}

Before giving the main result, we introduce a technical lemma from \cite[Lemma 5.6]{CafCab95} and a scaling-invariant property of non-homogeneous nonlocal equations.
\begin{lemma}
\label{pt : lem.tech.st}
Let three constants $0<\beta<1$, $0<\gamma\leq 1$ and $K>0$ be given. Assume that $g\in L^{\infty}([-1,1])$ with $\|g\|_{L^{\infty}([-1,1])}\leq K$. If 
\begin{align*}
\widetilde{g}_{\beta}\coloneqq\frac{\delta_{h}g}{|h|^{\beta}}\in C^{\gamma}([-1/2,1/2])
\end{align*}
with $\|\widetilde{g}_{\beta}\|_{C^{\gamma}([-1/2,1/2])}\leq K$ for any $0<|h|<{1}/{10}$, then we have the following: 
\begin{enumerate}
    \item if $\beta+\gamma<1$, $g\in C^{\beta+\gamma}([-1/2,1/2])$ with $\|g\|_{C^{\beta+\gamma}([-1/2,1/2])}\leq cK$, and
    \item if $\beta+\gamma>1$, $g\in C^{0,1}([-1/2,1/2])$ with $\|g\|_{C^{0,1}([-1/2,1/2])}\leq cK$
\end{enumerate}
for some constant $c=c(\beta,\gamma)$. Moreover, if $\widetilde{g}_{1}\in C^{\gamma}([-1/2,1/2])$ holds, then $g\in C^{1,\gamma}([-1/2,1/2])$ with the estimate $\|g\|_{C^{1+\gamma}([-1/2,1/2])}\leq cK$.
\end{lemma}

In the following straightforward lemma, we discuss the scaling properties the nonlocal equations we study.
\begin{lemma}
\label{el : lem.scale}
Let $f\in L^{\infty}(B_{R}(x_{0}))$ be given.
Let $w\in W^{s,2}(B_{R}(x_{0}))\cap L^{1}_{2s}(\bbR^{n})$ be a weak solution to 
\begin{equation*}
    \mathcal{L}w=f\quad\text{in }B_{R}(x_{0}).
\end{equation*}
Then $w_{R}(x)={w(R x+x_{0})}/{R^{s}}\in W^{s,2}(B_{1})\cap L^{1}_{2s}(\bbR^{n})$ is a weak solution to 
\begin{equation*}
    \mathcal{L}w_{R}=f_{R}\quad\text{in }B_{1},
\end{equation*}
where $f_{R}=R^{s}f(R x+x_{0})\in L^{\infty}(B_{1})$.
\end{lemma}

We now use a bootstrap argument to obtain H\"older continuity of the gradient of weak solutions.
\begin{lemma}
\label{el : lem.lochol}
Let $u\in W^{s,2}(B_{R}(x_{0}))\cap L^{1}_{2s}(\bbR^{n})$ be a weak solution to 
\begin{equation*}
    \mathcal{L}u=0\quad\text{in }B_{R}(x_{0}).
\end{equation*}
Then we have
\begin{equation*}
    \|\nabla u\|_{L^{\infty}(B_{R/2}(x_{0}))}+R^{\alpha_{0}}[\nabla u]_{C^{0,\alpha_{0}}(B_{R/2}(x_{0}))}\leq c\widetilde{E}(u/R;B_{R}(x_{0}))
\end{equation*}
for some constants $\alpha_0=\alpha_0(n,s_0,\Lambda)\in(0,1)$ and $c=c(n,s_0,\Lambda)$, where the constant $s_0$ is determined in \eqref{el : choi.s0.sec3}.
\end{lemma}
\begin{proof}
By Lemma \ref{el : lem.scale}, we may assume $x_{0}=0$ and $R=1$.
In light of Lemma \ref{pt : lem.hol.basis}, we have $u\in C_{\mathrm{loc}}^{\alpha_{0}}(B_{1})$ with 
\begin{equation}
\label{pt : defn.gamma0}
\alpha_{0}=\min\{\gamma_{0},\gamma_{1}\},
\end{equation}
where the constants $\gamma_{0}$ and $\gamma_{1}$ are determined in Lemma \ref{pt : lem.hol.basis} and Lemma \ref{pt : lem.hol.ind}, respectively.
In addition, if $B_{R}(x_{1})\subset B_{1}$, then we have
\begin{align}
\label{pt : thm.ineq.bas}
    R^{\alpha_{0}}[u]_{C^{0,\alpha_{0}}(B_{R/2}(x_{1}))}\leq c\widetilde{E}(u;B_{R}(x_{1}))
\end{align}
for some constant $c=c(n,s_0,\Lambda)$ Since $\gamma_{0}$ and $\gamma_{1}$ depend only on $n,s_0$ and $\Lambda$, there is a positive integer $i_{0}$ which depends only on $n,s_0$ and $\Lambda$  such that 
\begin{equation}
\label{el : choi.i0}
i_{0}\alpha_{0}\leq 1<(i_{0}+1)\alpha_{0}.
\end{equation}
We now fix $B_{10\rho}(x_{1})\subset B_{1}$ with $x_{1}$. Take a cut off function $\xi\in C_{c}^{\infty}\left(B_{4\rho}(x_{1})\right)$ such that $\xi\equiv 1$ on $B_{3\rho}(x_{1})$. By Lemma \ref{el : lem.loc}, we see that $u\xi\in C^{0,\alpha_{0}}(B_{5\rho}(x_{1}))$ with $u\xi\equiv 0$ on $\bbR^{n}\setminus B_{4\rho}(x_{1})$ is a weak solution to
\begin{equation*}
    \mathcal{L}(u\xi)=f\quad\text{in }B_{\rho}(x_{1}),
\end{equation*}
where 
{\small\begin{align*}
   f(x)= (1-s)\int_{\bbR^{n}\setminus B_{3\rho}(x_{1})}2\left[{{\Phi}}\left(\frac{(u\xi)(x)-(u\xi)(y)}{|x-y|^{s}}\right)-{{{\Phi}}\left(\frac{u(x)-u(y)}{|x-y|^{s}}\right)}\right]\frac{\,dy}{|x-y|^{n+s}}
\end{align*}}\\
\noindent
is in $C^{0,\alpha_{0}}\left(B_{2\rho}(x_{1})\right)$. Then $w(x)=\frac{(u\xi)(\rho x+x_{1})}{\rho^{s}}\in W^{s,2}(B_{5})\cap L^{1}_{2s}(\bbR^{n})\cap C^{0,\alpha_{0}}(B_{5})$ is a weak solution to
\begin{equation*}
    \mathcal{L}w=f_{\rho}\quad\text{in }B_{1},
\end{equation*}
where $f_{\rho}(x)=\rho^{s}f(\rho x+x_{1})\in C^{0,\alpha_{0}}\left(B_{2}\right)$ by Lemma \ref{el : lem.scale}.
In particular, we note from \eqref{el : lem.loc.est} that 
\begin{equation}
\label{pt : ineq.size.holf}
\begin{aligned}
    \sup_{x\in B_{1},h\in B_{\frac{1}{10}}}\frac{|f_{\rho}(x+h)-f_{\rho}(x)|}{|h|^{\alpha_{0}}}&=\sup_{x\in B_{\rho}(x_{1}),h\in B_{\frac{1}{10}}}\rho^{s}\frac{|f(x+\rho h)-f(x)|}{|h|^{\alpha_{0}}}\\
    &\leq c\rho^{\alpha_{0}-s}[u]_{C^{0,\alpha_{0}}(B_{3\rho}(x_{1}))}\\
    &\quad+ c\rho^{-s}\left(\|u\|_{L^{\infty}(B_{3\rho}(x_{1}))}+\mathrm{Tail}(u;B_{3\rho}(x_{1}))\right)
\end{aligned}
\end{equation}
for some constant $c=c(n,s_0,\Lambda)$.
By Lemma \ref{pt : lem.hol.ind}, we have $\frac{\delta_{h}{w}}{|h|^{\alpha_{0}}}\in C^{0,\alpha_{0}}_{\mathrm{loc}}(B_{1})$ with the estimate
\begin{align*}
   \left\|\frac{\delta_{h}{w}}{|h|^{\alpha_{0}}} \right\|_{C^{0,\alpha_{0}}(B_{1/2})}&\leq c\widetilde{E}\left(\frac{\delta_{h}{w}}{|h|^{\alpha_{0}}};B_{1}\right)+c\left\|\frac{\delta_{h}f_{\rho}}{|h|^{\alpha_{0}}}\right\|_{L^{\infty}(B_{1})}.
\end{align*}
By \eqref{pt : ineq.size.holf} along with Lemma \ref{pt : lem.hol.basis}, we get 
\begin{align*}
    \sup_{0<|h|<1/10}\left\|\frac{\delta_{h}{w}}{|h|^{\alpha_{0}}} \right\|_{C^{0,\alpha_{0}}(B_{1/2})}\leq c\rho^{-s}\widetilde{E}(u;B_{10\rho}(x_{1}))
\end{align*}
for some constant $c=c(n,s_0,\Lambda)$, where we have used the fact that
\begin{align*}
    \widetilde{E}\left(\frac{\delta_{h}{w}}{|h|^{\alpha_{0}}};B_{1}\right)  &\leq c\rho^{\alpha_0-s}[u\xi]_{C^{0,\alpha_{0}}(B_{5\rho}(x_{1}))}\leq c\rho^{-s}\widetilde{E}(u;B_{10\rho}(x_{1}))
\end{align*}
for some constant $c=c(n,s_0,\Lambda)$.
Thus by Lemma \ref{pt : lem.tech.st} along with the choice of $\alpha_{0}$ given in \eqref{pt : defn.gamma0}, we get $w\in C^{2\alpha_{0}}_{\mathrm{loc}}(B_{1})$ with the estimate 
\begin{align*}
    \|w\|_{C^{2\alpha_{0}}(B_{1/2})}&\leq c\rho^{-s}\widetilde{E}(u;B_{10\rho}(x_{1}))
\end{align*}
for some constant $c=c(n,s,\Lambda)$. Thus we have 
\begin{equation}
\label{pt : estimate1}
\begin{aligned}
    \left[u\right]_{C^{0,\alpha_{0}}(B_{\rho/2}(x_{1}))}&\leq c\rho^{-2\alpha_{0}}\widetilde{E}(u;B_{10\rho}(x_{1}))
\end{aligned}
\end{equation}
for some constant $c=c(n,s_0,\Lambda)$, where we have used  change of variables. Since the ball $B_{10\rho}(x_{1})$ is chosen arbitrarily, using standard covering arguments, we obtain that $u\in C^{2\alpha_{0}}_{\mathrm{loc}}(B_{1})$ with the estimate 
\begin{align*}
    [u]_{C^{2\alpha_{0}}(B_{\rho/2}(x_{1}))}\leq {c}{\rho^{-2\alpha_{0}}}\widetilde{E}(u;B_{\rho}(x_{1})),
\end{align*}
provided that $B_{\rho}(x_{1})\subset B_{1}$.
By following the above arguments with $\alpha_{0}$ replaced by $2\alpha_{0}$, we obtain 
\begin{align*}
    [u]_{C^{3\alpha_{0}}(B_{\rho/2}(x_{1}))}\leq {c}{\rho^{-3\alpha_{0}}}\widetilde{E}(u;B_{\rho}(x_{1})).
\end{align*}
By proceeding $i_{0}$ times, we obtain 
\begin{align*}
    \sup_{0<|h|<1/10}\left\|\frac{\delta_{h}{w}}{|h|^{i_{0}\alpha_{0}}} \right\|_{C^{0,\alpha_{0}}(B_{1/2})}\leq c\rho^{-s}\widetilde{E}(u;B_{10\rho}(x_{1}))
\end{align*}
for some constant $c=c(n,s_0,\Lambda)$, since the positive integer $i_{0}$ also depends only on $n,s_0$ and $\Lambda$.
Using Lemma \ref{pt : lem.tech.st} together with standard covering arguments, we obtain that $u\in C^{0,1}_{\mathrm{loc}}(B_{1})$ with the estimate 
\begin{align*}
    [u]_{C^{0,1}(B_{\rho/2}(x_{1}))}\leq {c}{\rho}^{-1}\widetilde{E}(u;B_{\rho}(x_{1}))
\end{align*}
for some constant $c=c(n,s_0,\Lambda)$, provided that $B_{\rho}(x_{1})\subset B_{1}$. We now repeat the above argument with $\beta=1$ to conclude that $u\in C^{1,\alpha_{0}}_{\mathrm{loc}}(B_{1})$ with the desired estimate
\begin{align*}
    \|\nabla u\|_{L^{\infty}(B_{\rho/2}(x_{1}))}+\rho^{-\alpha_{0}}[\nabla u]_{C^{0,\alpha_{0}}(B_{\rho/2}(x_{1}))}\leq {c}{\rho}^{-1}\widetilde{E}(u;B_{\rho}(x_{1}))
\end{align*}
for some constant $c=c(n,s_0,\Lambda)$, whenever $B_{\rho}(x_{1})\subset B_{1}$.
\end{proof}

Finally, we now employ Lemma \ref{el : lem.lochol} to prove our main result in the homogeneous case given by Theorem \ref{el : thm.main1}.

\textbf{Proof of Theorem \ref{el : thm.main1}.}
Let us fix $s_{0}\in(0,s]$ and $B_{R}(x_{0})\subset \Omega$. Then we observe that $u-(u)_{B_{R}(x_{0})}$ satisfies
\begin{equation*}
    \mathcal{L}(u-(u)_{B_{R}(x_{0})})=0\quad\text{in }B_{R}(x_{0}).
\end{equation*}
By Lemma \ref{el : lem.lochol} and \eqref{el : obs.e}, there is a constant $\alpha=\alpha(n,s_0,\Lambda)\in(0,1)$ such that
\begin{equation*}
    \|\nabla u\|_{L^{\infty}(B_{R/2}(x_{0}))}+R^{\alpha}[\nabla u]_{C^{0,\alpha}(B_{R/2}(x_{0}))}\leq c{E}(u/R;B_{R}(x_{0}))
\end{equation*}
for some constant $c=c(n,s_0,\Lambda)$, which completes the proof.
\qed

We conclude this section by proving an oscillation decay estimate for the first-order quotients of solutions.


\begin{lemma}
\label{pt : lem.hol.var}
Let $x_0 \in \mathbb{R}^n$, $R>0$ and suppose that $v \in W^{s,2}(B_{R}(x_0)) \cap L^1_{2s}(\mathbb{R}^n)$ is a weak solution to 
\begin{equation*}
    \setLcal v= 0\quad\text{in } B_{R}(x_{0}).
\end{equation*}
Then there are constants $\alpha_{0}=\alpha_{0}(n,s,\Lambda)\in(0,1)$ and $c=c(n,s,\Lambda)$ such that
\begin{align}
 \label{el : ineq.deosc}
		&\left\|\frac{\delta_{h}{v}}{|h|}-\left(\frac{\delta_{h}v}{|h|}\right)_{B_{R}(x_0)}\right\|_{L^{\infty}(B_{R/2}(x_{0}))}+R^{\alpha_0}\left [\frac{\delta_{h}{v}}{|h|}\right]_{C^{0,\alpha_0}(B_{R/2}(x_{0}))}\nonumber\\
		& \leq c E\left(\frac{\delta_{h}{v}}{|h|};B_{R}(x_{0})\right),
	\end{align}
 where $0<|h|<R/8$. In addition, for any fixed $s_0\in(0,1)$, the constants $\alpha_0$ and $c$ depend only on $n,s_0$ and $\Lambda$ when $s\geq s_0$.
\end{lemma}
\begin{proof}
Let us fix $s_0\in(0,s]$.
    	We first note that for any $h \in B_{R/8}\setminus\{0\}$, any $k \in \mathbb{R}$ and $j\in\{1,2,\ldots,n\}$, estimates \eqref{pt : ineq.cca.diffq} clearly also holds with $k$ replaced by $\left(\frac{\delta_{h}v}{|h|}\right)_{B_{R}(x_0)}+k$ and with $f=0$. This implies that \eqref{pt : ineq.cca.diffq} holds with $\widetilde{w}_{\beta}$ replaced by $\frac{\delta_{h}{v}}{|h|}-\left(\frac{\delta_{h}v}{|h|}\right)_{B_{R}(x_0)}$. Therefore Lemma \ref{pt : lem.hol.ind} along with \eqref{el : obs.e} yields the estimate \eqref{el : ineq.deosc} by taking  $\gamma_{1}=\alpha_0$, where the constant $\alpha_{0}$ is determined in Lemma \ref{el : lem.lochol}.
\end{proof}



\section{Oscillation decay and higher differentiability of the gradient}\label{sec:4}
The aim of this section is twofold. First of all, we prove gradient oscillation decay estimates on solutions to homogeneous nonlinear nonlocal equations which are consistent with given complement data. Moreover, we prove higher differentiability of the gradient of solutions to nonlinear nonlocal measure data problems and in particular Theorem \ref{HD}.

In order to obtain results that are stable as $s \to 1$, for the rest of this paper we fix some parameter $s_0 \in (0,1)$ and assume that
\begin{equation}
	\label{el : choi.s0.sec4}
	s\in [s_0,1).
\end{equation}

We start with the following observation.
\begin{remark}
\label{el : rmk.str}
In this remark, we will show that our nonlocal operator defined in \eqref{nonlocalop} can be rewritten as 
\begin{equation}
\label{el : eq,non}
    \mathcal{L}u(x)=(1-s)\mathrm{P}.\mathrm{V}.\int_{\RRn} \left({u(x)-u(y)}\right)K(x,y)\,dy,
\end{equation}
where the associated kernel $K:\bbR^{n}\times \bbR^{n}\to\bbR$ is measurable and satisfies 
\begin{equation} \label{Kcond}
    \frac{\Lambda^{-1}}{|x-y|^{n+2s}}\leq K(x,y) \leq\frac{\Lambda}{|x-y|^{n+2s}}\quad\text{for any }x,y\in\bbR^{n}, x\neq y.
\end{equation}
In fact, define
\begin{equation*}
    K(x,y):= \Phi\left(\frac{u(x)-u(y)}{|x-y|^{s}}\right) (u(x)-u(y))^{-1} |x-y|^{-n-s}\quad \forall x,y\in\bbR^{n},\,x\neq y.
\end{equation*}

In light of \eqref{assump}, we obtain that the kernel $K$ satisfies the conditions \eqref{Kcond} assumed in \cite{KuuMinSir15}, enabling us to apply the estimates obtained in \cite{KuuMinSir15}. 
\end{remark}

Using the previous remark, we next provide the following comparison estimate.
\begin{lemma}
\label{el : lem.comp1}
Let $\mu\in L^{\infty}(B_{R}(x_{0}))$ and let $u\in W^{s,2}(B_{R}(x_{0}))\cap L^{1}_{2s}(\bbR^{n})$ be a weak solution to 
\begin{equation*}
    \mathcal{L}u=\mu\quad\text{in }B_{R}(x_{0}).
\end{equation*}
Then there is a weak solution $v\in W^{s,2}(B_{R}(x_{0}))\cap L^{1}_{2s}(\bbR^{n})$ to
\begin{equation}
\label{el : eq.comp1}
\left\{
\begin{alignedat}{3}
\setLcal v&= 0&&\quad \mbox{in  $B_{R}(x_{0})$}, \\
v&=u&&\quad  \mbox{a.e. in }\RRn\setminus B_{R}(x_{0})
\end{alignedat} \right.
\end{equation}
such that
\begin{equation*}
    \dashint_{B_{R}(x_{0})}|u-v|\,dx\leq cR^{2s-n}|\mu|(B_{R}(x_{0}))
\end{equation*}
for some constant $c=c(n,\Lambda)$.
\end{lemma}
\begin{proof}
We first observe that the existence of the weak solution $v$ to \eqref{el : eq.comp1} can be proved in the same way as in e.g.\ \cite[Remark 3]{existence} or \cite[Appendix A]{Kyeongbae2}.

We next note from Remark \ref{el : rmk.str} that we are able to employ \cite[Lemma 3.2]{KuuMinSir15} with $q=1$, $p=2$ and $h=2s-1$. Thus we get
\begin{equation}\label{el : ineq.comp1}
\begin{aligned}
    &(2-2s)\int_{B_{2R}(x_{0})}\dashint_{B_{2R}(x_{0})}\frac{|(u-v)(x)-(u-v)(y)|}{|x-y|^{n+2s-1}}\,dx\,dy\\
    &\leq c\frac{(1-h)^{2}}{(1-s)(s-h)}R^{-s/2+2s-n}{|\mu|(B_{R}(x_{0}))}\\
    &\leq cR^{-s/2+2s-n}{|\mu|(B_{R}(x_{0}))}
\end{aligned}
\end{equation}
for some constant $c=c(n,\Lambda)$ independent of $s$, as $s>1/2$, which can be observed by following the proof of \cite[Lemma 3.2]{KuuMinSir15}. With the help of the Sobolev-Poincar\'e inequality given in \cite[Lemma 4.7]{CozziJFA}, we deduce
\begin{align*}
    \dashint_{B_{R}(x_{0})}|u-v|\,dx\leq c\dashint_{B_{2R}(x_{0})}|u-v|\,dx&\leq cR^{2s-n}{|\mu|(B_{R}(x_{0}))}\\
    &\leq cR^{2s-n}|\mu|(B_{R}(x_{0}))
\end{align*}
for some constant $c=c(n)$ independent of $s$, as $s>1/2$.
\end{proof}

We next provide a local higher Sobolev regularity estimate of weak solutions to a homogeneous problem with regular boundary data, which will be the essential ingredient to obtain suitable decay estimates at the gradient level. Note that in view of Lemma \ref{el : lem.graesti0} below and the localization argument from Lemma \ref{el : lem.loc}, we can always assume that any weak solution $u$ to \eqref{pt : eq.main} with $\mu\in L^{\infty}_{\mathrm{loc}}$ is in $W^{1,1}(\bbR^{n})$, as this assumption can always be removed in the end, see the arguments nearby \eqref{el : eq.hd}.
\begin{lemma}
\label{el : lem.bdry}
Let $\mu\in L^{\infty}(B_{R}(x_{0}))$, $q\in[1,\infty)$ and let $u\in W^{s,2}(B_{2R}(x_{0}))\cap W^{1,1}(\bbR^{n})$. Suppose $v\in W^{s,2}(B_{R}(x_{0}))\cap L^{1}_{2s}(\bbR^{n})$ is the weak solution to
\begin{equation*}
\left\{
\begin{alignedat}{3}
\setLcal v&= 0&&\quad \mbox{in  $B_{R}(x_{0})$}, \\
v&=u&&\quad  \mbox{a.e. in }\RRn\setminus B_{R}(x_{0}).
\end{alignedat} \right.
\end{equation*}
Let ${R}/{2}\leq r<\rho\leq {3R}/{4}$. Then there is a constant $\kappa=\kappa(n,s_0,\Lambda)\in(0,1)$ which is independent of $q$ such that  
\begin{align*}
     r^{-n/q+\kappa}[\nabla v]_{W^{\kappa,q}\left(B_{r}(x_{0})\right)}&\leq \frac{cR^{n+2s}}{(\rho-r)^{n+2s}}\left[E^{q}_{\loc}\left({\nabla v};B_{\rho}(x_{0})\right)+E\left({\nabla u};B_{\rho}(x_{0})\right)\right]\\
     &\quad+\frac{cR^{n+2s}}{(\rho-r)^{n+2s}}\frac{{|\mu|(B_{R}(x_{0}))}}{R^{n-2s+1}}
\end{align*}
holds with $c=c(n,s_0,\Lambda,q)$, where the constant $s_0$ is determined in \eqref{el : choi.s0.sec4}.
\end{lemma}

\begin{proof}
By Lemma \ref{el : lem.scale}, we may assume $x_{0}=0$, $R=1$ and ${1}/{2}\leq r<\rho\leq{3}/{4}$. Let us fix $q\in[1,\infty)$. By Lemma \ref{el : lem.comp1}, we have 
\begin{align}
\label{el : ineq1.bdry}
    \int_{B_{1}}|u-v|\,dx\leq c|\mu|(B_{1})
\end{align}
for some constant $c=c(n,\Lambda)$.
Fix $|h|\leq\frac{\rho-r}{10^{5}}\leq 1$ and choose $\beta=\beta(s_0)\in(0,1)$ such that
\begin{align}
\label{el : choi.beta}
    2s-\beta>2s_0\beta-1>0.
\end{align}
The choice of $s_0$ guarantees the existence of such $\beta$, where $s_{0}$ is determined in \eqref{el : choi.s0.sec4}. 
Then by Lemma \ref{el : lem.besi}, there is a covering $\{B_{|h|^{\beta}}(z_{i})\}_{i\in I}$ of $B_r$ such that $z_i\in B_r$, 
\begin{equation}
\label{el : ineq2.bdry1}
|I||h|^{n\beta}\leq c
\end{equation}
and
\begin{align}
\label{el : cov.ineq1}
\sup_{x \in \mathbb{R}^n} \sum_{i\in I}\bfchi_{B_{2^{k}|h|^{\beta}}(z_{i})} (x) \leq c2^{nk}
\end{align}
for some constant $c=c(n)$, where we denote $|I|$ the number of elements in the index set $I$. We now fix a positive integer $m_{0}$ such that 
\begin{equation}
\label{el : choi.m01}
\frac{1}{100}\left({\rho-r}\right)\leq2^{m_{0}+4}|h|^{\beta}<\frac{1}{50}\left(\rho-r\right)
\end{equation}
to see that 
\begin{equation}
\label{el : cov.ineq2}
2^{-2sm_{0}}|h|^{-1}\leq \frac{c}{(\rho-r)^{2s}}|h|^{2s\beta-1}\underset{\eqref{el : choi.beta}}{\leq} \frac{c}{(\rho-r)^{2s}},
\end{equation}
and 
\begin{equation}
\label{el : cov.ineq3}
    B_{2^{m_{0}+4}|h|^{\beta}}(z_{i})\subset B_{\frac{\rho-r}{50}}(z_{i})\subset B_{\frac{\rho-r}{10}}(z_{i}) \subset B_{\frac{r+\rho}{2}},
\end{equation}
since $z_i\in B_r$. Now for some constants $\alpha_0=\alpha_0(n,s_0,\Lambda)\in(0,1)$ and $c=c(n,s_0,\Lambda)$,
\begin{align}\label{el : cov.ineq4}
\begin{split}
    \dashint_{B_{|h|^{\beta}}(z_{i})}|\delta_{h}^{2}v|^{q}\,dx&= |h|^{q}\dashint_{B_{|h|^{\beta}}(z_{i})}\left|\delta_{h}\left(\frac{\delta_{h}v}{|h|}\right)\right|^{q}\,dx\\
&\leq c|h|^{q+q\alpha_{0}}\left[\frac{\delta_{h}v}{|h|}\right]^{q}_{C^{0,\alpha_{0}}\left(B_{2|h|^{\beta}}(z_{i})\right)}\\
&\leq c|h|^{q+q\alpha_{0}(1-\beta)}\left[{E}\left(\frac{\delta_{h}v}{|h|};B_{4|h|^{\beta}}(z_{i})\right)\right]^{q}
\end{split}
\end{align}
holds, where for the second line we have used the fact that 
\begin{equation*}
\dashint_{B_{|h|^{\beta}}(z_{i})}|\delta_{h}g|^q\,dx\leq \dashint_{B_{|h|^{\beta}}(z_{i})}|g(x+h)-g(x)|^q\,dx\leq |h|^{q\alpha_{0}}[g]^q_{C^{0,\alpha_{0}}(B_{2|h|^{\beta}}(z_{i}))}
\end{equation*}
for any $g\in C^{0,\alpha_{0}}(B_{2|h|^{\beta}}(z_{i}))$, and for the third line we have employed \eqref{el : ineq.deosc}. With the aid of Lemma \ref{el : lem.tail}, we further estimate the last term given in \eqref{el : cov.ineq4} as
\begin{align}\label{el : cov.ineq5}
\begin{split}
&{E}\left(\frac{\delta_{h}v}{|h|};B_{4|h|^{\beta}}(z_{i})\right)\\
&\quad\leq c\sum_{j=0}^{m_{0}+2}2^{-2sj}{E}_{\loc}\left(\frac{\delta_{h}v}{|h|};B_{2^{j+2}|h|^{\beta}}(z_{i})\right)\\
&\quad\quad+c2^{-2sm_{0}}\mathrm{Tail}\left(\frac{\delta_{h}v}{|h|}-\left(\frac{\delta_{h}v}{|h|}\right)_{B_{2^{m_{0}+4}|h|^{\beta}}(z_{i})};B_{2^{m_{0}+4}|h|^{\beta}}(z_i)\right):=I_1+I_2,
\end{split}
\end{align}
where $c=c(n)$ as $s>1/2$. For $I_1$, \eqref{eq:osc.est1} in Lemma \ref{el : thm.main4} yields
\begin{align}\label{el : cov.ineq6}
I_1\leq c\sum_{j=0}^{m_{0}+1}2^{-2sj}{E}_{\loc}\left(\nabla v;B_{2^{j+4}|h|^{\beta}}(z_{i})\right).
\end{align}
For $I_2$, using \cite[Lemma 2.3]{BLS}, \eqref{el : cov.ineq3}, $z_i\in B_r$ and \eqref{el : choi.m01}, we observe 
\begin{align*}
&\mathrm{Tail}\left(\frac{\delta_{h}v}{|h|}-\left(\frac{\delta_{h}v}{|h|}\right)_{B_{2^{m_{0}+4}|h|^{\beta}}(z_{i})};B_{2^{m_{0}+4}|h|^{\beta}}(z_i)\right)\\
&\quad \leq (2^{m_0+5}|h|^{\beta})^{2s}\left(\dfrac{3r+\rho}{3r+\rho-4|z_i|}\right)^{n+2s}\mathrm{Tail}\left(\frac{\delta_{h}v}{|h|}-\left(\frac{\delta_{h}v}{|h|}\right)_{B_{\frac{3r+\rho}{4}}};B_{\frac{3r+\rho}{4}}\right)\\
&\quad\quad+(2^{m_{0}+4}|h|^{\beta})^{-n}E_{\mathrm{loc}}\left(\frac{\delta_{h}v}{|h|};B_{\frac{3r+\rho}{4}}\right)\\
&\quad \leq c \left(\dfrac{1}{\rho-r}\right)^{n+2s}\mathrm{Tail}\left(\frac{\delta_{h}v}{|h|}-\left(\frac{\delta_{h}v}{|h|}\right)_{B_{\frac{3r+\rho}{4}}};B_{\frac{3r+\rho}{4}}\right)\\
&\quad\quad+c\left(\dfrac{1}{\rho-r}\right)^{n}E_{\mathrm{loc}}\left(\frac{\delta_{h}v}{|h|};B_{\frac{3r+\rho}{4}}\right)
\end{align*}
for some constant $c=c(n)$ independent of $s$, as $s>1/2$. Then \eqref{eq:osc.est1} in Lemma \ref{el : thm.main4} enable us to find that
\begin{align}\label{el : cov.ineq7}
\begin{split}
I_2&\leq\frac{c}{(\rho-r)^{n}}2^{-2sm_{0}}{E}_{\loc}\left(\nabla v;B_{\frac{7r+\rho}{8}}\right)\\
&\quad+\frac{c}{(\rho-r)^{n+2s}}2^{-2sm_{0}}\mathrm{Tail}\left(\frac{\delta_{h}v}{|h|}-\left(\frac{\delta_{h}v}{|h|}\right)_{B_{\frac{3r+\rho}{4}}};B_{\frac{3r+\rho}{4}}\right)
\end{split}
\end{align}
for some constant $c=c(n)$. For the last term in the above inequality, we note from \eqref{el : ineq1.bdry}, \eqref{el : choi.m01} and \eqref{el : cov.ineq2} that
\begin{equation}
\label{el : ineq2.bdry}
\begin{aligned}
   &2^{-2sm_{0}} \mathrm{Tail}\left(\frac{\delta_{h}v}{|h|}-\left(\frac{\delta_{h}v}{|h|}\right)_{B_{\frac{3r+\rho}{4}}};B_{\frac{3r+\rho}{4}}\right)\\
   &\quad\leq c2^{-2sm_{0}}(2^{m_{0}}|h|^{\beta})^{2s}|h|^{-1}\int_{B_{1}}|u-v|\,dx\\
   &\quad\quad+c2^{-2sm_{0}}\mathrm{Tail}\left(\frac{\delta_{h}u}{|h|}-\left(\frac{\delta_{h}u}{|h|}\right)_{B_{\frac{3r+\rho}{4}}};B_{\frac{3r+\rho}{4}}\right)\\
   &\quad\leq c|\mu|(B_{1})+c2^{-2sm_{0}}\mathrm{Tail}\left(\frac{\delta_{h}u}{|h|}-\left(\frac{\delta_{h}u}{|h|}\right)_{B_{\frac{3r+\rho}{4}}};B_{\frac{3r+\rho}{4}}\right)
\end{aligned}
\end{equation}
for some constant $c=c(n)$. Combining \eqref{el : cov.ineq5}--\eqref{el : ineq2.bdry} along with the fact that $(a+b+c+d)^{q}\leq 4^{q}(a^{q}+b^{q}+c^{q}+d^{q})$ for $a,b,c,d\geq 0$, we get 
\begin{align*}
    \left[{E}\left(\frac{\delta_{h}v}{|h|};B_{4|h|^{\beta}}(z_{i})\right)\right]^{q}&\leq c\left[\sum_{j=0}^{m_{0}+3}2^{-2sj}E_{\mathrm{loc}}(\nabla v;B_{{2^{j+2}|h|^{\beta}}}(z_{i}))\right]^{q}\\
    &\quad +\frac{c}{(\rho-r)^{qn}}2^{-2sm_{0}q}\left[E_{\mathrm{loc}}\left(\nabla v;B_{\frac{7r+\rho}{8}}\right)\right]^{q}\\
    &\quad +\frac{c}{(\rho-r)^{q(n+2s)}}\left[|\mu|(B_{1})\right]^{q}\\
    &\quad +\frac{c}{(\rho-r)^{q(n+2s)}}\left[\mathrm{Tail}\left(\frac{\delta_{h}u}{|h|}-\left(\frac{\delta_{h}u}{|h|}\right)_{B_{\frac{3r+\rho}{4}}};B_{\frac{3r+\rho}{4}}\right)\right]^{q}\\
    &\eqqcolon L_{1}(i)+L_{2}+L_{3}+L_{4}
\end{align*}
for each $i\in I$. As a result, together with \eqref{el : cov.ineq4} we have
\begin{align*}
    \sum_{i\in I}\dashint_{B_{|h|^{\beta}}(z_{i})}|\delta_{h}^{2}v|^{q}\,dx&\leq \sum_{i\in I}c|h|^{q+q\alpha_{0}(1-\beta)}\left(L_{1}(i)+L_{2}+L_{3}+L_{4}\right)\\
    &=: \widetilde{L}_1(i)+\widetilde{L}_2+\widetilde{L}_3+\widetilde{L}_4.
\end{align*}
We now estimate each term $\widetilde{L}_1(i)$, $\widetilde{L}_2$, $\widetilde{L}_3$, and $\widetilde{L}_4$.

\textbf{Estimate of $\widetilde{L}_1(i)$.}
We first note from H\"older's inequality that
\begin{equation*}
    \sum_{j}a_{j}b_{j}\leq \left(\sum_{j}a_{j}b_{j}^{q}\right)^{\frac{1}{q}}\left(\sum_{j}a_{j}\right)^{\frac{1}{q'}}
\end{equation*}
for any $a_{j},b_{j}\geq0$. Using this along with Fubini's theorem, Lemma \ref{el : lem.ave}, \eqref{el : cov.ineq1} and \eqref{el : cov.ineq3}, we estimate $\widetilde{L}_1(i)$ as 
\begin{align*}
    \widetilde{L}_1(i)&\leq c|h|^{q+q\alpha_{0}(1-\beta)}\sum_{i\in I}\left(\sum_{j=0}^{m_{0}+1}2^{-2sj}\dashint_{B_{2^{j+4}|h|^{\beta}}(z_{i})}\left|\nabla v-\left(\nabla v\right)_{B_{2^{j+4}|h|^{\beta}}(z_{i})}\right|^{q}\,dx\right)\\
    &\qquad\quad\qquad\qquad\qquad\times\underbrace{\left(\sum_{j=0}^{m_{0}+2}2^{-2sj}\right)^{\frac{q}{q'}}}_{\leq c(q)}\\
    &\leq c|h|^{q+q\alpha_{0}(1-\beta)}\sum_{j=0}^{m_{0}+2}2^{-2sj}\sum_{i\in I}\dashint_{B_{2^{j+2}|h|^{\beta}}(z_{i})}\left|\nabla v-\left(\nabla v\right)_{B_{\frac{r+\rho}{2}}}\right|^{q}\,dx\\
    &\leq c|h|^{q+q\alpha_{0}(1-\beta)-n\beta}\sum_{j=0}^{m_{0}}2^{-2sj}\dashint_{B_{\frac{r+\rho}{2}}}\left|\nabla v-\left(\nabla v\right)_{B_{\frac{r+\rho}{4}}}\right|^{q}\,dx\\
    &\leq c|h|^{q+q\alpha_{0}(1-\beta)-n\beta}\left[E^{q}_{\mathrm{loc}}\left(\nabla v;B_{\frac{r+\rho}{2}}\right)\right]^{q}
\end{align*}
for some constant $c=c(n,\Lambda,q)$. 

\textbf{Estimate of $\widetilde{L}_2$.}
We use \eqref{el : ineq2.bdry1} to see that
\begin{align*}
   \widetilde{L}_2&\leq \frac{c}{(\rho-r)^{qn}}|h|^{q+q\alpha_{0}(1-\beta)-n\beta}\left[E^{q}_{\mathrm{loc}}\left(\nabla v;B_{\frac{7r+\rho}{8}}\right)\right]^{q}
\end{align*}
for some constant $c=c(n,\Lambda)$.

\textbf{Estimate of $\widetilde{L}_3$.}
By \eqref{el : ineq2.bdry1}, we have 
\begin{equation*}
    \widetilde{L}_3\leq \frac{c}{(\rho-r)^{q(n+2s)}}|h|^{q+q\alpha_{0}(1-\beta)-n\beta}\left[|\mu|(B_{1})\right]^{q}
\end{equation*}
for some constant $c=c(n,\Lambda)$.

\textbf{Estimate of $\widetilde{L}_4$.}
With the aid of Lemma \ref{el : thm.main4} and \eqref{el : ineq2.bdry1}, we have 
\begin{align*}
    \widetilde{L}_4&\leq \frac{c}{(\rho-r)^{q(n+2s)}}|h|^{q+q\alpha_{0}(1-\beta)}\sum_{i\in I}\left[\mathrm{Tail}\left(\frac{\delta_{h}u}{|h|}-\left(\frac{\delta_{h}u}{|h|}\right)_{B_{\frac{3r+\rho}{4}}};B_{\frac{3r+\rho}{4}}\right)\right]^{q}\\
    &\leq \frac{c}{(\rho-r)^{q(n+2s)}}|h|^{q+q\alpha_{0}(1-\beta)-n\beta}\left[E(\nabla u;B_{\rho})\right]^{q}.
\end{align*}
Combining all the estimates, we obtain
\begin{align*}
    \int_{B_{r}}|\delta_{h}^{2}v|^{q}\,dx&\leq\sum_{i\in I}\int_{B_{|h|^{\beta}}(z_{i})}|\delta_{h}^{2}v|^{q}\,dx\\
    &\leq \frac{c|h|^{q(1+\gamma)}}{(\rho-r)^{q(n+2s)}}\left[|\mu|(B_{1})\right]^{q}+\frac{c|h|^{q(1+\gamma)}}{(\rho-r)^{qn}}\left[E^{q}_{\mathrm{loc}}\left(\nabla v;B_{\frac{7r+\rho}{8}}\right)\right]^{q}\\
    &\quad+\frac{c|h|^{q(1+\gamma)}}{(\rho-r)^{q(n+2s)}}\left[E(\nabla u;B_{\rho})\right]^{q},
\end{align*}
by taking 
\begin{align*}
    \gamma=\gamma(n,s_0,\Lambda):=\min\left\{2s_0\beta-1,\alpha_{0}(1-\beta)\right\}.
\end{align*}
In light of Lemma \ref{el : lem.emb}, we obtain 
\begin{align*}
    [\nabla v]_{W^{\kappa,q}(B_{r})}&\leq \frac{1}{\gamma^{1+2/q}(1-\gamma)}\frac{c}{(\rho-r)^{n+2s}}\left[E^{q}_{\loc}\left({\nabla v};B_{\rho}\right)+E\left({\nabla u};B_{\rho}\right)\right]\\
    &\quad+\frac{1}{\gamma^{1+2/q}(1-\gamma)}\frac{c}{(\rho-r)^{n+2s}}{{|\mu|(B_{1})}}
\end{align*}
for some constant $c=c(n,\Lambda,q)$ by taking $\kappa=\kappa(n,s_0,\Lambda):=\frac{\gamma}{2}$ which is independent of $q$. Since $\gamma$ depends only on $n,s_0$ and $\Lambda$, it completes the proof.
\end{proof}

By combining Lemma \ref{el : lem.emb.hol} and Lemma \ref{el : lem.bdry} with a bootstrap argument, we now obtain a crucial decay estimate for gradients of solutions to homogeneous problems.
\begin{lemma}
\label{el : lem.holosc}
    Let $u\in W^{s,2}(B_{2R}(x_{0}))\cap W^{1,1}(\bbR^{n})$ and let $v\in W^{s,2}(B_{R}(x_{0}))\cap L^{1}_{2s}(\bbR^{n})$ be a weak solution to
\begin{equation*}
\left\{
\begin{alignedat}{3}
\setLcal v&= 0&&\quad \mbox{in  $B_{R}(x_{0})$}, \\
v&=u&&\quad  \mbox{a.e. in }\RRn\setminus B_{R}(x_{0}).
\end{alignedat} \right.
\end{equation*}
Then there is a constant $\alpha_{1}=\alpha_{1}(n,s_0,\Lambda)\in(0,1)$ such that 
\begin{align*}
    R^{\alpha_{1}}[\nabla v]_{C^{0,\alpha_{1}}\left(B_{{R}/{4}}(x_{0})\right)}&\leq c\left[E_{\loc}\left(\nabla v;B_{{R}/{2}}(x_{0})\right)+E\left(\nabla u;B_{{R}/{2}}(x_{0})\right)\right]\\
    &\quad+c\frac{|\mu|(B_{R}(x_{0}))}{R^{n-2s+1}},
\end{align*}
where $c=c(n,s_0,\Lambda)$. In particular, this implies that for any $\rho\in\left(0,{1}/{4}\right]$,
\begin{align*}
    \osc_{B_{\rho R}(x_{0})}\nabla v\leq c\rho^{\alpha_{1}}\left[E_{\loc}\left(\nabla v;B_{{R}/{2}}(x_{0})\right)+E\left(\nabla u;B_{{R}/{2}}(x_{0})\right)+\frac{|\mu|(B_{R}(x_{0}))}{R^{n-2s+1}}\right]
\end{align*}
holds for some constant $c=c(n,s_0,\Lambda)$, where the constant $s_0$ is determined in \eqref{el : choi.s0.sec4}.
\end{lemma}
\begin{proof}
We may assume $x_{0}=0$ and $R=1$ by Lemma \ref{el : lem.scale}. First, consider the smallest natural number $l_{0}=l_{0}(n,s,\Lambda)$ such that 
    \begin{equation}\label{el : alpha0}
        \kappa-n\frac{n-l_{0}\kappa}{n}=(l_{0}+1)\kappa-n>0,
    \end{equation}
    where the number $\frac{n}{n-l_{0}\kappa}$ is $l_{0}$-th $\kappa$-fractional Sobolev conjugate number of $1$ and the constant $\kappa=\kappa(n,s_0,\Lambda)\in(0,1)$ is determined in Lemma \ref{el : lem.bdry}. 
    We write $1_{k}^{*}=\frac{n}{n-k\kappa}$ and $R_{k}={1}/{2}-\frac{k}{4l_{0}}$ for any integer $k$ with $0\leq k\leq l_0$. We point out that every constant introduced until now depends only on $n,s_0$ and $\Lambda$.
    We first note from Lemma \ref{el : lem.emb.poi} and Lemma \ref{el : lem.bdry} that
    \begin{align*}
        E^{ 1_{k+1}^{*}}_{\loc}(\nabla v;B_{R_{k+1}})&\leq c[\nabla v]_{W^{\kappa,1_{k}^{*}}(B_{R_{k+1}})}\\
        &\leq c\left(E^{1_{k}^{*}}_{\loc}\left({\nabla v};B_{R_{k}}\right)+E\left({\nabla u};B_{R_{k}}\right)\right)+c{|\mu|(B_{R_{k}}})
    \end{align*}
    for any integer $k$ such that $0\leq k\leq l_0$ with some constant $c=c(n,s_0,\Lambda)$. As a result, from Lemma \ref{el : lem.bdry} again we have 
    \begin{align*}
        [\nabla v]_{W^{\kappa,1_{l_{0}}^{*}}(B_{\frac{1}{4}})}&\leq c\left(E_{\loc}\left({\nabla v};B_{1/2}\right)+c\sum_{k=0}^{l_{0}}E\left({\nabla u};B_{R_{k}}\right)\right)+c\sum_{k=0}^{l_{0}}{|\mu|(B_{R_{k}}})\\
        &\leq cE_{\loc}\left({\nabla v};B_{1/2}\right)+cE\left(\nabla u;B_{1/2}\right)+c{|\mu|}(B_{1})
    \end{align*}
    for some constant $c=c(n,s_0,\Lambda)$, as the constants $l_{0}$ and $\kappa$ depend only on $n,s_0$ and $\Lambda$. In light of Lemma \ref{el : lem.emb.hol} along with the choice of the constant $l_{0}$, we obtain 
    \begin{align*}
        [\nabla v]_{C^{0,\alpha_{1}}(B_{{1}/{4}})}\leq cE_{\loc}\left({\nabla v};B_{1/2}\right)+cE\left(\nabla u;B_{1/2}\right)+c{|\mu|}(B_{1})
    \end{align*}
    for some constant $c=c(n,s_0,\Lambda)$ with the choice $\alpha_{1}=\kappa(l_{0}+1)-n>0$ by \eqref{el : alpha0}. 
\end{proof}

By employing Lemma \ref{lem:reg.cacc} and Lemma \ref{pt : lem.hol.var}, we prove the following one which will be used to obtain a borderline regularity.
\begin{lemma}
\label{el : lem.twoqu}
Let us fix $h\in B_{1}\setminus\{0\}$ and $\beta\in(0,1)$. Let $v\in W^{s,2}(B_{4|h|^{\beta}}(x_{0}))\cap L^{1}_{2s}(\bbR^{n})$ be a weak solution to 
\begin{equation*}
    \mathcal{L}v=0\quad\text{in }B_{4|h|^{\beta}}(x_{0}).
\end{equation*}
 Then we have 
\begin{align*}
   \dashint_{B_{|h|^{\beta}}(x_{0})}{|\delta_{h}^{2}v|}\,dx\leq c|h|^{1+s(1-\beta)}E\left(\frac{\delta_{h}v}{|h|};B_{4|h|^{\beta}}(x_{0})\right)
\end{align*}
for some constant $c=c(n,\Lambda)$.
\end{lemma}
\begin{proof}
We may assume $x_{0}=0$ and fix $h\in B_{1}\setminus \{0\}$. We first note from \cite[Proposition 2.6]{BL} with $\alpha=s$, $p=1$, $\psi=\frac{\delta_{h}v}{|h|}-\left(\frac{\delta_{h}v}{|h|}\right)_{B_{2|h|^{\beta}}}$, $R=|h|^{\beta}$ and $h_{0}=|h|^{\beta}$,
\begin{equation}
\label{el : ineq1.lemtt}
\begin{aligned}
     \int_{B_{|h|^{\beta}}}{|\delta_{h}^{2}v|}\,dx&\leq|h|^{s}\sup_{0<\widetilde{h}<|h|}\int_{B_{|h|^{\beta}}}\frac{|\delta_{\widetilde{h}}(\delta_{h}v)|}{|\widetilde{h}|^{s}}\,dx\\
     &\leq |h|^{s+1}\sup_{0<\widetilde{h}<|h|^{\beta}}\int_{B_{|h|^{\beta}}}\left|\frac{\delta_{\widetilde{h}}}{|\widetilde{h}|^{s}}\left(\frac{\delta_{h}v}{|h|}-\left(\frac{\delta_{h}v}{|h|}\right)_{B_{2|h|^{\beta}}}\right)\right|\,dx\\
     &\leq |h|^{s+1-n\beta}\left(\sup_{0<\widetilde{h}<|h|^{\beta}}\dashint_{B_{|h|^{\beta}}}\left|\frac{\delta_{\widetilde{h}}}{|\widetilde{h}|^{s}}\left(\frac{\delta_{h}v}{|h|}-\left(\frac{\delta_{h}v}{|h|}\right)_{B_{2|h|^{\beta}}}\right)\right|^2\,dx\right)^{\frac12}\\
     &\leq c(1-s)^{\frac12}|h|^{s+1-n\beta/2}\left[\frac{\delta_{h}v}{|h|}\right]_{W^{s,2}(B_{2|h|^{\beta}})}\\
     &\quad+c|h|^{1+s(1-\beta)}E^2_{\mathrm{loc}}\left(\frac{\delta_{h}v}{|h|};B_{2|h|^{\beta}}\right)\eqqcolon I,
\end{aligned}
\end{equation}
where we have used the fact that $|h|\leq|h|^{\beta}$. We first note that for any $g\in L^{1}$, $x,y\in \bbR^{n}$ and $k\in\bbR$, we directly observe
\begin{equation*}
    |g(x)-g(y)|^{2}\leq 4|(g-k)_{+}(x)-(g-k)_{+}(y)|^{2}+4|(g-k)_{-}(x)-(g-k)_{-}(y)|^{2}.
\end{equation*}
Using this and Lemma \ref{lem:reg.cacc} with $x_{1}=0$, $r=2|h|^{\beta}$, $\rho=3|h|^{\beta}$, $f=0$, $\widetilde{w}_{\beta}=\frac{\delta_{h}v}{|h|}$ and $k=\left(\frac{\delta_{h}v}{|h|}\right)_{B_{3|h|^{\beta}}}$ to get 
\begin{align*}
    &(1-s)|h|^{-n\beta}\left[\frac{\delta_{h}v}{|h|}\right]^{2}_{W^{s,2}(B_{2|h|^{\beta}})}\\
    &\leq c|h|^{-2s\beta}\left[E_{\mathrm{loc}}^{2}\left(\frac{\delta_{h}v}{|h|};B_{3|h|^{\beta}}\right)\right]^{2}\\
    &\quad+c|h|^{-2s\beta}E_{\mathrm{loc}}\left(\frac{\delta_{h}v}{|h|};B_{3|h|^{\beta}}\right)\mathrm{Tail}\left(\frac{\delta_{h}v}{|h|}-\left(\frac{\delta_{h}v}{|h|}\right)_{B_{3|h|^\beta}};B_{3|h|^{\beta}}\right)
\end{align*}
for some constant $c=c(n,\Lambda)$. By employing this, Cauchy's inequality and Lemma \ref{pt : lem.hol.var}, we further estimate $I$ as 
\begin{equation}
\label{el : ineq2.lemtt}
\begin{aligned}
    I&\leq 
    c|h|^{1+s(1-\beta)}\left[E_{\mathrm{loc}}^{2}\left(\frac{\delta_{h}v}{|h|};B_{3|h|^{\beta}}\right) +{\mathrm{Tail}}\left(\frac{\delta_{h}v}{|h|}-\left(\frac{\delta_{h}v}{|h|}\right)_{B_{3|h|^\beta}};B_{3|h|^{\beta}}\right)\right]\\
    &\leq  c|h|^{1+s(1-\beta)}E\left(\frac{\delta_{h}v}{|h|};B_{4|h|^{\beta}}\right).
\end{aligned}
\end{equation}
Combining \eqref{el : ineq1.lemtt} and \eqref{el : ineq2.lemtt}, we obtain the desired result.
\end{proof}

We next provide a borderline regularity result for the gradient of solutions to measure data problems. We point out that obtaining a small increment of fractional differentiability can be achieved to the proof given in Lemma \ref{el : lem.bdry}. However, to reach the optimal amount of differentiability, we use Lemma \ref{el : lem.twoqu} and a bootstrap argument inspired by \cite{MinCZMD,AKMARMA}.
\begin{lemma}
\label{el : lem.graesti}
Let $\mu\in L^{\infty}(B_{R}(x_{0}))$ and let $u\in W^{s,2}(B_{2R}(x_{0}))\cap W^{1,1}(\bbR^{n})$ be a weak solution to
\begin{equation*}
    \mathcal{L}u=\mu\quad\text{in }B_{R}(x_{0}).
\end{equation*}
Then for any $\sigma\in(0,2s_{0}-
1)$, we have 
\begin{align*}
   R^{-n+\sigma} [\nabla u]_{W^{\sigma,1}\left(B_{{R}/{2}}(x_{0})\right)}\leq cE\left({\nabla u};B_{R}(x_{0})\right)+c\frac{{|\mu|(B_{R}(x_{0}))}}{R^{n-2s+1}}
\end{align*}
for some constant $c=c(n,s_{0},\Lambda,\sigma)$, where the constant $s_{0}$ is determined in \eqref{el : choi.s0.sec4}.
\end{lemma}

\begin{proof}
By Lemma \ref{el : lem.scale}, we may assume $x_{0}=0$ and $R=1$.
Fix $|h|\leq\frac{1}{100}$, $\sigma\in(0,2s_{0}-1)$, and let $\beta=\beta(s_{0},\sigma)\in(0,1)$ satisfying 
\begin{equation}
\label{el : choi.beta1}
2s_{0}\beta-1>\frac{\sigma+1}{2}>\sigma.
\end{equation}
By Lemma \ref{el : lem.besi}, there is a covering $\{B_{|h|^{\beta}}(z_{i})\}_{i\in I}$ of $B_{1/2}$, such that $z_{i}\in B_{1/2}$, $|I||h|^{n\beta}\leq c$ and
\begin{align}
\label{el : cov.ineq}
\sup_{x \in \mathbb{R}^n} \sum_{i\in I}\bfchi_{B_{2^{k}|h|^{\beta}}(z_{i})}(x) \leq c2^{nk}
\end{align}
for some constant $c=c(n)$, where we denote $|I|$ the number of elements in the index set $I$. We now fix a positive integer $m_{0}$ such that 
\begin{equation}
\label{el : choi.m0}
1/8\leq2^{m_{0}+4}|h|^{\beta}<1/4.
\end{equation}
By Lemma \ref{el : lem.comp1}, there is the weak solution $v_{i}\in W^{s,2}(B_{4|h|^{\beta}}(z_{i}))\cap L^{1}_{2s}(\bbR^{n})$ to
\begin{equation*}
\left\{
\begin{alignedat}{3}
\setLcal v_{i}&= 0&&\quad \mbox{in  $B_{4|h|^{\beta}}(z_{i})$}, \\
v_{i}&=u&&\quad  \mbox{a.e. in }\RRn\setminus B_{4|h|^{\beta}}(z_{i})
\end{alignedat} \right.
\end{equation*}
such that
\begin{equation}
\label{el : comp.ineq2}
    \dashint_{B_{4|h|^{\beta}}(z_{i})}|u-v_{i}|\,dx\leq c|h|^{(-n+2s)\beta}|\mu|(B_{4|h|^{\beta}}(z_{i}))
\end{equation}
for some constant $c=c(n,\Lambda)$.
We now note that
\begin{align}\label{el : comp.ineq3}
\begin{split}
\dashint_{B_{|h|^{\beta}}(z_{i})}|\delta_{h}^{2}u|\,dx&\leq\dashint_{B_{|h|^{\beta}}(z_{i})}|\delta_{h}^{2}(u-v_{i})|\,dx+\dashint_{B_{|h|^{\beta}}(z_{i})}|\delta_{h}^{2}v_i|\,dx\eqqcolon J_{1}+J_{2}.
\end{split}
\end{align}
By \eqref{el : comp.ineq2}, we have 
\begin{align}\label{el : comp.ineq4}
    J_{1}\leq c|h|^{(-n+2s)\beta}|\mu|(B_{4|h|^{\beta}}(z_{i}))
\end{align}
for some constant $c=c(n,\Lambda)$. We next estimate $J_{2}$ as 
\begin{align}\label{el : comp.ineq5}
\begin{split}
J_{2}&\leq c|h|^{s(1-\beta)+1}{E}\left(\frac{\delta_{h}v_{i}}{|h|};B_{4|h|^{\beta}}(z_{i})\right)\\
&\leq c|h|^{s(1-\beta)}\dashint_{B_{4|h|^{\beta}}(z_{i})}|u-v_{i}|\,dx+ c|h|^{s(1-\beta)+1}{E}\left(\frac{\delta_{h}u}{|h|};B_{4|h|^{\beta}}(z_{i})\right)\\
&\leq  c|h|^{s(1-\beta)
	+(2s-n)\beta}|\mu|(B_{4|h|^{\beta}}(z_{i}))+c|h|^{s(1-\beta)+1}{E}\left(\frac{\delta_{h}u}{|h|};B_{4|h|^{\beta}}(z_{i})\right)
\end{split}
\end{align}
for some constant $c=c(n,\Lambda)$, where we have used Lemma \ref{el : lem.twoqu} and \eqref{el : comp.ineq2}. We further estimate the last term given in the above inequality as follows: 
\begin{align}\label{el : comp.ineq6}
\begin{split}
{E}\left(\frac{\delta_{h}u}{|h|};B_{4|h|^{\beta}}(z_{i})\right)&\leq c\sum_{j=0}^{m_{0}+2}2^{-2sj}E_{\mathrm{loc}}\left(\frac{\delta_{h}u}{|h|};B_{2^{j+2}|h|^{\beta}}(z_{i})\right)\\
&+c2^{-2sm_{0}}\mathrm{Tail}\left(\frac{\delta_{h}u}{|h|}-\left(\frac{\delta_{h}u}{|h|}\right)_{B_{2^{m_{0}+4}|h|^{\beta}}(z_{i})};B_{2^{m_{0}+4}|h|^{\beta}}(z_{i})\right)\\
&\leq c\sum_{j=0}^{m_{0}+1}2^{-2sj}E_{\mathrm{loc}}(\nabla u;B_{2^{j+4}|h|^{\beta}}(z_{i}))+c2^{-2sm_{0}}E\left(\nabla u;B_{{3}/{4}}\right)
\end{split}
\end{align}
for some constant $c=c(n)$, where we have used Lemma \ref{el : lem.tail} and Lemma \ref{el : thm.main4}. As a result, merging \eqref{el : comp.ineq3}--\eqref{el : comp.ineq6} we have 
\begin{align*}
    \sum_{i\in I}\dashint_{B_{|h|^{\beta}}(z_{i})}|\delta_{h}^{2}u|\,dx&\leq \sum_{i\in I}c|h|^{(2s-n)\beta}|\mu|(B_{4|h|^{\beta}}(z_{i}))\\
    &\quad+c|h|^{s(1-\beta)+1}\sum_{i\in I}\sum_{j=0}^{m_{0}+1}2^{-2sj}E_{\mathrm{loc}}(\nabla u;B_{2^{j+4}|h|^{\beta}}(z_{i}))\\
    &\quad +c|h|^{s(1-\beta)+1}\sum_{i\in I}2^{-2sm_{0}}E\left(\nabla u;B_{{3}/{4}}\right)\eqqcolon L_{1}+L_{2}+L_{3}.
\end{align*}
We now estimate each term $L_{1},L_{2}$ and $L_{3}$.

\textbf{Estimate of $L_{1}$.} Using the fact that $|I||h|^{n\beta}\leq c$ for some constant $c=c(n)$, we get
\begin{align*}
    L_{1}\leq c|h|^{(2s-n)\beta}|\mu|(B_{1}),
\end{align*}
where $c=c(n,\Lambda)$.

\textbf{Estimate of $L_{2}$.}
We use Fubini's theorem, Lemma \ref{el : lem.ave} and \eqref{el : cov.ineq} to see that
\begin{equation}\label{el : comp.ineq7}
\begin{aligned}
L_{2}&\leq c|h|^{1+s(1-\beta)}\sum_{j=0}^{m_{0}}2^{-2sj}\sum_{i\in I}\dashint_{B_{2^{j+4}|h|^{\beta}}(z_{i})}\left|\nabla u-\left(\nabla u\right)_{B_{{3}/{4}}}\right|\,dx\\
&\leq c|h|^{1+s(1-\beta)-n\beta}\sum_{j=0}^{m_{0}}2^{-2sj}\int_{B_{{3}/{4}}}\left|\nabla u-\left(\nabla u\right)_{B_{{3}/{4}}}\right|\,dx\\
&\leq c|h|^{1+s(1-\beta)-n\beta}E_{\mathrm{loc}}(\nabla u;B_{3/4})
\end{aligned}
\end{equation}
for some constant $c=c(n,\Lambda)$ independent of $s$, as $s>1/2$.

\textbf{Estimate of $L_{3}$.}
With the aid of \eqref{el : choi.m0} and \eqref{el : cov.ineq} along with the fact that $|I||h|^{n\beta}\leq c$ for some constant $c=c(n)$, we have 
\begin{align*}
    L_{3}&\leq c|h|^{s(1-\beta)+1-n\beta}\sum_{i\in I}2^{-2sm_{0}}|h|^{n\beta}E\left(\nabla u;B_{{3}/{4}}\right)\\
    &\leq c|h|^{s(1-\beta)+1-n\beta}2^{-2sm_{0}}E\left(\nabla u;B_{{3}/{4}}\right)\\
    &\leq c|h|^{s(1-\beta)+1-n\beta+2s\beta}E\left(\nabla u;B_{{3}/{4}}\right).
\end{align*}
Combining all the estimates, we obtain
\begin{align}\label{el : comp.ineq8}
\begin{split}
\int_{B_{{1}/{2}}}|\delta_{h}^{2}u|\,dx&\leq\sum_{i\in I}\int_{B_{|h|^{\beta}}(z_{i})}|\delta_{h}^{2}u|\,dx\\
&\leq c|h|^{2s\beta}|\mu|(B_{1})+c|h|^{1+s(1-\beta)}E_{\mathrm{loc}}\left(\nabla u;B_{{3}/{4}}\right)\\
&\quad+c|h|^{1+2s\beta+s(1-\beta)}E\left(\nabla u;B_{{3}/{4}}\right)
\end{split}
\end{align}
for some constant $c=c(n,\Lambda)$. 
In light of Lemma \ref{el : lem.emb} along with the choice of 
\begin{equation*}
    \sigma_{0}\coloneqq{\min\left\{2s_{0}\beta-1,\frac{\sigma+1}{2}s(1-\beta)\right\}}
\end{equation*}
and a standard covering argument,
we obtain $\nabla u\in W^{\sigma_{0},1}_{\mathrm{loc}}(B_{1})$ with the estimate
\begin{align}
\label{el : comp.ineqsigma01}
[\nabla u]_{W^{\sigma_{0},1}\left(B_{{7}/{8}}\right)}\leq cE\left({\nabla u};B_1\right)+c{{|\mu|(B_{1})}}
\end{align}
with $c=c(n,s_{0},\Lambda,\sigma)$. By the choice of $\beta$ given in \eqref{el : choi.beta1}, we observe $\sigma_{0}=\frac{\sigma+1}{2}s(1-\beta)$.  We now use a bootstrap argument to increase the differentiability. To do this, we aim to obtain a more refined estimate of the term $L_{2}$ using the information  $\nabla u\in W^{\sigma_{0},1}_{\mathrm{loc}}(B_{1})$. To this end, we have
\begin{align}\label{el : comp.ineq9}
\begin{split}
L_{2}&\leq c|h|^{1+s(1-\beta)}\sum_{j=0}^{m_{0}}2^{-2sj}\sum_{i\in I}E_{\mathrm{loc}}(\nabla u;B_{2^{j+4}|h|^{\beta}}(z_i))\\
&\leq c|h|^{1+s(1-\beta)+\sigma_{0}\beta-n\beta}\sum_{j=0}^{m_{0}}2^{(-2s+\sigma_{0})j}\sum_{i\in I}[\nabla u]_{W^{\sigma_{0},1}(B_{2^{j+4}|h|^{\beta}}(z_{i}))}\\
&\leq c|h|^{1+s(1-\beta)+\sigma_{0}\beta-n\beta}\sum_{j=0}^{m_{0}}2^{(-2s+\sigma_{0})j}[\nabla u]_{W^{\sigma_{0},1}\left(B_{{7}/{8}}\right)}\\
&\leq c|h|^{1+s(1-\beta)+\sigma_{0}\beta-n\beta}\left[E(\nabla u;B_{1})+|\mu|(B_{1})\right]
\end{split}
\end{align}
for some constant $c=c(n,s_{0},\Lambda,\sigma)$, where we have used Poincar\'e inequality and \eqref{el : comp.ineqsigma01}. Plugging \eqref{el : comp.ineq9} into \eqref{el : comp.ineq8} instead of \eqref{el : comp.ineq7}, we obtain $\nabla u\in W_{\mathrm{loc}}^{\sigma_{1},1}(B_{1})$ with the estimate
\begin{equation*}
    [\nabla u]_{W^{\sigma_{1},1}(B_{7/8})}\leq cE(\nabla u;B_{1})+c|\mu|(B_{1}),
\end{equation*}
where 
\begin{equation*}
    \sigma_{1}\coloneqq{\min\left\{2s_{0}\beta-1,\sigma_{0}\beta+\frac{\sigma+1}{2}s(1-\beta)\right\}}
\end{equation*}
and $c=c(n,s_{0},\Lambda,\sigma)$. Similarly, we get $\nabla u\in W_{\mathrm{loc}}^{\sigma_{2},1}(B_{1})$ with the estimate
\begin{equation*}
    [\nabla u]_{W^{\sigma_{2},1}(B_{7/8})}\leq cE(\nabla u;B_{1})+c|\mu|(B_{1}),
\end{equation*}
where 
\begin{equation*}
    \sigma_{2}\coloneqq{\min\left\{2s_{0}\beta-1,\sigma_{1}\beta+\frac{\sigma+1}{2}s(1-\beta)\right\}}=\sum_{i=0}^{2}\beta^{i}\frac{\sigma+1}{2}s(1-\beta).
\end{equation*} 
Indeed, there is the smallest positive integer $N=N(n,s_0,\sigma)$ such that 
\begin{equation*}
    \sigma<\sum_{i=0}^{N}\beta^{i}\frac{\sigma+1}{2}s(1-\beta)<\frac{\sigma+1}{2}.
\end{equation*}
By proceeding $N-2$ times as in the above, we have 
$\nabla u\in W_{\mathrm{loc}}^{\sigma_{N},1}(B_{1})$ with the estimate
\begin{equation*}
    [\nabla u]_{W^{\sigma_{N},1}(B_{7/8})}\leq cE(\nabla u;B_{1})+c|\mu|(B_{1}),
\end{equation*}
where 
\begin{equation*}
    \sigma_{N}=\sum_{i=0}^{N}\beta^{i}\frac{\sigma+1}{2}s(1-\beta)
\end{equation*} 
and $c=c(n,s_{0},\Lambda,\sigma)$. As in \eqref{el : ineq2.embdi} with $q=1$, ${\gamma}=\sigma$ and $\gamma+\epsilon=\sigma_{N}$, we get 
\begin{equation*}
     [\nabla u]_{W^{\sigma,1}(B_{1/2})}\leq cE(\nabla u;B_{1})+c|\mu|(B_{1})
\end{equation*}
for some constant $c=c(n,s_{0},\Lambda,\sigma)$.
\end{proof}

Next, we prove that any weak solution to \eqref{pt : eq.main} with $\mu\in L^{\infty}_{\mathrm{loc}}$ belongs to $W^{1+t,1}_{\mathrm{loc}}$ for some small exponent $t$. 
\begin{lemma}
\label{el : lem.graesti0}
Let $\mu\in L^{\infty}(B_{R}(x_{0}))$ and let $u\in W^{s,2}(B_{R}(x_{0}))\cap L^{1}_{2s}(\bbR^{n})$ be a weak solution to
\begin{equation}
\label{el : eq.graesti0}
    \mathcal{L}u=\mu\quad\text{in }B_{R}(x_{0}).
\end{equation}
Then there is a constant $t=t(n,s_0,\Lambda)\in(0,1)$ such that
\begin{align*}
    R^{-n}\|\nabla u\|_{L^{1}(B_{R/2}(x_{0}))}+R^{-n+t}[\nabla u]_{W^{t,1}\left(B_{{R}/{2}}(x_{0})\right)}&\leq cE(u/R;B_{R}(x_{0}))\\
    &\quad+c\frac{{|\mu|(B_{R}(x_{0}))}}{R^{n-2s+1}}
\end{align*}
for some constant $c=c(n,s_0,\Lambda)$, where the constant $s_0$ is determined in \eqref{el : choi.s0.sec4}.
\end{lemma}
\begin{proof}
By Lemma \ref{el : lem.scale}, we may assume $x_0=0$ and $R=1$. We now choose $h\in B_{1/100}\setminus \{0\}$ and $\beta\in(0,1)$ such that $2s_0\beta-1>0$. Then we have a covering
$\{B_{|h|^{\beta}}(z_{i})\}_{i\in I}$ of $B_{1/2}$, such that $z_{i}\in B_{1/2}$, $|I||h|^{n\beta}\leq c$. By Lemma \ref{el : lem.graesti}, we have
\begin{align*}
     \dashint_{B_{|h|^{\beta}}(z_{i})}|\delta_{h}^{2}u|\,dx&\leq c|h|^{
    (2s-n)\beta}|\mu|(B_{4|h|^{\beta}}(z_{i}))\\
    &\quad+c|h|^{s(1-\beta)+1}{E}\left(\frac{\delta_{h}u}{|h|};B_{4|h|^{\beta}}(z_{i})\right)
\end{align*}
for some constant $c=c(n,s_0,\Lambda)$, since we did not use the assumption $u\in W^{1,1}(\bbR^{n})$ for the estimates of $J_{1}$ and $J_{2}$ given in Lemma \ref{el : lem.graesti}. We now follow the same iterative scheme given in \cite[Lemma 4.2]{KuuSimYan22} to obtain the following estimate
\begin{align*}
    \sup_{0<|h|<\frac{1}{100}}\int_{B_{\frac{1}{2}}}\frac{|\delta_{h}^{2}u|}{|h|^{1+t_{0}}}\,dx\leq c\widetilde{E}(u;B_{1})+c|\mu|(B_{1})
\end{align*}
for some constants $c=c(n,s_0,\Lambda)$ and $t_{0}=t_{0}(n,s_0,\Lambda)\in(0,1)$, since the number of iteration depends on the constant $\alpha_{0}$ which is determined in \eqref{el : thm.main1}. We recall the notation $\widetilde{E}(u;\cdot)$ which is defined in \eqref{el : defn.widee}.  Then, by \cite[Proposition 2.6]{BL} with $R=1$, $\alpha=2s-1$ and $p=1$, and \eqref{el : ineq.comp1} with $R=1$, we observe that
\begin{align*}
     \sup_{0<|h|<\frac{1}{100}}\int_{B_{{1}/{2}}}\frac{|\delta_{h}u|}{|h|^{2s-1}}\,dx&\leq c(1-s)[u]_{W^{2s-1,1}(B_{3/4})}+c\|u\|_{L^{1}(B_{3/4})}
\end{align*}
for some constant $c=c(n,s_0,\Lambda)$. We consider the weak solution $v\in W^{s,2}(B_1)\cap L^{1}_{2s}(\bbR^{n})$ to \eqref{el : eq.comp1} with $R=1$ and $x_0=0$. Since $v$ satisfies \eqref{el : ineq.degiorene} which follows from Lemma \ref{pt : lem.hol.basis}, a few simple calculations along with \eqref{el : res.difrac} and \eqref{el : ineq.comp1} with $R=1$ yield
\begin{align*}
    (1-s)[u]_{W^{2s-1,1}(B_{3/4})}&\leq (1-s)[u-v]_{W^{2s-,1}(B_{3/4})}+(1-s)[v]_{W^{2s-,1}(B_{3/4})}\\
    &\leq c|\mu|(B_{1})+c\widetilde{E}(v;B_{1})\\
    &\leq c|\mu|(B_{1})+c\widetilde{E}(u;B_{1})
\end{align*}
for some constant $c=c(n,\Lambda)$.
Thus, combining the above two inequalities, we obtain
\begin{align*}
     \sup_{0<|h|<\frac{1}{100}}\int_{B_{{1}/{2}}}\frac{|\delta_{h}u|}{|h|^{2s-1}}\,dx\leq c|\mu|(B_{1})+c\widetilde{E}(u;B_{1})
\end{align*}
for some constant $c=c(n,s_0,\Lambda)$. We now choose $t=\min\{t_0,2s-1\}$.
Using the fact that $u-(u)_{B_{R}(x_{0})}$ is also a weak solution to \eqref{el : eq.graesti0}, Lemma \ref{el : lem.emb1} and Lemma \ref{el : lem.emb}, we get the desired estimate. 
\end{proof}

We now employ the localization argument given in Lemma \ref{el : lem.loc}, Lemma \ref{el : lem.graesti} and Lemma \ref{el : lem.graesti0} to prove Theorem \ref{HD}.

\textbf{Proof of Theorem \ref{HD}.} Let us fix $s_{0}\in(0,s]$ and $B_{\rho}(x_{1})\subset\Omega$.
Let us take sequences $\{u_j\}_{j=1}^{\infty}\subset W^{s,2}(\bbR^{n})$,  $\{g_j\}_{j=1}^{\infty}\subset C_{0}^{\infty}(\bbR^{n})$ and $\{\mu_j\}_{j=1}^{\infty}\subset C_{0}^{\infty}(\bbR^{n})$ satisfying (1)-(7) given in Definition \ref{SOLA}. Then we first observe from Lemma \ref{el : lem.graesti0} along with a standard covering argument that $u_j\in W^{1,1}_{\mathrm{loc}}(\Omega)$ with the estimate
\begin{equation}
\label{el : ineq1.hd}
\begin{aligned}
    \rho^{-n}\|\nabla u_j\|_{L^{1}(B_{4\rho/5}(x_{1}))}+\rho^{-n+t}[\nabla u_j]_{W^{t,1}(B_{4\rho/5}(x_{1}))}&\leq cE(u_j/\rho;B_{\rho}(x_{1}))\\
    &\quad+c\frac{|\mu_j|(B_{\rho}(x_{1}))}{\rho^{n-2s+1}}
\end{aligned}
\end{equation}
for some constants $t=t(n,s_0,\Lambda)\in(0,1)$ and $c=c(n,s_0,\Lambda)$, whenever $B_{\rho}(x_{1})\subset \Omega$. By (5)-(7) from Definition \ref{SOLA}, we have that the right-hand side of the above inequality is bounded independently of $j$. Therefore, using a standard compactness argument based on e.g.\ \cite[Theorem 7.1]{DinPalVal12}, we observe that up to passing to a subsequence we have
\begin{equation*}
    \nabla u_j\to 
\nabla u\quad\text{in } L^{1}_{\mathrm{loc}}(\Omega) \text{ as } j\to\infty
\end{equation*}
and 
\begin{equation}
\label{el : conv.hd}
    \nabla u_j(x)\to \nabla u(x)\quad\text{a.e. }x\in\Omega\text{ as }j\to\infty.
\end{equation}
We now localize the equation
\begin{equation}
\label{el : eq.lochd}
    \mathcal{L}u_j=\mu_j\quad\text{in }\Omega.
\end{equation}
Let us take $\sigma$ and $q$ satisfying \eqref{el : choi.sigq} with $s=s_0$. We next choose $\widetilde{\sigma}\in(0,2s_0-1)$ satisfying
\begin{equation}
\label{el : choi.s0hd}
\widetilde{\sigma}-n\geq \sigma-n/q.
\end{equation}
This is always possible, because $2s_0-1-n>\sigma-n/q$.
Let us take a cut off function $\xi\in C_{c}^{\infty}\left(B_{{4\rho}/{5}}(x_{1})\right)$ with $\xi\equiv 1$ on $B_{{3\rho}/{5}}(x_{1})$. Then by Lemma \ref{el : lem.loc}, we have that $w_{j}\coloneqq u_{j}\xi\in W^{s,2}(B_{{4\rho}/{5}}(x_{1}))\cap W^{1,1}(\bbR^{n})$ is a weak solution to 
\begin{equation}
\label{el : eq.hd}
    \mathcal{L}w_j=\mu_{j}+f_{j}\quad\text{in }B_{{\rho}/{5}}(x_{1}), 
\end{equation}
where 
\begin{align*}
    f_j(x)&=2(1-s)\int_{\mathbb{R}^{n}\setminus B_{{3\rho}/{5}}(x_{1})}{\Phi}\left(\frac{w_j(x)-w_j(y)}{|x-y|^{s}}\right)\frac{\,dy}{|x-y|^{n+s}}\\
    &\quad-2(1-s)\int_{\mathbb{R}^{n}\setminus B_{{3\rho}/{5}}(x_{1})}{\Phi}\left(\frac{u_j(x)-u_j(y)}{|x-y|^{s}}\right)\frac{\,dy}{|x-y|^{n+s}}
\end{align*}
and $f_j\in L^{\infty}\left(B_{{2\rho}/{5}}(x_{1})\right)$ by \eqref{el : ineq.tail}.
By taking $\nu_j\coloneqq\mu_j+f_j$, we observe that $w_j\in W^{s,2}(B_{{4\rho}/{5}}(x_{1}))\cap W^{1,1}(\bbR^{n})$ is a weak solution to 
\begin{equation*}
    \mathcal{L}w_j=\nu_j\quad\text{in }B_{{\rho}/{5}}(x_{1}).
\end{equation*}
By Lemma \ref{el : lem.graesti}, we have 
\begin{align*}
   \rho^{-n+\widetilde{\sigma}} [\nabla w_{j}]_{W^{\widetilde{\sigma},1}(B_{\rho/10}(x_{1}))}\leq cE(\nabla w_j;B_{\rho/5}(x_{1}))+c\frac{|\nu_{j}|(B_{\rho/5}(x_{1}))}{\rho^{n-2s+1}}.
\end{align*}
By Lemma \ref{el : lem.embdiff} with $\gamma=\sigma$ and \eqref{el : choi.s0hd}, we have
\begin{align*}
  &\rho^{-n/q}\|\nabla w_{j}\|_{L^{q}(B_{\rho/10}(x_{1}))} +\rho^{-n/q+\sigma} [\nabla w_{j}]_{W^{\sigma,q}(B_{\rho/10}(x_{1}))}\\
  &\leq c\rho^{-n}\|\nabla w_j\|_{L^{1}(B_{\rho/10}(x_{1}))}+ cE(\nabla w_j;B_{\rho/5}(x_{1}))+c\frac{|\nu_{j}|(B_{\rho/5}(x_{1}))}{\rho^{n-2s+1}}.
\end{align*}
Using the fact that $\xi\equiv 0$ on $\mathbb{R}^{n}\setminus B_{4\rho/5}(x_{1})$ and \eqref{el : ineq.tail}, we note
\begin{equation}
\label{el : ineq3.hd}
\begin{aligned}
    E(\nabla w_j;B_{\rho/5}(x_{1}))\leq c\dashint_{B_{4\rho}/{5}(x_{1})}|\nabla w_j|\,dx&\leq \widetilde{E}(\nabla u_j;B_{4\rho/5}(x_{1}))\\
    &\quad+{c}\widetilde{E}(u_j/\rho;B_{4\rho/5}(x_{1}))
\end{aligned}
\end{equation}
and 
\begin{align*}
    {|\nu_{j}|(B_{\rho/5}(x_{1}))}&\leq {|\mu_{j}|(B_{\rho/5}(x_{1}))}+c\rho^{n}\|f_{j}\|_{(B_{\rho/5}(x_{1}))}\\
    &\leq {|\mu_{j}|(B_{\rho/5}(x_{1}))}+c\rho^{n-2s}\mathrm{Tail}(u_j;B_{3\rho/5}(x_{1}))
\end{align*}
for some constant $c=c(n,\Lambda)$.
After a few simple calculations along with the above two inequalities and the fact that $w_j=u_j$ on $B_{\rho/10}(x_1)$, we have 
\begin{align*}
    &\rho^{-n/q}\|\nabla u_{j}\|_{L^{q}(B_{\rho/10}(x_{1}))} +\rho^{-n/q+\sigma} [\nabla u_{j}]_{W^{\sigma,q}(B_{\rho/10}(x_{1}))}\\
    &\leq \widetilde{E}(\nabla u_j;B_{4\rho/5}(x_{1}))+{c}\widetilde{E}(u_j/\rho;B_{4\rho/5}(x_{1}))+c\frac{|\mu_{j}|(B_{\rho/5}(x_{1}))}{\rho^{n-2s+1}}
\end{align*}
for some constant $c=c(n,s_0,\Lambda,q)$.
We now use \eqref{el : ineq1.hd} and then employ the fact that $u_j-(u_j)_{B_{4R/5}(x_{0})}$ is also solution to \eqref{el : eq.lochd}, in order to deduce
\begin{align*}
    \rho^{-n/q}\|\nabla u_{j}\|_{L^{q}(B_{\rho/10}(x_{1}))} +\rho^{-n/q+\sigma} [\nabla u_{j}]_{W^{\sigma,q}(B_{\rho/10}(x_{1}))}&\leq cE(u_j/\rho;B_{\rho}(x_{1}))\\
    &\quad+c\frac{|\mu_{j}|(B_{\rho}(x_{1})}{\rho^{n-2s+1}}
\end{align*}
for some constant $c=c(n,s_0,\Lambda,q)$.
Using the dominated convergence theorem together with \eqref{el : conv.hd}, (5)-(7) given in Definition \ref{SOLA} and a standard covering argument, we obtain the desired result.
\qed

We end this section by giving the proof of Corollary \ref{el : cor.sin}.
\begin{proof}[Proof of Corollary \ref{el : cor.sin}.]
By Theorem \ref{HD}, we have $\nabla u\in W^{\sigma,1}_{\mathrm{loc}}(\Omega,R^{n})$ for any $\sigma<2s-1$. Thus the desired result directly follows from Lemma \ref{el : lem.sin}.
\end{proof}

\section{Potential estimates at the gradient level} \label{sec:5}
In this section, we prove pointwise estimates of the gradient of solutions to non-homogeneous nonlinear nonlocal equations with general measure data.

In order to upgrade our comparison estimate to the gradient level, we shall utilize the following classical iteration lemma, see e.g.\ \cite[Lemma 6.1]{GiustiBook}.
\begin{lemma}\label{iterlem}
	Let $h\colon [1/2,3/4]\to \mathbb{R}$ be a non-negative and bounded function, and let $a$ and $M$ be non-negative numbers. Assume that the inequality 
	$ 
	h(\delta_1)\le  (1/2) h(\delta_2)+(\delta_2-\delta_1)^{-M}a
	$
	holds, whenever $1/2\le \delta_1<\delta_2\le 3/4$. Then there is a constant $c=c(M)$ such that
    $
	h(1/2)\le ca.
	$ 
\end{lemma}

We next prove comparison estimates on the gradient level by employing an interpolation argument.
\begin{corollary} \label{cor:gradcomp} 
	Let $\mu \in L^\infty(B_{R}(x_{0}))$ and assume that $u \in W^{s,2}(B_{2R}(x_0)) \cap L^1_{2s}(\mathbb{R}^n)\cap W^{1,1}(\bbR^{n})$ is a weak solution of $\mathcal{L} u = \mu$ in $B_{R}(x_{0})$. Moreover, consider the weak solution $v \in W^{s,2}(B_{R}(x_0)) \cap L^1_{2s}(\mathbb{R}^n)$ of
	\begin{equation*}
		\begin{cases} \normalfont
			\mathcal{L} v = 0 & \text{ in } B_{R}(x_0), \\
			v = u & \text{ a.e. in } \mathbb{R}^n \setminus B_{R}(x_0).
		\end{cases}
	\end{equation*}
	Then we have the comparison estimate
	\begin{align*}
		& \dashint_{B_{R/2}(x_0)} |\nabla u-\nabla v| dx \\ & \leq  c \frac{|\mu|(B_R(x_0))}{R^{n-2s+1}}+\left( \frac{|\mu|(B_R(x_0))}{R^{n-2s+1}}\right)^{1-\theta}E(\nabla u;B_{R}(x_{0}))^{\theta},
	\end{align*}
	where $\theta \in (0,1)$ and $c$ depend only on $n,s_0$ and $\Lambda$. Here, the constant $s_0$ is determined in \eqref{el : choi.s0.sec4}.
\end{corollary}

\begin{proof}
	We prove the result for $R=1$ and $x_0=0$, since the general case can then be deduced by a straightforward scaling argument as in Lemma \ref{el : lem.scale}. Let $\kappa=\kappa(n,s_0,\Lambda) \in (0,1)$ be given by Lemma \ref{el : lem.bdry}. Let us fix ${1}/{2}\leq r<\rho\leq{3}/{4}$. 
	 By \cite{BrezisM}, we have the Gagliardo-Nirenberg-type inequality
	\begin{align*}
		\norm{\nabla u- \nabla v}_{L^1(B_{r})}& \leq c \norm{u-v}_{L^1(B_{r})}^{1-\theta}\left(\|u-v\|_{L^{1}(B_{r})}+ \|\nabla u-\nabla v\|_{L^{1}(B_{r})}\right)^{\theta}\\
  &\quad+c\norm{u-v}_{L^1(B_{r})}^{1-\theta}[\nabla u- \nabla v]_{W^{\kappa,1}(B_{r})}^\theta,
	\end{align*}
	where $c$ and $\theta:=\tfrac{1}{1+\kappa} \in (0,1)$ depend only on $n,s_0$ and $\Lambda$, as $\kappa$ depends only on the aforementioned parameters. By applying Young's inequality to the first term in the right-hand side of the above inequality, we have
   \begin{align*}
     \norm{\nabla u- \nabla v}_{L^1(B_{r})}\leq c \norm{u-v}_{L^1(B_{r})}+ c\norm{u-v}_{L^1(B_{r})}^{1-\theta}[\nabla u- \nabla v]_{W^{\kappa,1}(B_{r})}^\theta
   \end{align*}
	for some constant $c=c(n,\Lambda,s_0)$.
	Together with Lemma \ref{el : lem.bdry} and Lemma \ref{el : lem.graesti}, we obtain
	{\small\begin{align*}
		\norm{\nabla u- \nabla v}_{L^1(B_{r})} & \leq c|\mu|(B_1)+c |\mu|(B_1)^{1-\theta}\left ([\nabla u]_{W^{\kappa,1}(B_{r})} + [\nabla v]_{W^{\kappa,1}(B_{r})} \right )^\theta \\
		& \leq c|\mu|(B_1)^{1-\theta} \left (E\left(\nabla u;B_1\right)+|\mu|(B_{1}) \right )^\theta\\
  &\quad +\frac{c}{(\rho-r)^{n+2s}}|\mu|(B_1)^{1-\theta} \left (E_{\mathrm{loc}}(\nabla v;B_\rho)+E\left(\nabla u;B_\rho\right)+|\mu|(B_{1}) \right )^\theta\\
  & \leq c|\mu|(B_1)^{1-\theta} \left (E\left(\nabla u;B_1\right)+|\mu|(B_{1}) \right )^\theta\\
  &\quad +\frac{c}{(\rho-r)^{n+2s}}|\mu|(B_1)^{1-\theta} \left (E_{\mathrm{loc}}(\nabla u-\nabla v;B_\rho) \right )^\theta\\
  &\quad +\frac{c}{(\rho-r)^{n+2s}}|\mu|(B_1)^{1-\theta} \left (E\left(\nabla u;B_1\right)+|\mu|(B_{1}) \right )^\theta
	\end{align*}}\\
 \noindent
for some constant $c=c(n,s_0,\Lambda)$. By applying Young's inequality to the fifth line of the above inequality along with the fact that
\begin{equation*}
        E_{\mathrm{loc}}(\nabla u-\nabla v;B_\rho)\leq c\norm{\nabla u- \nabla v}_{L^1(B_{\rho})}
\end{equation*}
for some constant $c=c(n)$, we have 
\begin{align*}
  \norm{\nabla u- \nabla v}_{L^1(B_{r})}&\leq \frac{1}{2} \norm{\nabla u- \nabla v}_{L^1(B_{\rho})}+\frac{c}{(\rho-r)^{M}}|\mu|(B_{1})\\
  &\quad+\frac{c}{(\rho-r)^{M}}E\left(\nabla u;B_1\right)^{\theta}|\mu|(B_{1})^{1-\theta},
\end{align*}
where $M=M(n,s_0,\Lambda)$ and $c=c(n,s_0,\Lambda)$.
By using Lemma \ref{iterlem}, we obtain the desired result. 
\end{proof}

\begin{lemma}[Excess decay] \label{lem:OscDec} 
	Let $\mu \in L^\infty(B_{R}(x_{0}))$. Moreover, denote by $\alpha_1=\alpha_1(n,s_0,\Lambda) \in(0,1)$ the small exponent given in Lemma \ref{el : lem.holosc}. Then for any $\rho \in (0,1]$ and any weak solution $u \in W^{s,2}(B_{R}(x_0)) \cap W^{1,1}(\mathbb{R}^n)$ of $\mathcal{L} u = \mu$ in $B_{R}(x_0)$, we have
	\begin{equation*} 
		\begin{aligned}
			E(\nabla u;B_{\rho R}(x_{0})) & \leq c \rho^{\alpha_1} E(\nabla u;B_{R}(x_{0}))\\
			& \quad + c \rho^{-n} \left(\frac{|\mu|(B_R(x_0))}{R^{n-2s+1}}\right)^{1-\theta} E\left(\nabla u;B_{R}(x_{0})\right)^\theta \\ & \quad + c \rho^{-n} \frac{|\mu|(B_R(x_0))}{R^{n-2s+1}}
		\end{aligned}
	\end{equation*}
	for some constant $c=c(n,s_0,\Lambda)$,  where the constant $s_0$ is determined in \eqref{el : choi.s0.sec4}. 
\end{lemma}

\begin{proof}
	If $\rho \geq 2^{-6}$, then by an elementary computation similar to \cite[Lemma 2.4]{KuuMinSir15}, we have
	\begin{equation} \label{el : basineq.osc}
		E(\nabla u;B_{\rho R}(x_0)) \leq c \rho^{\alpha_0} E(\nabla u;B_{R}(x_{0})),
	\end{equation}
	where $c$ does not depend on $\rho$.
	
	Next, assume that $\rho \in (0,2^{-6})$. Then there is a natural number $N_\rho$ such that $2^{-6}\leq 2^{N_{\rho}}\rho<2^{-5}$. 
	Consider the weak solution $v \in W^{s,2}(B_{R/4}(x_0)) \cap L^{1}_{2s}(\mathbb{R}^n)$ of
	$$
	\begin{cases} \normalfont
		\mathcal{L} v = 0 & \text{ in } B_{R/4}(x_0) \\
		v = u & \text{ a.e. in } \mathbb{R}^n \setminus B_{R/4}(x_0).
	\end{cases}
	$$
	In view of Corollary \ref{cor:gradcomp} along with \eqref{el : basineq.osc}, we have
	\begin{align*}
		E_{\loc}(\nabla u-\nabla v;B_{\rho R}(x_0)) & \leq c \rho^{-n} R^{-n} \int_{B_{\rho R}(x_{0})} |\nabla u-\nabla v| dx \\
		& \leq c \rho^{-n} \left(\frac{|\mu|(B_R(x_0))}{R^{n-2s+1}}\right)^{1-\theta} E\left(\nabla u;B_{R}(x_{0})\right)^\theta \\ & \quad + c \rho^{-n} \frac{|\mu|(B_R(x_0))}{R^{n-2s+1}}.
	\end{align*}
	
	On the other hand, as in Lemma \ref{el : lem.tail} with $g$ replaced by $\nabla u$,  splitting into annuli along with \eqref{el : basineq.osc} yields
\begin{align*}
		&\mathrm{Tail}(\nabla u-(\nabla u)_{B_{\rho R}(x_{0})};B_{\rho R}(x_0)) \\
  & \leq c\sum_{i=0}^{N_{\rho}}2^{-2si}\dashint_{B_{2^{i}\rho R}(x_0)} |\nabla u- (\nabla u)_{B_{2^{i}\rho R}(x_0)}|\, dx + c2^{-2sN_{\rho}}E(\nabla u;B_{2_{i}N_{\rho}\rho R}(x_{0})) \\
		& \leq c\sum_{i=0}^{N_{\rho}}2^{-2si}\dashint_{B_{2^{i}\rho R}(x_0)} |\nabla u- \nabla v|\,dx\\
		&\quad + c\sum_{i=0}^{N_{\rho}}2^{-2si}\dashint_{B_{2^{i}\rho R}(x_0)} |\nabla v- (\nabla v)_{B_{2^{i}\rho R}(x_0)}| \,dx +c2^{-2sN_{\rho}}E(\nabla u;B_{R}(x_{0}))\\
  & \eqqcolon J_{1}+J_{2}+J_{3}.
	\end{align*}
 We now estimate each term $J_{1},J_{2}$ and $J_{3}$.

 \textbf{Estimate of $J_{1}$.} By Corollary \ref{cor:gradcomp} together with \eqref{el : basineq.osc}, we have
 \begin{align*}
     J_{1}&\leq c\sum_{i=0}^{N_{\rho}}\left[2^{-2si} (2^{i}\rho)^{-n} \left(\frac{|\mu|(B_R(x_0))}{R^{n-2s+1}}\right)^{1-\theta} E\left(\nabla u;B_{R}(x_{0})\right)^\theta+ (2^{i}\rho)^{-n} \frac{|\mu|(B_R(x_0))}{R^{n-2s+1}}\right]\\
     &\leq c\rho^{-n} \left(\frac{|\mu|(B_R(x_0))}{R^{n-2s+1}}\right)^{1-\theta} E\left(\nabla u;B_{R}(x_{0})\right)^\theta+ \rho^{-n} \frac{|\mu|(B_R(x_0))}{R^{n-2s+1}}
 \end{align*}
 for some constant $c=c(n,s_0,\Lambda)$.

 \textbf{Estimate of $J_{2}$.} By Lemma \ref{el : lem.holosc}, Corollary \ref{cor:gradcomp} and \eqref{el : basineq.osc}, we get
 \begin{align*}
     J_{2}&\leq c\sum_{i=0}^{N_{\rho}}2^{-2si}\dashint_{B_{2^{i}\rho R}(x_0)} |\nabla v- (\nabla v)_{B_{2^{i}\rho R}(x_0)}| dx\\
     &\leq c\sum_{i=0}^{N_{\rho}}2^{-2si}(2^{i}\rho)^{\alpha_{1}}\left[E_{\mathrm{loc}}\left(\nabla v;B_{{R}/{8}}(x_0)\right)+E\left(\nabla u;B_{{R}/{8}}(x_0)\right)+\frac{|\mu|(B_{R}(x_{0}))}{R^{n-2s+1}}\right]\\
     &\leq c\rho^{\alpha_{1}}\left[E_{\mathrm{loc}}\left(\nabla v;B_{{R}/{8}}(x_0)\right)+E\left(\nabla u;B_{{R}/{8}}(x_0)\right)+\frac{|\mu|(B_{R}(x_{0}))}{R^{n-2s+1}}\right]\\
     &\leq c\rho^{-n} \left(\frac{|\mu|(B_R(x_0))}{R^{n-2s+1}}\right)^{1-\theta} E\left(\nabla u;B_{R}(x_{0})\right)^\theta+ c\rho^{-n} \frac{|\mu|(B_R(x_0))}{R^{n-2s+1}}\\
     &\quad+ c\rho^{\alpha_{1}}E\left(\nabla u;B_{{R}}(x_0)\right)+c\frac{|\mu|(B_{R}(x_{0}))}{R^{n-2s+1}}
 \end{align*}
 for some constant $c=c(n,s_0,\Lambda)$. Finally, by the choice of the constant $N_{\rho}$, we have 
\begin{equation*}
    J_{3}\leq c\rho^{2s}E(\nabla u;B_{R}(x_{0}))\leq c\rho^{\alpha_1}E(\nabla u;B_{R}(x_{0}))
\end{equation*}
for some constant $c=c(n,s_0,\Lambda)$.
Combining all the estimates $J_{1},J_{2}$ and $J_{3}$, we obtain the desired result.
\end{proof}
Next, we prove that the averages of $\nabla u$ on arbitrarily small scales can be uniformly controlled by the Riesz potential of the data.
\begin{proposition}
\label{el : prop,main}
Let $\mu \in L^\infty(B_{R}(x_{0}))$ and let $u\in W^{s,2}(B_{2R}(x_{0}))\cap W^{1,1}(\bbR^{n})$ be a weak solution to 
\begin{equation*}
    \mathcal{L}u=\mu\quad\text{in }B_{R}(x_{0}).
\end{equation*}
Then there is a positive integer $m=m(n,s_0,\Lambda)$ such that
\begin{align*}
    \left|(\nabla u)_{B_{2^{-mi}R}(x_{0})}\right|\leq c \widetilde E(\nabla u;B_{R}(x_{0}))+cI^{|\mu|}_{2s-1}(x_{0},R)
\end{align*}
holds for any nonnegative integer $i$, where the constant $c=c(n,s_0,\Lambda)$ is independent of $i$. Here the constant $s_0$ is determined in \eqref{el : choi.s0.sec4}.
\end{proposition}
\begin{proof}
Let $m \geq 6$ to be chosen large enough, set $\rho_k:=2^{-km}$ and define the sequence of radii $R_k:=\rho_k R$ for any nonnegative integer $k$. Applying Lemma \ref{lem:OscDec} with $\rho=\rho_k \in (0,1)$ yields
\begin{align*}
	E(\nabla u;B_{R_{k+1}}(x_0)) & \leq c 2^{-\alpha_1m} E(\nabla u;B_{R_{k}}(x_{0})) \\
	& \quad + c 2^{nm} \left(\frac{|\mu|(B_{R_k}(x_0))}{R_k^{n-2s+1}}\right)^{1-\theta} E\left(\nabla u;B_{R_k}(x_{0}) \right)^\theta \\ & \quad + c 2^{nm} \frac{|\mu|(B_{R_k}(x_0))}{R_k^{n-2s+1}} ,
\end{align*}
where the constant $\alpha_1$ is determined in Lemma \ref{lem:OscDec} and the constant $c$ depends only on $n,s_0$ and $\Lambda$. Now we choose $m=m(n,s_0,\Lambda)$ large enough such that $c 2^{-\alpha_1m} \leq {1}/{2}$, so that for any $i \geq 1$, summing over $k$ leads to
\begin{align*}
	\sum_{k=1}^{i} E(\nabla u;B_{R_{k}}(x_{0})) & \leq \frac{1}{2} \sum_{k=0}^{i-1} E(\nabla u;B_{R_{k}}(x_{0})) \\
	& \quad + c\sum_{k=0}^{i-1} \left(\frac{|\mu|(B_{R_k}(x_0))}{R_k^{n-2s+1}}\right)^{1-\theta} E\left(\nabla u;B_{R_{k}}(x_{0}) \right)^\theta \\ & \quad + c\sum_{k=0}^{i-1} \frac{|\mu|(B_{R_k}(x_0))}{R_k^{n-2s+1}} ,
\end{align*}
where $c=c(n,s_0,\Lambda)$. We now apply Young's inequality to the second term given in the right-hand side of the above inequality in order to see that
\begin{align*}
	\sum_{k=1}^{i} E(\nabla u;B_{R_{k}}(x_{0})) & \leq \frac{3}{4} \sum_{k=0}^{i-1} E(\nabla u;B_{R_{k}}(x_{0})) + c\sum_{k=0}^{i-1} \frac{|\mu|(B_{R_k}(x_0))}{R_k^{n-2s+1}} ,
\end{align*}
for some constant $c=c(n,s_0,\Lambda)$.
Reabsorbing the first term on the left-hand side yields
\begin{align*}
	\sum_{k=0}^{i} E(\nabla u;B_{R_{k}}(x_{0})) & \leq c E(\nabla u;B_{R_0}(x_0))+ c\sum_{k=0}^{i-1} \frac{|\mu|(B_{R_k}(x_0))}{R_k^{n-2s+1}},
\end{align*}
where $c=c(n,s_0,\Lambda)$. We now obtain
\begin{align*}
   \left|(\nabla u)_{B_{2^{-mi}R}(x_{0})}\right|&\leq  \sum_{k=1}^{i} E(\nabla u;B_{R_{k}}(x_{0}))+\widetilde E(\nabla u;B_{R_0}(x_0))\\
   &\leq c \widetilde E(\nabla u;B_{R_0}(x_0))+ c\sum_{k=0}^{i-1} \frac{|\mu|(B_{R_k}(x_0))}{R_k^{n-2s+1}} \\
   &\leq c \widetilde E(\nabla u;B_{R_0}(x_0))+cI^{|\mu|}_{2s-1}(x_{0},R)
\end{align*}
for some constant $c=c(n,s_0,\Lambda)$ which is independent of $i$.
\end{proof}

We now prove our main theorem concerning gradient potential estimates for SOLA to \eqref{pt : eq.main} in bounded domains, which in particular involves using the localization argument given in Lemma \ref{el : lem.loc}.

\textbf{Proof of Theorem \ref{el : thm.main2}.} Fix $s_{0}\in(1/2,s]$. Let us take sequences $\{u_j\}_{j=1}^{\infty}\subset W^{s,2}(\bbR^{n})$,  $\{g_j\}_{j=1}^{\infty}\subset C_{0}^{\infty}(\bbR^{n})$ and $\{\mu_j\}_{j=1}^{\infty}\subset C_{0}^{\infty}(\bbR^{n})$ satisfying (1)-(7) given in Definition \ref{SOLA}. By following the same lines as in the proof of Theorem \ref{HD}, we get $u_j\in W^{1,1}_{\mathrm{loc}}(\Omega)$ with the estimate \eqref{el : ineq1.hd} and the convergence result \eqref{el : conv.hd}.
Let us fix $B_{R}(x_{0})\subset \Omega$ and take a cut off function $\xi\in C_{c}^{\infty}\left(B_{{4R}/{5}}(x_{0})\right)$ with $\xi\equiv 1$ on $B_{{3R}/{5}}(x_{0})$. Then by Lemma \ref{el : lem.loc}, we have that $w_{j}\coloneqq u_{j}\xi\in W^{s,2}(B_{{4R}/{5}}(x_{0}))\cap W^{1,1}(\bbR^{n})$ is a weak solution to 
\begin{equation*}
    \mathcal{L}w_j=\mu_{j}+f_{j}\quad\text{in }B_{{R}/{5}}(x_{0}), 
\end{equation*}
where 
\begin{align*}
    f_j(x)&=2(1-s)\int_{\mathbb{R}^{n}\setminus B_{{3R}/{5}}(x_{0})}{\Phi}\left(\frac{w_j(x)-w_j(y)}{|x-y|^{s}}\right)\frac{\,dy}{|x-y|^{n+s}}\\
    &\quad-2(1-s)\int_{\mathbb{R}^{n}\setminus B_{{3R}/{5}}(x_{0})}{\Phi}\left(\frac{u_j(x)-u_j(y)}{|x-y|^{s}}\right)\frac{\,dy}{|x-y|^{n+s}}
\end{align*}
and $f_j\in L^{\infty}\left(B_{{2R}/{5}}(x_{0})\right)$ by \eqref{el : ineq.tail}.
By taking $\nu_j\coloneqq\mu_j+f_j$, we observe that $w_j\in W^{s,2}(B_{{4R}/{5}}(x_{0}))\cap W^{1,1}(\bbR^{n})$ is a weak solution to 
\begin{equation*}
    \mathcal{L}w_j=\nu_j\quad\text{in }B_{{R}/{5}}(x_{0}).
\end{equation*}
Here, we note from \eqref{el : ineq3.hd} with $x_{1}=x_{0}$ and $\rho=R$ that
\begin{equation}
\label{el : ineq.thm}
    E\left(\nabla w_j;B_{{2R}/{5}}(x_{0})\right)\leq c\widetilde{E}(\nabla u_j;B_{4R/5}(x_{0}))+c\widetilde{E}(u_j/R;B_{4R/5}(x_{0}))
\end{equation}
for some constant $c=c(n,s_0,\Lambda)$.
We next observe from \eqref{el : ineq.tail} that
\begin{equation*}
    \|f_j\|_{L^{\infty}\left(B_{{2R}/{5}}(x_{0})\right)}\leq cR^{-2s}\mathrm{Tail}\left(u_j;B_{{3R}/{5}}(x_{0})\right)
\end{equation*}
for some constant $c=c(n,\Lambda)$. Thus, we have
\begin{align*}
    I^{|f_j|}_{2s-1}\left(x_{0},{2R}/{5}\right)&\leq \int_{0}^{R}\frac{1}{t^{n-2s+1}}t^{n}\|f_j\|_{L^{\infty}\left(B_{{2R}/{5}}(x_{0})\right)}\frac{\,dt}{t}\\
    &\leq {c}\widetilde{E}\left(u_j/{R};B_{R}(x_{0})\right),
\end{align*}
which implies 
\begin{align*}
    I^{|\nu_j|}_{2s-1}\left(x_{0},{2R}/{5}\right)&\leq I^{|\mu_j|}_{2s-1}\left(x_{0},{2R}/{5}\right)+I^{|f|}_{2s-1}\left(x_{0},{2R}/{5}\right)\\
    &\leq I^{|\mu_j|}_{2s-1}\left(x_{0},{2R}/{5}\right)+\widetilde{E}\left(u_j/{R};B_{R}(x_{0})\right).
\end{align*}
Using this along with Proposition \ref{el : prop,main}, \eqref{el : ineq.thm} and \eqref{el : ineq1.hd} with $x_{1}=x_{0}$ and $\rho=R$, we obtain that for any nonnegative integer $i$
\begin{align*}
    \left|(\nabla u_j)_{B_{2^{-mi}{R}/{5}}(x_{0})}\right|&=\left|(\nabla w_j)_{B_{2^{-mi}{R}/{5}}(x_{0})}\right|\\
    &\leq \widetilde E\left(\nabla w_j;B_{{2R}/{5}}(x_{0})\right)+I^{|\mu_j|}_{2s-1}\left(x_{0},{2R}/{5}\right)\\
    &\quad+c\widetilde{E}\left(u_j/{R};B_{R}(x_{0})\right)\\
    &\leq cI^{|\mu_j|}_{2s-1}\left(x_{0},R\right)+c\widetilde{E}\left(u_j/R;B_{R}(x_{0})\right)
\end{align*}
holds, where $c=c(n,s_0,\Lambda)$ and the constant $m=m(n,s_0,\Lambda)$ is determined in Proposition \ref{el : prop,main}. Since $u_j-(u_j)_{B_{R}(x_{0})}$ satisfies
\begin{equation*}
    \mathcal{L}(u_j-(u_j)_{B_{R}(x_{0})})=\mu_j\quad\text{in }\Omega,
\end{equation*}
we obtain 
\begin{equation*}
    \left|(\nabla u_j)_{B_{2^{-mi}{R}/{5}}(x_{0})}\right|\leq cE(u_j/R;B_{R}(x_{0}))+cI^{|\mu_j|}_{2s-1}\left(x_{0},R\right)
\end{equation*}
for some constant $c=c(n,s_0,\Lambda)$. By taking the limit $j\to \infty$ along with \eqref{el : conv.hd} and (5)-(7) given in Definition \ref{SOLA}, we get
\begin{equation*}
    \left|(\nabla u)_{B_{2^{-mi}{R}/{5}}(x_{0})}\right|\leq cE(u/R;B_{R}(x_{0}))+cI^{|\mu|}_{2s-1}\left(x_{0},R\right)
\end{equation*}
for some constant $c=c(n,s_0,\Lambda)$. Since the constant $c$ is independent of $i$, we obtain the desired estimate by taking $i\to\infty$.
\qed

Finally, we are in the position to prove Theorem \ref{el : thm.main3}.

\begin{proof}[Proof of Theorem \ref{el : thm.main3}]
	Fix $s_0 \in (1/2,s]$. Since any weak distributional solution $u \in W^{s,2}(\mathbb{R}^n)$ of \eqref{pt : eq.main} clearly is a SOLA to \eqref{NonlocalDir} for any bounded domain $\Omega \subset \mathbb{R}^n$ with respect to $g=u$, by Theorem \ref{el : thm.main2}, for any fixed $x_0 \in \mathbb{R}^n$ and any $R>0$ we obtain the estimate
	\begin{equation} \label{Deq}
		|\nabla u(x_0)| \leq c E(u/R;B_{R}(x_{0}))+cI^{|\mu|}_{2s-1}\left(x_{0},R\right),
	\end{equation}
	where $c=c(n,s_0,\Lambda)$.
	Now observe that for $c=c(n,s_0)$, by H\"older's inequality
	\begin{align*}
		E(u/R;B_{R}(x_{0})) & \leq c R^{-1} \left (E_{\loc}^2 (u;B_{R}(x_{0})) + R^{s} \left ( \int_{\mathbb{R}^n \setminus B_{R}(x_{0})} \frac{|u(y)|^2}{|x_0-y|^{n+2s}} \right )^\frac{1}{2} \right ) \\
		& \leq c R^{-\frac{n}{2}-1} ||u||_{L^2(\mathbb{R}^n)} \to 0 \quad \text{as } R \to \infty.
	\end{align*}
	Therefore, taking into account \eqref{truncated} and letting $R \to \infty$ in \eqref{Deq} yields the estimate \eqref{PERn} for some $c=c(n,s_0,\Lambda)$. The proof is finished.
\end{proof}

\printbibliography

\end{document}